\documentclass{mcom-l}
\usepackage{xcolor}
\colorlet{siaminlinkcolor}{green!50!black}
\colorlet{siamexlinkcolor}{red!100!black}
\usepackage{amssymb}
\usepackage{bbm}
\usepackage{amsmath,mathtools}
\usepackage{comment}
\usepackage{empheq}
\usepackage{mathrsfs} 

\usepackage{graphicx,float}

\makeatletter
\@namedef{subjclassname@2020}{\textup{2020} Mathematics Subject Classification}
\makeatother

\makeatletter
\newcommand{\thickhline}{%
	\noalign {\ifnum 0=`}\fi \hrule height 1.2pt
	\futurelet \reserved@a \@xhline
}
\makeatother

\newcommand{\vertiii}[1]{{\left\vert\kern-0.25ex\left\vert\kern-0.25ex\left\vert #1 
    \right\vert\kern-0.25ex\right\vert\kern-0.25ex\right\vert}}

\newcommand{\lvertiii}[1]{{\left\vert\kern-0.25ex\left\vert\kern-0.25ex\left\vert #1 
    \right.\kern-0.25ex\right.\kern-0.25ex\right.}}

\newcommand{\rvertiii}[1]{{\left.\kern-0.25ex\left.\kern-0.25ex\left. #1 
    \right\vert\kern-0.25ex\right\vert\kern-0.25ex\right\vert}}

\DeclareFontFamily{U}{mathx}{}
\DeclareFontShape{U}{mathx}{m}{n}{<-> mathx10}{}
\DeclareSymbolFont{mathx}{U}{mathx}{m}{n}
\DeclareMathAccent{\widecheck}{0}{mathx}{"71}

\usepackage{amsopn}

\usepackage[hidelinks]{hyperref}
\usepackage[capitalize,nameinlink,noabbrev]{cleveref}

\newtheorem{theorem}{Theorem}[section]
\newtheorem{lemma}[theorem]{Lemma}
\newtheorem{proposition}[theorem]{Proposition}

\theoremstyle{definition}

\theoremstyle{remark}
\newtheorem{remark}[theorem]{Remark}

\numberwithin{equation}{section}
\numberwithin{figure}{section}

\crefname{example}{Example}{Examples}
\crefname{hypothesis}{Hypothesis}{Hypotheses}
\crefname{conj}{Conjecture}{Conjectures}
\crefformat{equation}{\textup{#2(#1)#3}}
\crefrangeformat{equation}{\textup{#3(#1)#4--#5(#2)#6}}
\crefmultiformat{equation}{\textup{#2(#1)#3}}{ and \textup{#2(#1)#3}}
{, \textup{#2(#1)#3}}{, and \textup{#2(#1)#3}}
\crefrangemultiformat{equation}{\textup{#3(#1)#4--#5(#2)#6}}%
{ and \textup{#3(#1)#4--#5(#2)#6}}{, \textup{#3(#1)#4--#5(#2)#6}}{, and \textup{#3(#1)#4--#5(#2)#6}}

\Crefformat{equation}{#2Equation~\textup{(#1)}#3}
\Crefrangeformat{equation}{Equations~\textup{#3(#1)#4--#5(#2)#6}}
\Crefmultiformat{equation}{Equations~\textup{#2(#1)#3}}{ and \textup{#2(#1)#3}}
{, \textup{#2(#1)#3}}{, and \textup{#2(#1)#3}}
\Crefrangemultiformat{equation}{Equations~\textup{#3(#1)#4--#5(#2)#6}}%
{ and \textup{#3(#1)#4--#5(#2)#6}}{, \textup{#3(#1)#4--#5(#2)#6}}{, and \textup{#3(#1)#4--#5(#2)#6}}

\hypersetup{
	allcolors = siaminlinkcolor,
	linkcolor = siamexlinkcolor,
}
\usepackage{cite}
\begin{document}

\title[Convergence of ESFEM for the EKS model of tumour growth]{Convergence of finite elements for the Eyles-King-Styles model of tumour growth}
%\title{}


%    author one information
%\author[short version for running head]{name for top of paper}
%\author[W. Bao]{Weizhu Bao}
%\address{Department of Mathematics, National University of Singapore, Singapore 119076}
%\curraddr{}
%\email{matbaowz@nus.edu.sg}
%\thanks{This work was partially supported by the Ministry of Education of Singapore under its AcRF Tier 2 funding MOE-T2EP20122-0002 (A-8000962-00-00). Part of the work was done when the authors were visiting the Institute of Mathematical Science at the National University of Singapore in 2023.}

%    author two information
\author[Y. Li]{Yifei Li}
\address{Mathematisches Institut, Universit\"at T\"{u}bingen, Auf der Morgenstelle 10., 72076, T\"{u}bingen, Germany}
%\curraddr{}
\email{yifei.li@mnf.uni-tuebingen.de}
%\thanks{}

%    \subjclass is required.
\date{}
%\dedicatory{}

%    Abstract is required.
\begin{abstract}
    This paper presents a convergence analysis of the evolving surface finite element method (ESFEM) applied to the original Eyles-King-Styles model of tumour growth. The model consists of a Poisson equation in the bulk, a forced mean curvature flow on the surface, and a coupled velocity law between bulk and surface. Due to the non-trivial bulk-surface coupling, all previous analyses required an additional regularization term. By introducing a $H^{1/2}(\Gamma)$ energy estimates theory, we develop an essentially new theoretical framework that addresses the intrinsic bulk-surface coupling. Based on this framework, we provide the first rigorous convergence proof for the original model without regularization.
\end{abstract}

\maketitle


\section{Introduction}
In this paper, we present a rigorous convergence  analysis for the evolving surface finite element method (ESFEM) applied to the original Eyles-King-Styles model of tumour growth \cite{EKS19}. In the Eyles-King-Styles model, the tumour growth is described by both the bulk domain $\Omega(t) \subset \mathbb{R}^3$ and its free boundary $\Gamma(t)$ evolving under the following coupled bulk--surface system:
\begin{subequations}
   \label{eq:EKS19}
   \begin{alignat}{3}
       \label{eq:EKS19 poisson}
       &-\Delta u = -1 &&\quad\text{in } &&\Omega(t), \\
       \label{eq:EKS19 robin}
       &\partial_{n} u + \alpha u = \beta H + Q &&\quad\text{on } &&\Gamma(t), \\
       \label{eq:EKS19 mcf}
       &v_{\Gamma} = Vn, \qquad \text{with } V = - \beta H + \alpha u &&\quad\text{on } &&\Gamma(t).
   \end{alignat}
\end{subequations}
Here $u(x, t)$ denotes the tissue pressure, $n(x, t)$ is the outward unit normal vector on $\Gamma(t)$, $H(x, t)$ is the mean curvature of $\Gamma(t)$, $V(x, t)$, $v_{\Gamma}(x, t)$ are the normal velocity and the velocity field on $\Gamma(t)$, $Q(x, t)$ is a given surface source term, and $\alpha>0, \beta>0$ are given constants. This coupled bulk--surface system consists of a Poisson equation \eqref{eq:EKS19 poisson} for $u$ with a Robin boundary condition \eqref{eq:EKS19 robin}, which in turn affects the evolution of the surface $\Gamma(t)$ in the forced mean curvature flow \eqref{eq:EKS19 mcf}. This non-trivial coupling of the velocity field $v_{\Gamma}$ and the tissue pressure $u$ makes the numerical analysis of this model particularly challenging.

The evolving surface finite element method is an important numerical approach for solving partial differential equations on evolving domains \cite{Dziuk88,DziukElliott_ESFEM,DziukElliott_SFEM,Elliott2010,ElliottRanner_unified}, especially the bulk-surface coupled system \cite{ElliottRanner_bulksurface,EKL24}. In fact, Edelmann, Kov\'{a}cs and Lubich have developed and analyzed an evolving bulk-surface finite element method (which is also denoted as ESFEM) for the Eyles-King-Styles model and proved its convergence \cite{EKL24}. However, their analysis was achieved by replacing \eqref{eq:EKS19 robin} with a generalized Robin boundary condition with an additional surface Laplacian term $\Delta_\Gamma u$. Although this modification was necessary in their convergence analysis in controlling the boundary error in $H^1(\Gamma)$ norm, the numerical experiments suggested the expected convergence rates can be achieved for the original tumour model without the modification. This gap between the theoretical analysis and numerical observations motivates us to develop a rigorous convergence analysis for the original tumour model without regularization. Indeed, providing such a rigorous convergence analysis for the original bulk-surface coupled model is an open problem.

The fundamental challenge in analyzing the original tumour model arises from its intrinsic bulk-surface coupling structure. In equation \eqref{eq:EKS19 mcf}, the $H^1(\Gamma)$ estimates for $V$ requires $H^1(\Gamma)$ boundary estimates for $u$, which cannot be obtained from the PDE. Previous ESFEM analyzes avoided this difficulty through the addition of regularization terms. The first approach replaces $V$ with $-\Delta_\Gamma V + V$ in \eqref{eq:EKS19 mcf}, so that $H^1(\Gamma)$ estimates for $V$ can be obtained using only $L^2(\Gamma)$ estimates for $u$, resulting in a $\left(H^1(\Omega); L^2(\Gamma)\right)$ estimate for $u$ and its boundary trace $\gamma(u):= u|_{\Gamma}$ \cite{KovacsPower_quasilinear,KLLP2017,ElliottRanner_bulksurface}. The second approach replaces $u$ with $-\Delta_\Gamma u + u$ in \eqref{eq:EKS19 robin}, enabling a $H^1(\Gamma)$ estimates for $\gamma(u)$ by energy estimates, thus resulting in a $\left(H^1(\Omega); H^1(\Gamma)\right)$ estimate for $u$ and $\gamma(u)$ \cite{EKL24}. However, the trace inequality $\left\|\gamma(u)\right\|_{H^{1/2}(\Gamma)} \leq C \left\|u\right\|_{H^1(\Omega)}$ indicates that intrinsic bulk-surface coupling should have a mixed regularity estimate $\left(H^1(\Omega); H^{1/2}(\Gamma)\right)$ for $u$ and $\gamma(u)$. Therefore, neither regularization can capture the intrinsic coupling. Furthermore, although various ESFEMs have been proposed for different equations, their estimates still fall into these two cases, and hence are also insufficient for intrinsic bulk-surface coupling. To address the intrinsic bulk-surface coupling, we develop a new analytical framework, which can establish a $\left(H^1(\Omega); H^{1/2}(\Gamma)\right)$ estimate.

To this end, we introduce a comprehensive framework consisting of four components: using $H^1(\Omega)$ norm to control $H^{1/2}(\Gamma)$ norm; using $H^{1/2}(\Gamma)$ norm to control $H^1(\Omega)$ norm; $H^1(\Omega)$ energy estimates in the bulk; and $H^{1/2}(\Gamma)$ energy estimates on the surface. The first relies on the trace inequality. The second has been established in Proposition 6 of \cite{EKL24}, and we are the first to apply it to intrinsic coupling problems. The third has also been well-established in ESFEM analysis, which is based on $L^2$-$L^2$-$L^\infty$ estimates and a derivative trick. The challenge appears in the last, which requires developing corresponding $H^{1/2}$-$H^{-1/2}$-$W^{1/2,\infty}$ estimates and a fractional order derivative trick. This creates two new difficulties: for a finite element function $u_h$, the $W^{1/2,\infty}$ norm of its gradient $\nabla_\Gamma u_h$ is not well-defined; and fractional order derivative trick requires estimating the $L^2$ norm of $\frac{d}{dt}(I - \Delta_\Gamma)^{1/2}$. These difficulties are fundamentally caused by the nature of fractional order and cannot be handled in existing ESFEM analysis.

We overcome these challenges by developing a $H^{1/2}(\Gamma)$ energy estimates theory. As a result, we present the first complete error analysis for the ESFEM applied to the original tumour model without regularization, addressing the intrinsic bulk-surface coupling that previous ESFEM works avoided. Our numerical method is based on the ESFEM for the regularized tumour model in \cite{EKL24}. Our error analysis is composed of stability analysis and consistency analysis, the stability analysis is derived by the framework, and the consistency analysis follows the results in \cite{EKL24}. The main contributions of this paper include:
\begin{itemize}
    \item \textbf{Novel methods for handling fractional order operators.} Beyond inverse estimates and interpolation, we introduce a new approach using Sylvester equation to handle fractional order operators \eqref{eq:explicit representation of solution to Sylvester equation}.
    \item \textbf{A comprehensive $H^{1/2}(\Gamma)$ energy estimates framework.} We establish the corresponding $H^{1/2}$-$H^{-1/2}$-$W^{1/2,\infty}$ estimates (Theorem \ref{theorem:estimate of the multilinear form}) and develop fractional order derivative tricks (Lemma \ref{lemma:time derivative of discrete Laplacian of order 1} --- \ref{lemma:time derivative of discrete Laplacian of order -1/2}), enabling rigorous stability analysis for intrinsic bulk-surface coupling problems.
    \item \textbf{First rigorous convergence analysis for the original tumour growth model.} We apply the framework to prove optimal convergence rates in the mixed regularity space $\left(H^1(\Omega), H^{1/2}(\Gamma)\right)$ for finite element spaces of polynomial degree $k \geq 3$.
 \end{itemize}

This paper is organized as follows: Section 2 reviews the evolving surface finite element method (ESFEM) for the tumour model. Including the continuous formulation and variational formulation of the model, and its spatial semi-discretization by ESFEM in the matrix-vector form. Section 3 establishes the $H^{1/2}(\Gamma)$ energy estimates framework, including the basic properties of the $H^{1/2}(\Gamma)$ spaces, the $H^{1/2}$-$H^{-1/2}$-$W^{1/2,\infty}$ estimates and the fractional order derivative trick. Finally, Section 4 presents a comprehensive error analysis for the ESFEM applied to the original tumour model without regularization, with particular emphasis on applying the $H^{1/2}(\Gamma)$ energy estimates to stability analysis.

\section{The convergent evolving bulk-surface finite element method}

\subsection{Basic notions and notations}
Suppose the initial domain $\Omega^0$ is a bounded domain in $\mathbb{R}^3$ with a smooth boundary $\Gamma^0$. Consider a smooth flow map $X(\cdot, t): \Omega^0  \to \mathbb{R}^3$ defined on $\Omega^0$, where $X(\cdot, t)$ constitutes an embedding for each $t \in [0, T]$. The time-dependent domain $\Omega(t)$ and its boundary $\Gamma(t)$ are parameterized by $X$ as follows:

\begin{equation*}
\Omega(t) = \left\{ X(q,t) \ | \ q \in \Omega^0 \right\}, \qquad \Gamma(t) = \left\{ X(q,t) \ | \ q \in \Gamma^0 \right\}.
\end{equation*}
In particular, the flow map $X(\cdot, 0)$ is the identity map on $\Omega^0$, i.e., $X(q, 0) = q$ for all $q \in \Omega^0$. Therefore, we have $\Omega(0) = \Omega^0$ and $\Gamma(0) = \Gamma^0$.

We further assume that $X(\cdot, t)$ is a homeomorphism from $\Omega^0$ to $\Omega(t)$, i.e., $X(\cdot, t)$ is a bijective map from $\Omega^0$ to $\Omega(t)$. Therefore, for any function $f(x, t)$ defined on $\Omega(t)$ or $\Gamma(t)$, we have a corresponding function $f(X(q, t), t)$ defined on $\Omega^0$ or $\Gamma^0$ via the pullback by $X(\cdot, t)$. Moreover, we introduce its material derivative \cite{barrett2020parametric}:
\begin{equation}\label{eq:material derivative}
    \partial^{\bullet}{f}(x,t) := \frac{d}{dt} f(X(q,t),t) = \nabla f(x,t) \cdot v(x,t) + \partial_t f(x,t), \quad \text{where } x = X(q,t).
\end{equation}

Specifically, we consider the material derivative of the identity map $\text{id}(x, t)$ on $\Omega(t)$, which gives the velocity field on $\Omega(t)$ as:
\begin{equation}\label{eq:velocity field}
    v(X(q, t), t) = \partial^{\bullet}{\text{id}}(x,t) = \frac{d}{dt} \text{id}(X(q,t),t) = \partial_t X(q, t), \quad \text{for } q \in \Omega^0.
\end{equation}
On the other hand, for a given velocity field $v$ on $\Omega(t)$, we can define the corresponding flow map $X(\cdot, t)$ by solving the initial value problem. The same relationship applies to the boundaries, and the surface flow map and surface velocity field are given by
\begin{equation*}
    X_{\Gamma}(\cdot, t) = X(\cdot, t)|_{\Gamma^0}, \quad \text{and} \quad v_{\Gamma}(\cdot, t) = v(\cdot, t)|_{\Gamma(t)}.
\end{equation*}
For the tumour model \eqref{eq:EKS19}, only the surface velocity field $v_{\Gamma}$ is given. The bulk velocity field $v$ is obtained by harmonic extension of the surface velocity field $v_{\Gamma}$:
\begin{equation*}
    - \Delta v = 0 \quad \text{in } \Omega(t), \quad v = v_{\Gamma} \quad \text{on } \Gamma(t).
\end{equation*}

Denoting the outward unit normal vector of $\Gamma(t)$ as $n$, we employ the following surface differential operators \cite{DeckelnickDE2005}: 
\begin{alignat*}{3} % Using alignat* ensures no line numbers
    &\nabla_{\Gamma} f &&= \nabla f - (\nabla f \cdot n) n, \quad &&\text{surface gradient for } f, \\
    &\nabla_{\Gamma} \boldsymbol{f} &&= (\nabla_{\Gamma} f_1, \nabla_{\Gamma} f_2, \nabla_{\Gamma} f_3)^T, \quad &&\text{surface Jacobian for } \boldsymbol{f} = (f_1, f_2, f_3)^T, \\
    &\nabla_{\Gamma} \cdot \boldsymbol{f} &&= \text{tr}(\nabla_{\Gamma} \boldsymbol{f}), \quad &&\text{surface divergence for } \boldsymbol{f}, \\
    & \Delta_{\Gamma} f &&= \nabla_{\Gamma} \cdot (\nabla_{\Gamma} f), \quad &&\text{surface Laplace-Beltrami for } f.
\end{alignat*}
We refer to \cite{Walker2015} for more details on the surface differential operators.

The Weingarten map $A = \nabla_{\Gamma} n$ (or the second fundamental form), being a symmetric $3 \times 3$ matrix, has $n$ as an eigenvector corresponding to the eigenvalue $0$ \cite{Walker2015}. The remaining two eigenvalues, denoted by $\kappa_1$ and $\kappa_2$, are termed the principal curvatures. We thus introduce the important geometric quantities:
\begin{equation*}%\label{eq:mean curvature and weingarten map}
    H = \text{tr}(A) = \kappa_1 + \kappa_2, \qquad |A|^2 = \kappa_1^2 + \kappa_2^2.
\end{equation*}
Here $H$ is the mean curvature of the surface $\Gamma(t)$, and $|A|^2$ is the squared Frobenius norm of the Weingarten map. From \cite[Lemma 16]{barrett2020parametric}, the surface gradient of $H$ is
\begin{equation*}%\label{eq:surface gradient of mean curvature}
    \nabla_{\Gamma} H = - \left( \Delta_{\Gamma} + |A|^2 \right) n.
\end{equation*}
The material derivative of $n$ and $H$ satisfies the following Huisken's identity \cite{Huisken1984}
\begin{equation*}%\label{eq:huisken identity}
    \partial^{\bullet} H = (\Delta_{\Gamma} + |A|^2) V, \qquad \partial^{\bullet} n = -\nabla_{\Gamma} V.
\end{equation*}

Using the above notations, the tumour growth model \eqref{eq:EKS19} can be written as
\begin{subequations}\label{eq:EKS19 alt form}
\begin{alignat}{3}
    &\left\{
    \begin{array}{l@{\hspace{2em}}l}
    -\Delta u = -1 & \qquad\qquad\qquad\qquad\qquad\text{in }\Omega(t), \\
    \partial_{n} u + \alpha u = \beta H + Q & \qquad\qquad\qquad\qquad\qquad \text{on }\Gamma(t);
    \end{array}
    \right.
    \label{eq:robin-bvp}
    \\[3mm]
    &\left\{
    \begin{array}{l@{\hspace{2em}}l}
    \partial^{\bullet} n = \beta \Delta_\Gamma n + \beta |A|^2 n - \alpha \nabla_\Gamma u & \text{on }\Gamma(t), \\[1mm]
    \partial^{\bullet} H = \beta \Delta_\Gamma H + \beta |A|^2 H - \alpha \Delta_\Gamma u - \alpha |A|^2 u & \text{on }\Gamma(t), \\[1mm]
    v_\Gamma = Vn \quad\text{with }V=- \beta H + \alpha u & \text{on }\Gamma(t);
    \end{array}
    \right.
    \label{eq:forced-mcf}
    \\[3mm]
    &\left\{
    \begin{array}{l@{\hspace{2em}}l}
    -\Delta v = 0 & \qquad\qquad\qquad\qquad\qquad\qquad\qquad\quad \text{in }\Omega(t), \\
    v = v_\Gamma & \qquad\qquad\qquad\qquad\qquad\qquad\qquad\quad \text{on }\Gamma(t);
    \end{array}
    \right.
    \label{eq:harmonic extension}
    \\[3mm]
    &\begin{array}{l@{\hspace{2em}}l}
    \partial_t X = v\circ X & \qquad\qquad\qquad\qquad\qquad\qquad\quad \quad\,\,\text{on }\Omega^0 \cup \Gamma^0.
    \end{array}
    \label{eq:ODE for flow map}
    \end{alignat}
\end{subequations}

Here \eqref{eq:robin-bvp} is the Robin boundary value problem, \eqref{eq:forced-mcf} is the forced mean curvature flow, \eqref{eq:harmonic extension} is the harmonic extension of the surface velocity field $v_{\Gamma}$, and \eqref{eq:ODE for flow map} is the ODE for obtaining the flow map $X(\cdot, t), X_{\Gamma}(\cdot, t)$.

\subsection{The weak formulation}
For notational convenience, we omit the time $t$ when the context is clear. We introduce the trace operator $\gamma$, such that $\gamma(f) = f|_{\Gamma}$.  Using this, the weak formulation of the governing equations \eqref{eq:EKS19 alt form} can thus be written as
\begin{subequations}\label{eq:weak form}
\begin{alignat}{3}
    &\begin{array}{ll}
    \int_{\Omega} \nabla u \cdot \nabla \phi^u + \alpha \int_{\Gamma} \gamma(u) \gamma(\phi^u) = -\int_{\Omega} \phi^u + \int_{\Gamma} (\beta H + Q)\gamma(\phi^u) &  \forall \phi^u \in H^1(\Omega);
    \end{array}
    \label{eq:robin-bvp, weak form}
    \\[3mm]
    &\left\{
    \begin{array}{ll}
    \int_{\Gamma} \partial^{\bullet} n \cdot \phi^{n} + \beta \int_{\Gamma} \nabla_{\Gamma} n \cdot \nabla_{\Gamma} \phi^{n} = \beta \int_{\Gamma} |A|^2 n \cdot \phi^{n} - \alpha \int_{\Gamma} \nabla_{\Gamma} \gamma(u) \cdot \phi^{n}, \\[1mm]
    \int_{\Gamma} \partial^\bullet H \phi^H + \beta \int_{\Gamma} \nabla_\Gamma H \cdot \nabla_\Gamma \phi^H = - \beta \int_\Gamma |A|^2 V \phi^H + \alpha \int_{\Gamma} \nabla_\Gamma \gamma(u) \cdot \nabla_\Gamma \phi^H, \\[1mm]
    v_\Gamma = Vn, \quad\quad\quad\text{with }V=- \beta H + \alpha u, \quad\quad\quad\quad\,\, \forall \phi^{n} \in [H^1(\Gamma)]^3, \phi^H \in H^1(\Gamma);
    \end{array}\right.
    \label{eq:forced-mcf, weak form}
    \\[3mm]
    &\begin{array}{ll}
    \int_{\Omega} \nabla v \cdot \nabla \phi^v = 0, \qquad \text{with } \gamma(v) = v_\Gamma, & \qquad\qquad\quad  \forall \phi^v \in [H_0^1(\Omega)]^3;
    \end{array}
    \label{eq:harmonic extension, weak form}
    \\[3mm]
    &\begin{array}{ll}
    \partial_t X = v\circ X, & \qquad\qquad\qquad\qquad\qquad\qquad\qquad\qquad \text{on }\Omega^0 \cup \Gamma^0.
    \end{array}
    \label{eq:ODE for flow map, weak form}
\end{alignat}
\end{subequations}

\subsection{Evolving surface finite element discretization}
\subsubsection{Evolving surface finite elements}
We discretize the initial bulk-surface domain $\Omega^0 \cup \Gamma^0$ using a tetrahedral mesh $\Omega_h^0 \cup \Gamma_h^0$ of degree $k$ with maximal element diameter $h$. This mesh is quasi-uniform and satisfies the shape regularity as detailed in \cite{DziukElliott_ESFEM,ElliottRanner_bulksurface}.

The nodes of the triangulation, denoted by $q_j$, are organized as two separate sets:
\begin{equation*}
    \mathbf{x}^0 = \begin{pmatrix} \mathbf{x}_\Gamma^0 \\ \mathbf{x}_\Omega^0 \end{pmatrix} \in \bigl(\mathbb{R}^3\bigr)^N,
\end{equation*}
where $\mathbf{x}_\Gamma^0 = (q_j)_{j=1}^{N_\Gamma}$ are the boundary nodes on $\Gamma^0$, and $\mathbf{x}_\Omega^0 = (q_j)_{j=N_\Gamma+1}^{N}$ contains the remaining $N_\Omega = N - N_\Gamma$ interior nodes.

We define the Lagrangian finite element basis functions:
\begin{align*}
    \varphi_i[\mathbf{x}^0] &: \Omega_h^0[\mathbf{x}^0] \to \mathbb{R}, \quad i = 1, \ldots, N, \\
    \psi_i[\mathbf{x}^0] &: \Gamma_h^0[\mathbf{x}^0] \to \mathbb{R}, \quad i = 1, \ldots, N_\Gamma,
\end{align*}
which are polynomials of degree $k$ on each element and satisfy the nodal properties:
\begin{align*}
    \varphi_i[\mathbf{x}^0](q_j) &= \delta_{ij}, \quad 1 \leq i,j \leq N, \\
    \psi_i[\mathbf{x}^0](q_j) &= \delta_{ij}, \quad 1 \leq i, j \leq N_\Gamma.
\end{align*}
Moreover, since $\Gamma_h^0 = \partial \Omega_h^0$, these basis functions satisfy the trace property:
\begin{equation} \label{eq:trace property}
    \psi_j[\mathbf{x}^0] = \gamma_h(\varphi_j[\mathbf{x}^0]), \quad 1 \leq j \leq N_\Gamma,
\end{equation}
where $\gamma_h$ is the trace operator on $\Gamma_h[\mathbf{x}^0]$.

The finite element spaces are then defined as:
\begin{subequations}\label{eq:finite element spaces}
\begin{align}
    \mathscr{V}_h[\mathbf{x}^0] &= \mathrm{span} \{ \varphi_1[\mathbf{x}^0],\ldots,\varphi_N[\mathbf{x}^0] \}, \\
    \mathscr{V}_h^0[\mathbf{x}^0] &= \mathrm{span} \{ \varphi_{N_\Gamma+1}[\mathbf{x}^0],\ldots,\varphi_N[\mathbf{x}^0] \} = \{ u_h \in \mathscr{V}_h[\mathbf{x}^0] : \gamma_h(u_h) = 0 \}, \\
    \mathscr{S}_h[\mathbf{x}^0] &= \mathrm{span} \{ \psi_1[\mathbf{x}^0],\ldots,\psi_{N_\Gamma}[\mathbf{x}^0] \}.
\end{align}
\end{subequations}

Suppose the temporal evolution of these nodes $\mathbf{x}^0$ is given by:
\begin{equation*}
    \mathbf{x}(t) = \begin{pmatrix} \mathbf{x}_\Gamma(t) \\ \mathbf{x}_\Omega(t) \end{pmatrix}  \in \bigl(\mathbb{R}^3\bigr)^N, \quad \text{with } \mathbf{x}(0) = \mathbf{x}^0.
\end{equation*}
This characterizes the discrete flow map $X_h(\cdot,t) : \Omega_h^0\cup \Gamma_h^0 \to \mathbb{R}^3$ as:
\begin{equation}
    X_h(q_h,t) = \sum_{j=1}^N \mathbf{x}_j(t) \varphi_j[\mathbf{x}^0](q_h), \quad q_h \in \Omega_h^0\cup \Gamma_h^0,
\end{equation}
where $X_h(q_j,t) = \mathbf{x}_j(t)$ at the nodes $q_j$ for $j=1,\dots,N$. And $X_h(q_h,0) = q_h$ for all $q_h \in \Omega_h^0$, i.e. $X_h(\cdot, 0)$ is an identity map on $\Omega_h^0\cup \Gamma_h^0$. The discrete evolving domain $\Omega_h(t)$ and its boundary $\Gamma_h(t)$ are thus determined:
\begin{equation*}
    \Omega_h(t) = \left\{ p_h\in \mathbb{R}^3 : p_h = \sum_{i=1}^N X_h(q_i,t) \varphi_i[\mathbf{x}^0](q_h), q_h \in \Omega_h^0 \right\}, \quad \Gamma_h(t) = \partial \Omega_h(t).
\end{equation*}

Moreover, we assume that $X_h(\cdot, t)$ is a homeomorphism from $\Omega_h^0 \cup \Gamma_h^0$ to $\Omega_h(t) \cup \Gamma_h(t)$. Therefore, it has a regular inverse $X_h^{-1}(\cdot, t) : \Omega_h(t) \cup \Gamma_h(t) \to \Omega_h^0 \cup \Gamma_h^0$. The Lagrangian basis functions on $\Omega_h(t)$ and $\Gamma_h(t)$ can be defined through the pullback using $X_h^{-1}(\cdot, t)$ as follows:
\begin{subequations}\label{eq:pullback basis functions}
\begin{align}
\varphi_i[\mathbf{x}_j(t)](q_h) &= \varphi_i[\mathbf{x}^0](X_h^{-1}(q_h, t)), \quad i = 1, \ldots, N, \\
\psi_i[\mathbf{x}_j(t)](q_h) &= \psi_i[\mathbf{x}^0](X_h^{-1}(q_h, t)), \quad i = 1, \ldots, N_\Gamma,
\end{align}
\end{subequations}
where $q_h \in \Omega_h(t)\cup \Gamma_h(t)$. The relation in \eqref{eq:pullback basis functions} is also known as the transport property of the basis functions \cite{MCF}. It is easy to verify that this pullback preserves the trace property \eqref{eq:trace property}:
\begin{equation}
\psi_j[\mathbf{x}_j(t)] = \gamma_h(\varphi_j[\mathbf{x}_j(t)]), \quad 1 \leq j \leq N_\Gamma.
\end{equation}
The corresponding finite element spaces $\mathscr{V}_h[\mathbf{x}_j(t)], \mathscr{V}_h^0[\mathbf{x}_j(t)], \mathscr{S}_h[\mathbf{x}_j(t)]$ on the evolving domain are given similarly as in \eqref{eq:finite element spaces}.

Denoting the outward unit normal vector on $\Gamma_h(t)$ by $n_h(\cdot, t) \in [\mathscr{S}_h(\mathbf{x}(t))]^3$, for a finite element function $w_h(\cdot, t) \in \mathscr{S}_h[\mathbf{x}(t)]$, its discrete surface gradient is
\begin{equation*}
    \nabla_{\Gamma_h} w_h = \nabla w_h - (\nabla w_h \cdot n_h) n_h.
\end{equation*}
The discrete version of the surface Jacobian, surface divergence, and surface Laplace-Beltrami can be defined similarly.

For any finite element function $u_h(\cdot,t) \in \mathscr{V}_h[\mathbf{x}(t)]$ with nodal values $(u_j(t))_{j=1}^N$, its discrete material derivative at $q_h(t) = X_h(q_h^0, t) \in \Omega_h(t)$ is \cite{EKL24}
\begin{equation}
    \begin{aligned}
        \partial_h^{\bullet} u_h(q_h,t) &= \frac{d}{dt} u_h(X_h(q_h^0,t)) = \sum_{j=1}^N \dot{u}_j(t) \varphi_j[\mathbf{x}(t)](q_h),
    \end{aligned}
\end{equation}
where $\dot{u}_j(t)$ denote the time derivative $\frac{d}{dt} u_j(t)$. Similarly, for any finite element function $w_h(\cdot,t) \in \mathscr{S}_h[\mathbf{x}(t)]$ with nodal values $(w_j(t))_{j=1}^{N_\Gamma}$, its discrete material derivative at $q_h(t) = X_h(q_h^0, t) \in \Gamma_h(t)$ is
\begin{equation}
    \partial_h^{\bullet} w_h(q_h,t) = \frac{d}{dt} w_h(X_h(q_h^0,t)) = \sum_{j=1}^{N_\Gamma} \dot{w}_j(t) \psi_j[\mathbf{x}(t)](q_h).
\end{equation}
By definition, we have $\partial_h^{\bullet} ( \gamma_h (u_h )) = \gamma_h ( \partial_h^{\bullet} u_h )$.

The discrete velocity field $v_h(\cdot,t) \in [\mathscr{V}_h(\mathbf{x}(t))]^3$ for $\Omega_h(t)$ is given by the material derivative of the identity map on $\Omega_h(t)$. For notational convenience, we denote this identity map as $x_h(\cdot, t)$, which leads to:
\begin{equation}\label{eq:discrete velocity field}
    \begin{aligned}
        v_h(X_h(q^0_h,t),t) = \partial_h^{\bullet} x_h(X_h(q^0_h,t),t) = \sum_{j=1}^N \mathbf{v}_j(t) \varphi_j[\mathbf{x}^0](q_h^0),
    \end{aligned}
\end{equation}
where $\mathbf{v}_j(t) = \dot{\mathbf{x}}_j(t)$ are the nodal velocities. The discrete surface velocity field $v_{\Gamma_h}(\cdot,t)$ for $\Gamma_h(t)$ is then defined as the trace of the bulk velocity: $v_{\Gamma_h}(\cdot,t) = \gamma_h(v_h(\cdot, t)) \in [\mathscr{V}_h(\mathbf{x}(t))]^3$.

\subsubsection{Matrix–vector formulation} 
Following \cite{EKL24}, we collect the nodal values of the finite element functions $u_h(\cdot, t), v_h(\cdot, t), n_h(\cdot, t), H_h(\cdot, t)$ into the nodal vectors 
\begin{equation*}
    \mathbf{u} = (u_j) \in \mathbb{R}^N, \quad \mathbf{v} = (v_j) \in \mathbb{R}^{3N}, \quad \mathbf{n} = (n_j) \in \mathbb{R}^{3N_\Gamma}, \quad \mathbf{H} = (H_j) \in \mathbb{R}^{N_\Gamma},
\end{equation*}
respectively. Similar to the partition of $\mathbf{x}$ into $\mathbf{x}_\Gamma$ and $\mathbf{x}_\Omega$, we partition $\mathbf{u}$, $\mathbf{v}$ as:
\begin{equation*}
    \mathbf{u} = \begin{pmatrix} \mathbf{u}_\Gamma \\ \mathbf{u}_\Omega \end{pmatrix}, \quad \mathbf{v} = \begin{pmatrix} \mathbf{v}_\Gamma \\ \mathbf{v}_\Omega \end{pmatrix}.
\end{equation*}
Using this, the nodal vector of the normal velocity $V_h$ can be written as
\begin{equation*}
    \mathbf{V} = - \beta \mathbf{H} + \alpha \mathbf{u}_\Gamma.
\end{equation*}
The trace operator $\gamma_h$ corresponds to the following matrix:
\begin{equation}
    \boldsymbol{\gamma} = \begin{pmatrix} \mathbf{I}_{N_\Gamma} & \mathbf{0}  \end{pmatrix} \in \mathbb{R}^{N_\Gamma \times N},
\end{equation}
where $\mathbf{I}_{N_\Gamma}$ is the identity matrix of size $N_\Gamma \times N_\Gamma$. It is easy to verify that $\boldsymbol{\gamma} \mathbf{u} = \mathbf{u}_\Gamma$.

The bulk-dependent mass matrix $\mathbf{M}_{\bar{\Omega}}$ and stiffness matrix $\mathbf{A}_{\bar{\Omega}}$ are determined by $\mathbf{x}$ as follows:
\begin{equation}
    \mathbf{M}_{\bar{\Omega}}|_{ij} = \int_{\Omega_h[\mathbf{x}]} \varphi_i[\mathbf{x}] \varphi_j[\mathbf{x}], \, \mathbf{A}_{\bar{\Omega}}|_{ij} = \int_{\Omega_h[\mathbf{x}]} \nabla \varphi_i[\mathbf{x}] \cdot \nabla \varphi_j[\mathbf{x}], \,\forall i,j = 1, \ldots, N.
\end{equation}
The mass matrix $\mathbf{M}_{\bar{\Omega}}$ and stiffness matrix $\mathbf{A}_{\bar{\Omega}}$ can be further decomposed as
\begin{equation}
    \mathbf{M}_{\bar{\Omega}} = \begin{pmatrix} \mathbf{M}_{\Gamma\Gamma} & \mathbf{M}_{\Gamma\Omega} \\ \mathbf{M}_{\Omega\Gamma} & \mathbf{M}_{\Omega\Omega} \end{pmatrix}, \quad \mathbf{A}_{\bar{\Omega}} = \begin{pmatrix} \mathbf{A}_{\Gamma\Gamma} & \mathbf{A}_{\Gamma\Omega} \\ \mathbf{A}_{\Omega\Gamma} & \mathbf{A}_{\Omega\Omega} \end{pmatrix}.
\end{equation}
Similarly, the surface-dependent mass matrix $\mathbf{M}$ and stiffness matrix $\mathbf{A}$ are:
\begin{equation}
    \mathbf{M}|_{ij} = \int_{\Gamma_h[\mathbf{x}]} \psi_i[\mathbf{x}] \psi_j[\mathbf{x}], \, \mathbf{A}|_{ij} = \int_{\Gamma_h[\mathbf{x}]} \nabla_\Gamma \psi_i[\mathbf{x}] \cdot \nabla_\Gamma \psi_j[\mathbf{x}], \,\forall i,j = 1, \ldots, N_\Gamma.
\end{equation}
Moreover, for any dimension $d$, we have
\begin{equation*}
    \mathbf{M}_{\bar{\Omega}}^{[d]} = I_d \otimes \mathbf{M}_{\bar{\Omega}}, \quad \mathbf{M}^{[d]} = I_d \otimes \mathbf{M}, \quad \mathbf{A}_{\bar{\Omega}}^{[d]} = I_d \otimes \mathbf{A}_{\bar{\Omega}}, \quad \mathbf{A}^{[d]} = I_d \otimes \mathbf{A},
\end{equation*}
and we omit the superscript $d$ if there is no ambiguity of the dimension.

The tangential gradient matrix $\mathbf{D}(\mathbf{x}) \in \mathbb{R}^{3N_\Gamma \times N_\Gamma}$ is defined as
\begin{equation}
    \mathbf{D}(\mathbf{x})|_{i+(\ell-1)N_\Gamma,j} = \int_{\Gamma_h[\mathbf{x}]} \phi_i[\mathbf{x}] \underline{D}_{h, \ell} \psi_j[\mathbf{x}], \,\forall i,j = 1, \ldots, N_\Gamma, \ell = 1, \ldots, 3,
\end{equation}
where $\underline{D}_{h, \ell}$ is the $\ell$-th component of $\nabla_{\Gamma_h[\mathbf{x}]}$. Moreover, we introduce the following matrix for the Robin boundary condition:
\begin{equation}
    \mathbf{L}(\mathbf{x}) = \mathbf{A}_{\bar{\Omega}}(\mathbf{x}) + \alpha \boldsymbol{\gamma}^T \mathbf{M}(\mathbf{x}) \boldsymbol{\gamma}.
\end{equation}

The nonlinear functions $\mathbf{f}_u(\mathbf{x}, \mathbf{H}) \in \mathbb{R}^{N}, \mathbf{f}_n(\mathbf{x}, \mathbf{n}) \in \mathbb{R}^{3N_\Gamma}, \mathbf{f}_H(\mathbf{x}, \mathbf{n}, \mathbf{V}) \in \mathbb{R}^{N_\Gamma}$ are defined as
\begin{subequations}
    \begin{align}
        &\mathbf{f}_u(\mathbf{x}, \mathbf{H})|_{i} = -\int_{\Omega_h[\mathbf{x}]} \phi_i[\mathbf{x}] + \int_{\Gamma_h[\mathbf{x}]} (\beta H_h + Q_h) \gamma_h(\phi_i[\mathbf{x}]), \,\forall i = 1, \ldots, N, \\
        &\mathbf{f}_n(\mathbf{x}, \mathbf{n})|_{i+(\ell-1)N_\Gamma} = \beta\int_{\Gamma_h[\mathbf{x}]} |A_h|^2 (n_h)_\ell \psi_i[\mathbf{x}], \,\forall i = 1, \ldots, N_\Gamma, \ell = 1, 2, 3, \\
        &\mathbf{f}_H(\mathbf{x}, \mathbf{n}, \mathbf{V})|_{i} = -\int_{\Gamma_h[\mathbf{x}]} |A_h|^2 V_h \psi_i[\mathbf{x}], \,\forall i = 1, \ldots, N_\Gamma.
    \end{align}
\end{subequations}

Using the above notations, the semi-discretization of the weak form \eqref{eq:weak form} can be written as the following matrix-vector form:
\begin{subequations}\label{eq:semi-discretization}
    \begin{align}
        \mathbf{L}(\mathbf{x}) \mathbf{u} & = \mathbf{f}_u(\mathbf{x}, \mathbf{H}), \label{eq:robin-bvp, semi-discretization} \\
        \mathbf{M}(\mathbf{x}) \dot{\mathbf{n}} + \beta \mathbf{A}(\mathbf{x}) \mathbf{n} & = \mathbf{f}_n(\mathbf{x}, \mathbf{n}) -\alpha \mathbf{D}(\mathbf{x}) \boldsymbol{\gamma} \mathbf{u}, \label{eq:forced-mcf, semi-discretization, n} \\
        \mathbf{M}(\mathbf{x}) \dot{\mathbf{H}} + \beta \mathbf{A}(\mathbf{x}) \mathbf{H} & = \mathbf{f}_H(\mathbf{x}, \mathbf{n}, \mathbf{V}) + \alpha \mathbf{A}(\mathbf{x}) \boldsymbol{\gamma} \mathbf{u}, \label{eq:forced-mcf, semi-discretization, H} \\
        \mathbf{V} & = - \beta \mathbf{H} + \alpha \mathbf{u}_\Gamma,\label{eq:forced-mcf, semi-discretization, V} \\
        \gamma_h(v_h) & = \sum_{j=1}^{N_\Gamma} (V_j n_j) \psi_j[\mathbf{x}],\label{eq:forced-mcf, semi-discretization, v} \\
        \mathbf{A}_{\Omega\Omega}(\mathbf{x}) \mathbf{v}_\Omega & = - \mathbf{A}_{\Omega\Gamma}(\mathbf{x}) \mathbf{v}_\Gamma,\label{eq:harmonic extension, semi-discretization} \\ 
        \dot{\mathbf{x}} & = \mathbf{v}.\label{eq:ODE for flow map, semi-discretization}
    \end{align}
\end{subequations}

\subsection{Interpolated surface and lifts}
To compare functions on the numerical discrete bulk/surface $\Omega_h[\mathbf{x}(t)]/\Gamma_h[\mathbf{x}(t)]$ with functions on the exact bulk/surface $\Omega(t)/\Gamma(t)$, we introduce an interpolated discrete bulk/surface $\Omega_h[\mathbf{x}^*(t)]/\Gamma_h[\mathbf{x}^*(t)]$ with lift operators.

The evolving nodal vector $\mathbf{x}^*=(q_j^*(t))_{j=1}^{N}$ is defined via the exact flow map $X$ as $q_j^*(t) = X(q_j, t)$ for all $j = 1, \ldots, N$. This nodal vector determines the interpolated discrete bulk/surface $\Omega_h[\mathbf{x}^*]/\Gamma_h[\mathbf{x}^*]$ and the corresponding discrete finite element spaces $\mathscr{V}_h[\mathbf{x}^*]$, $\mathscr{V}_h^0[\mathbf{x}^*]$, and $\mathscr{S}_h[\mathbf{x}^*]$.

For any point $q_h^* \in \Omega_h[\mathbf{x}^*(t)]$ in the interpolated discrete bulk, as constructed in detail in \cite{ElliottRanner_bulksurface}, there exists a unique corresponding point $q\in \Omega(t)$ \cite{DziukElliott_ESFEM}. This correspondence induces a lift operator $\,^\ell$, where for any $w_h$ defined on $\Omega_h[\mathbf{x}^*(t)]$, its lift $w_h^\ell$ is a function defined on $\Omega(t)$ given by \cite{DziukElliott_ESFEM}
\begin{align}\label{eq:lift operator to exact domain}
w_h^\ell(q) = w_h(q_h^*), \quad q \in \Omega(t).
\end{align}
And the inverse lift is $\,^{-\ell}$. Specifically, we denote the lifted finite element space $(\mathscr{V}_h[\mathbf{x}^*(t)])^\ell \subset H^1(\Omega(t))$ by $\mathscr{V}_h(t)$.

Next, we consider the lift for finite element functions on $\Omega_h[\mathbf{x}(t)]$. Let $w_h = \sum w_j \phi_j[{\bf x}(t)] \in \mathscr{V}_h[{\bf x}(t)]$ be a bulk finite element function. Its lift to the interpolated discrete bulk $\Omega_h[\mathbf{x}^*(t)]$ is given by \cite{MCF}
\begin{align}\label{eq:lift operator to interpolated bulk}
    \widehat{w}_h := \sum w_j \phi_j[{\bf x}^*] \in \mathscr{V}_h[{\bf x}^*(t)].
\end{align}
Its inverse is denoted as $\widecheck{w}_h := \sum w_j \phi_j[{\bf x}(t)] \in \mathscr{V}_h[{\bf x}(t)]$. The composed lift from $\Omega_h[\mathbf{x}^*]$ to $\Omega(t)$ is $w_h^L := (\widehat{w}_h)^\ell$ \cite{MCF}. All these lift operators extend naturally to surface functions, and we refer the properties of these lift operators in \cite{ElliottRanner_unified}.

\subsection{Main results}
Under sufficient regularity assumptions on the exact solution we establish error bounds for the ESFEM semi-discretization of the original tumour model by ESFEM \eqref{eq:semi-discretization} using finite elements of polynomial degree $k \geq 3$.
\begin{theorem}\label{thm:ESFEM_error_bound}
    Suppose that the exact solution $(u,X,v,n,H)$ to the original tumour model \eqref{eq:EKS19 alt form} is sufficiently differentiable on the time interval $t \in [0,T]$, and the flow map $X(\cdot,t)$ is a homeomorphism. Then there exists a constant $\bar{h} > 0$ such that for all mesh sizes $h \leq \bar{h}$, such that the solution $(u_h,X_h,v_h,n_h,H_h)$ of the semi-discretization \eqref{eq:semi-discretization} satisfies the following error bounds:
    \begin{align}
    \left\|u^L_h(\cdot,t)-u(\cdot,t)\right\|_{H^1(\Omega(t)} &\leq Ch^{k}, \\
    \left\|x^L_h(\cdot,t)-\text{id}\right\|_{H^1(\Omega(t))} &\leq Ch^{k}, \\
    \left\|\gamma(x^L_h)(\cdot,t)-\text{id}\right\|_{H^{1/2}(\Gamma(t))} &\leq Ch^{k}, \\
    \left\|v^L_h(\cdot,t)-v(\cdot,t)\right\|_{H^1(\Omega(t))} &\leq Ch^{k}, \\
    \left\|\gamma(v^L_h)(\cdot,t)-\gamma(v)(\cdot,t)\right\|_{H^{1/2}(\Gamma(t))} &\leq Ch^{k}, \\
    \left\|n^L_h(\cdot,t)-n(\cdot,t)\right\|_{H^{1/2}(\Gamma(t))} &\leq Ch^{k}, \\
    \left\|H^L_h(\cdot,t)-H(\cdot,t)\right\|_{H^{1/2}(\Gamma(t))} &\leq Ch^{k}.
    \end{align}
    Where $x_h$ is defined in \eqref{eq:discrete velocity field}. Furthermore, for the flow map $X$ defined on $\Omega_0 \cup \Gamma_0$, we have $\left\|X^{\ell}_h(\cdot,t) - X(\cdot,t)\right\|_{H^1(\Omega_0)} \leq Ch^{k}$.

    The constant $C$ is independent of $h$ and $t$, but depends on the Sobolev norms of the exact solution $(u,X,v,n,H)$ and $T$.
\end{theorem}
\begin{remark}
    Here we employ the functional space pairing $\left(H^1(\Omega); H^{1/2}(\Gamma)\right)$ with mixed regularity to correctly reflect the boundary-domain coupling, rather than the uniform regularity pairing $\left(H^1(\Omega); H^{1}(\Gamma)\right)$. The trace inequality $\left\|\gamma(u)\right\|_{H^{1/2}(\Gamma)} \leq C \left\|u\right\|_{H^1(\Omega)}$ provides a domain-to-boundary control, while PDE regularity theory \cite[Proposition 6]{EKL24} establishes the converse boundary-to-domain estimate.
 \end{remark}

The proof of Theorem \ref{thm:ESFEM_error_bound} is obtained from stability and consistency. Its consistency has already been studied in \cite{EKL24}. The main difficulty is the stability analysis, which we present in Section 4.2. One important step of the stability analysis is a comprehensive $H^{1/2}(\Gamma)$ energy estimates theory, which we develop in the following section.

\section{\texorpdfstring{$H^{1/2}(\Gamma)$ energy estimates theory}{H(1/2) energy estimates framework}}
To establish this stability result, we first develop the necessary $H^{1/2}(\Gamma)$ energy estimates theory. We first recall the definitions of the fractional order spaces and discuss their basic properties. Then we relate the different fractional order spaces on different surfaces. Finally, we extend the derivative trick to the fractional order setting.

The analysis of fractional order typically relies on inverse estimates and the interpolation theory (Lemma \ref{lem:equivalence of norms, continuous and discrete}). Beyond these standard techniques, we introduce a new approach that exploits the explicit representation of solutions to Sylvester equations of the form $A X + X B = Y$ for the critical case $s = 1/2$ (Lemma \ref{lemma:norm equivalence for Hhs on different surfaces}, \ref{lemma:time derivative of discrete Laplacian of order 1/2}). When $A$ and $B$ are symmetric and positive definite, its solution $X$ can be explicitly written as \cite[Theorem VII.2.3]{bhatia2013matrix}
\begin{equation}\label{eq:explicit representation of solution to Sylvester equation}
    X = \int_0^\infty e^{-At} Y e^{-Bt}\, dt.
\end{equation}

Throughout, $C$ indicates a positive constant independent of $h$ and $t$. We consider three different surfaces: a regular surface $\Gamma$, its interpolation $\Gamma_h[\mathbf{x}^*]$ and a discrete surface $\Gamma_h[\mathbf{x}]$.

\subsection{Fractional order spaces and their properties}
Let $\Gamma$ be a regular surface with the Laplace-Beltrami operator $\Delta_{\Gamma}$. For any $s\geq 0$, the fractional Sobolev space $H^{s}(\Gamma)$ is defined as the domain of $(-\Delta_\Gamma)^s$, equipped with the norm
\begin{equation}
    \left\| u \right\|_{H^{s}(\Gamma)} = \left\| (I-\Delta_\Gamma)^{s/2} u \right\|_{L^2(\Gamma)}, \qquad \forall u \in H^s(\Gamma).
\end{equation}
And $H^{-s}(\Gamma)$ is the dual space of $H^s(\Gamma)$, with the dual norm
\begin{equation}
    \left\| u \right\|_{H^{-s}(\Gamma)} = \sup_{0\neq v \in H^s(\Gamma)} \frac{\int_{\Gamma} u v}{\left\| v \right\|_{H^s(\Gamma)}}, \qquad \forall u \in H^{-s}(\Gamma).
\end{equation}

Suppose $\Gamma$ be either a regular surface or a discrete surface with a finite element space $\mathscr{S}_h\subset H^1(\Gamma)$. The discrete Laplace-Beltrami operator $\Delta_{h}: \mathscr{S}_h \to \mathscr{S}_h$ is 
    \begin{equation*}
        \int_{\Gamma} \Delta_{h} u_h\, \phi_h = -\int_{\Gamma} \nabla_{\Gamma} u_h \nabla_{\Gamma} \phi_h, \qquad \forall u_h, \phi_h \in \mathscr{S}_h.
    \end{equation*}
    For any $s$, the discrete fractional Sobolev norm $\left\|\cdot\right\|_{H_h^s}$ is 
    \begin{equation}\label{eq:definition of discrete fractional Sobolev norm}
        \left\| u_h \right\|_{H_h^s} = \left\| (I-\Delta_{h})^{s/2} u_h \right\|_{L^2(\Gamma_h)}, \qquad \forall u_h \in \mathscr{S}_h.
    \end{equation}

By definition, we have $\left\| u_h \right\|_{H_h^{s_1}(\Gamma)} \leq \left\| u_h \right\|_{H_h^{s_2}(\Gamma)}, \,\forall\,  s_1 \leq s_2$. Conversely, we have the following inverse estimate, which we refer to Lemma 4.1 in \cite{hu2022optimal} for the proof.
    \begin{lemma}\label{lem:inverse estimate}
        Suppose $\Gamma$ is either a regular surface or a discrete surface. Then for any $u_h \in \mathscr{S}_h$,
        \begin{equation}\label{eq:inverse estimate}
            \left\| u_h \right\|_{H_h^{s_2}(\Gamma)} \leq C h^{s_1 - s_2} \left\| u_h \right\|_{H_h^{s_1}(\Gamma)}, \qquad \forall\,  s_1 \leq s_2.
        \end{equation}
    \end{lemma}

Next, we investigate the relationship between $H_h^s(\Gamma)$ and $H^s(\Gamma)$ for different $s$: (1)\, $0\leq s \leq 1$; (2)\, $-1\leq s \leq 0$; and (3)\, $1\leq s \leq 2$.

For $0\leq s \leq 1$, the two norms are equivalent for a finite element function $u_h \in \mathscr{S}_h$.
\begin{lemma}\label{lem:equivalence of norms, continuous and discrete}
    Suppose $\Gamma$ is a regular surface. Then there are constants $C_1, C_2 > 0$ such that for any $u_h \in \mathscr{S}_h$,
    \begin{equation}\label{eq:equivalence of norms, continuous and discrete}
        C_1 \left\| u_h \right\|_{H^s(\Gamma)} \leq \left\| u_h \right\|_{H_h^s(\Gamma)} \leq C_2 \left\| u_h \right\|_{H^s(\Gamma)}, \qquad \forall\,  0\leq s \leq 1.
    \end{equation}
\end{lemma}
\begin{proof}
    The equivalence of the two norms is straightforward to verify for the boundary cases $s = 0$ and $s = 1$. For $0 < s < 1$, we employ the interpolation argument. 
    
    We first show that the interpolation of $H_h^0(\Gamma)$ and $H_h^1(\Gamma)$ is equivalent to $H_h^s(\Gamma)$. Let $X = H^1(\Gamma)$ and $Y = L^2(\Gamma)$ be our continuous spaces, and define the discrete spaces $X_h = (\mathscr{S}_h, \|\cdot\|_{X}) = H_h^1(\Gamma)$ and $Y_h = (\mathscr{S}_h, \|\cdot\|_{Y}) = H_h^0(\Gamma)$. It is well-known that $[X, Y]_\theta = H^{1-\theta}(\Gamma)$ for $0 \leq \theta \leq 1$. To characterize $[X_h, Y_h]_\theta$, we observe that:
    \begin{equation*}
        \int_{\Gamma} u_h v_h + \int_{\Gamma} \nabla_{\Gamma} u_h \nabla_{\Gamma} v_h = \int_{\Gamma} u_h \left((I-\Delta_h) v_h\right), \qquad \forall u_h, v_h \in \mathscr{S}_h.
    \end{equation*}
    By \cite[equation (3.2)]{arioli2009discrete}, the discrete interpolation space is 
    \begin{equation*}
        [X_h, Y_h]_\theta = (\mathscr{S}_h, \|\cdot\|_{[X_h, Y_h]_\theta}),
    \end{equation*}
    where the norm is defined as
    \begin{equation*}
        \left\|u_h\right\|_{[X_h, Y_h]_\theta} = \left(\left\|u_h\right\|_{Y}^2 + \left\|(I-\Delta_h)^{(1-\theta)/2} u_h\right\|_{Y}^2\right)^{1/2}.
    \end{equation*}
    One can readily verify that $\left\|\cdot\right\|_{[X_h, Y_h]_\theta}$ and $\left\|\cdot\right\|_{H_h^{1-\theta}(\Gamma)}$ are equivalent, therefore
    \begin{equation*}
         [X_h, Y_h]_\theta = H_h^{1-\theta}(\Gamma), \qquad \forall \,0\leq \theta \leq 1.
    \end{equation*}

    Next, denote $i: \mathscr{S}_h \to L^2(\Gamma)$ as the identity map. We have already proved that $i: \mathcal{L}(X_h; X) \cap \mathcal{L}(Y_h; Y)$. Using the interpolation and norm equivalence, we have
        \begin{equation*}
            \left\| u_h \right\|_{[X, Y]_\theta} = \left\| i(u_h) \right\|_{[X, Y]_\theta} \leq C \left\| u_h \right\|_{[X_h, Y_h]_\theta}, \qquad \forall u_h \in \mathscr{S}_h, \, 0\leq \theta \leq 1.
        \end{equation*}

On the other hand, there exists an operator $\widetilde{I}_{h, 1}^{SZ}: L^2(\Gamma) \to \mathscr{S}_h$ such that
\begin{equation}\label{eq:some results of Scott-Zhang interpolation}
    \begin{aligned}
    & \left\| \widetilde{I}_{h, 1}^{SZ} u \right\|_{Y_h} \leq C \left\| \widetilde{I}_{h, 1}^{SZ} u \right\|_{L^2(\Gamma)} \leq C \left\| u \right\|_{L^2(\Gamma)}, \qquad \forall u \in L^2(\Gamma), \\
    & \left\| \widetilde{I}_{h, 1}^{SZ} u \right\|_{X_h} \leq C \left\| \widetilde{I}_{h, 1}^{SZ} u \right\|_{H^1(\Gamma)} \leq C \left\| u \right\|_{H^1(\Gamma)}, \qquad \forall u \in H^1(\Gamma), \\
    & \left\| \widetilde{I}_{h, 1}^{SZ} u - u \right\|_{L^2(\Gamma)} \leq C h \left\| u \right\|_{H^1(\Gamma)}, \qquad \forall u \in H^1(\Gamma).
    \end{aligned}
\end{equation}
$\widetilde{I}_{h, 1}^{SZ}$ is a composition of an inverse lift, a Scott-Zhang interpolation, and a lift. For the details, we refer to \cite{elliott2024sfem}.

%let $\Gamma_h^{(1)}$  a polyhedron given in []. The lift from $\Gamma_h^{(1)}$ to $\Gamma$ and its inverse are denoted as $\,^{\ell, 1}$ and $\,^{-\ell, 1}$ respectively. Let $I_{h, 1}^{SZ}: L^2(\Gamma_h^{(1)}) \to \mathscr{S}_h[{\mathbf{x}^*}]$ be the Scott-Zhang interpolation, and its lift is $\widetilde{I}_{h, 1}^{SZ}: L^2(\Gamma) \to \mathscr{S}_h$ such that
% \begin{equation*}
%     \widetilde{I}_{h, 1}^{SZ} u = (I_{h, 1}^{SZ}(u^{-\ell, 1}))^{\ell, 1}.
% \end{equation*}
% From [] and the norm equivalence, we know that 

Since $\widetilde{I}_{h, 1}^{SZ} \in \mathcal{L}(X; X_h)\cap \mathcal{L}(Y; Y_h)$, the interpolation gives
\begin{equation*}
    \begin{aligned}
    &\left\| \widetilde{I}_{h, 1}^{SZ} u \right\|_{H_h^{1-\theta}(\Gamma)} \leq C \left\| \widetilde{I}_{h, 1}^{SZ} u \right\|_{[X_h, Y_h]_\theta} \leq C \left\| u \right\|_{[X, Y]_\theta} \leq C \left\| u \right\|_{H^{1-\theta}(\Gamma)},  \\
    &\left\| \widetilde{I}_{h, 1}^{SZ} u - u \right\|_{L^2(\Gamma)} \leq C h^{1-\theta} \left\| u \right\|_{H^{1-\theta}(\Gamma)}, \qquad \forall u \in H^{1-\theta}(\Gamma).
    \end{aligned}
\end{equation*}
Finally, using the inverse estimate, we deduce that
\begin{equation*}
    \begin{aligned}
    \left\| u_h \right\|_{H_h^{1-\theta}(\Gamma)} &\leq \left\| u_h - \widetilde{I}_{h, 1}^{SZ} u_h \right\|_{H_h^{1-\theta}(\Gamma)} + \left\| \widetilde{I}_{h, 1}^{SZ} u_h \right\|_{H_h^{1-\theta}(\Gamma)}\\
    & \leq C h^{-(1-\theta)} \left\| u_h - \widetilde{I}_{h, 1}^{SZ} u_h \right\|_{L^2(\Gamma)} + C \left\| u_h \right\|_{H^{1-\theta}(\Gamma)}\leq C \left\| u_h \right\|_{H^{1-\theta}(\Gamma)}.
    \end{aligned}
\end{equation*}
The desired equivalence \eqref{eq:equivalence of norms, continuous and discrete} follows immediately
\end{proof}

For the dual space $-1\leq s \leq 0$, we have the following one side result.
\begin{lemma}\label{lem:equivalence of norms, continuous and discrete, negative order}
    Suppose $\Gamma$ is a regular surface. There is constant $C > 0$ such that for any $u_h \in \mathscr{S}_h$,
    \begin{equation}\label{eq:equivalence of norms, continuous and discrete, negative order}
        \left\| u_h \right\|_{H^{-s}(\Gamma)} \leq C \left\| u_h \right\|_{H_h^{-s}(\Gamma)}, \qquad \forall\,  0\leq s \leq 1.
    \end{equation}
\end{lemma}
\begin{proof}
    For $w \in H^s(\Gamma)$, equation \eqref{eq:some results of Scott-Zhang interpolation} gives us
    \begin{equation*}
        \left\| \widetilde{I}_{h, 1}^{SZ} w \right\|_{H^{s}(\Gamma)} \leq C \left\| w \right\|_{H^{s}(\Gamma)}, \quad \left\| \widetilde{I}_{h, 1}^{SZ} w - w \right\|_{L^2(\Gamma)} \leq C h^{s} \left\| w \right\|_{H^{s}(\Gamma)}.
    \end{equation*}

    Therefore, for any $u_h \in \mathscr{S}_h$, the dual norm gives
    \begin{align*}
        \left\|u_h\right\|_{H^{-s}(\Gamma)} & = \sup_{0\neq w \in H^s(\Gamma)} \frac{1}{\left\| w \right\|_{H^s(\Gamma)}}\int_{\Gamma} u_h w\\
         &= \sup_{0\neq w \in H^s(\Gamma)} \frac{1}{\left\| w \right\|_{H^s(\Gamma)}}\left(\int_{\Gamma} u_h \widetilde{I}_{h, 1}^{SZ} w + \int_{\Gamma} u_h (w - \widetilde{I}_{h, 1}^{SZ} w) \right)\\
          &\leq \sup_{0\neq w \in H^s(\Gamma)} \frac{1}{\left\| w \right\|_{H^s(\Gamma)}} \left\| u_h \right\|_{H_h^{-s}(\Gamma)} \left\| \widetilde{I}_{h, 1}^{SZ} w \right\|_{H_h^s(\Gamma)} \\
        &\quad + \sup_{0\neq w \in H^s(\Gamma)} \frac{1}{\left\| w \right\|_{H^s(\Gamma)}} \left\| u_h \right\|_{L^2(\Gamma)} \left\| w - \widetilde{I}_{h, 1}^{SZ} w \right\|_{L^2(\Gamma)} \\
        &\leq C \left\|u_h\right\|_{H_h^{-s}(\Gamma)} + C h^s \left\|u_h\right\|_{L^2(\Gamma)}.
    \end{align*}
    Since $\left\|u_h\right\|_{L^2(\Gamma)} \leq C h^{-s} \left\|u_h\right\|_{H_h^{-s}(\Gamma)}$ by the inverse estimate \eqref{eq:inverse estimate}, we conclude the desired inequality \eqref{eq:equivalence of norms, continuous and discrete, negative order}.
\end{proof}

When $1 < s \leq 2$, the $H^s(\Gamma)$ norm is not well-defined for a finite element function $u_h\in \mathscr{S}_h$, while we can construct its suitable approximation $u\in H^s(\Gamma)$, whose $H^s(\Gamma)$ norm can be controlled by the discrete $H_h^s(\Gamma)$ norm of $u_h$.
\begin{lemma}\label{lemma:PDE result}
    Suppose $\Gamma$ is a regular surface. For any $u_h \in \mathscr{S}_h$ and $1\leq s \leq 2$, there is a constant $C > 0$ such that
    \begin{equation}\label{eq:PDE result}
        \left\| (I - \Delta_\Gamma)^{-1}(I-\Delta_h) u_h \right\|_{H^{s}(\Gamma)} \leq C \left\| u_h \right\|_{H_h^{s}(\Gamma)}.
    \end{equation}
\end{lemma}
\begin{proof}
    Denote $w_h = (I - \Delta_h) u_h \in \mathscr{S}_h$ and $u = (I - \Delta_\Gamma)^{-1} w_h$. We know that $u$ is a solution of the following weak formulation
    \begin{equation*}
        \int_\Gamma u \phi + \int_\Gamma \nabla_\Gamma u \cdot \nabla_\Gamma \phi = \int_\Gamma w_h \phi, \quad \forall \phi \in H^1(\Gamma).
    \end{equation*}
    Using the PDE regularity result Theorem A.1 in \cite{KovacsPower_quasilinear} and taking $\phi = u$, we have
    \begin{equation*}
        \left\|u\right\|_{H^2(\Gamma)} \leq C \left\|w_h\right\|_{L^2(\Gamma)} + C \left\|u\right\|_{L^2(\Gamma)} \leq C \left\|w_h\right\|_{L^2(\Gamma)} \leq C \left\|(I - \Delta_h)^{1/2} u_h \right\|_{H_h^1(\Gamma)}.
    \end{equation*}

    On the other hand, we find $u_h$ is a solution of the following weak formulation
    \begin{equation*}
        \int_{\Gamma} u_h \phi_h + \int_{\Gamma} \nabla_{\Gamma} u_h \cdot \nabla_{\Gamma} \phi_h = \int_{\Gamma} w_h \phi_h, \quad \forall \phi_h \in \mathscr{S}_h.
    \end{equation*}
    We know that $u_h = R_h u$, where $R_h$ is the Ritz projection. The property of the Ritz projection gives us \cite{highorderESFEM}
    \begin{equation*}
        \left\|u_h - u\right\|_{H^1(\Gamma)} = \left\|R_h u - u\right\|_{H^1(\Gamma)} \leq C h \left\|u\right\|_{H^2(\Gamma)} \leq C h \left\|u_h\right\|_{H_h^2(\Gamma)}.
    \end{equation*}
    This, together with the inverse estimate $\left\|u_h\right\|_{H_h^2(\Gamma)} \leq C h^{-1} \left\|u_h\right\|_{H_h^1(\Gamma)}$ yields
    \begin{equation*}
        \left\|u\right\|_{H^1(\Gamma)} \leq \left\|u_h\right\|_{H^1(\Gamma)} + \left\|u_h - u\right\|_{H^1(\Gamma)} \leq \left\|u_h\right\|_{H_h^1(\Gamma)} \leq \left\|(I - \Delta_h)^{1/2} u_h \right\|_{H_h^0(\Gamma)}.
    \end{equation*}
    \eqref{eq:PDE result} follows by the interpolation argument with the fact $\left\|(I - \Delta_h)^{1/2} u_h \right\|_{H_h^s(\Gamma)} = \left\|u_h \right\|_{H_h^{s+1}(\Gamma)}$.
\end{proof}

\subsection{Estimates for integrations over discrete surfaces}
To motivate the results, we consider the following multilinear form defined on a discrete surface $\Gamma_h[\mathbf{x}]$:
\begin{equation*}
    \int_{\Gamma_h[\mathbf{x}]} \nabla_{\Gamma_h[{\mathbf{x}}]} u_h \cdot \nabla_{\Gamma_h[{\mathbf{x}}]} w_h \, \nabla_{\Gamma_h[{\mathbf{x}}]} \cdot v_h,
\end{equation*}
where $u_h, w_h \in \mathscr{S}_h[\mathbf{x}], v_h \in [\mathscr{S}_h[\mathbf{x}]]^3$. Corollary to the $L^{2}$-$L^{\infty}$-$L^2$ estimate by H\"older's inequality in \cite{MCF}, we want to get the following $H^{1/2}$-$W^{1/2, \infty}$-$H^{-1/2}$ estimate:
\begin{equation*}
    \begin{aligned}
        &\int_{\Gamma_h[\mathbf{x}]} \nabla_{\Gamma_h[{\mathbf{x}}]} u_h \cdot \nabla_{\Gamma_h[{\mathbf{x}}]} w_h \, \nabla_{\Gamma_h[{\mathbf{x}}]} \cdot v_h\\
         &\leq C \left\| \nabla_{\Gamma_h[\mathbf{x}]} u_h \right\|_{H^{1/2}(\Gamma_h)} \left\| \nabla_{\Gamma_h[\mathbf{x}]} w_h \right\|_{H^{1/2}(\Gamma_h)} \left\| \nabla_{\Gamma_h[\mathbf{x}]} v_h \right\|_{H^{-1/2}(\Gamma_h)} \\
    & \leq C \left\| u_h \right\|_{H^{3/2}(\Gamma_h)} \left\| w_h \right\|_{W^{3/2, \infty}(\Gamma_h)} \left\| v_h \right\|_{H^{1/2}(\Gamma_h)}.
    \end{aligned}
\end{equation*}
Since the $H^{3/2}$ and $W^{3/2, \infty}$ norms are only well-defined for regular functions on a regular surface $\Gamma$, we employ a three-step approach:
\begin{enumerate}
    \item Move the integration to a regular surface $\Gamma$ via two comparisons:
    \begin{equation*}
        \begin{aligned}
        &\int_{\Gamma_h[\mathbf{x}]} \nabla_{\Gamma_h[\mathbf{x}]} u_h \cdot \nabla_{\Gamma_h[\mathbf{x}]} w_h \, \nabla_{\Gamma_h[\mathbf{x}]}\cdot v_h - \int_{\Gamma_h[\mathbf{x}^*]} \nabla_{\Gamma_h[\mathbf{x}^*]} \widehat{u}_h \cdot \nabla_{\Gamma_h[\mathbf{x}^*]} \widehat{w}_h \, \nabla_{\Gamma_h[\mathbf{x}^*]} \cdot \widehat{v}_h, \\
        &\int_{\Gamma_h[\mathbf{x}^*]} \nabla_{\Gamma_h[\mathbf{x}^*]} \widehat{u}_h \cdot \nabla_{\Gamma_h[\mathbf{x}^*]} \widehat{w}_h \, \nabla_{\Gamma_h[\mathbf{x}^*]} \cdot \widehat{v}_h - \int_{\Gamma} \nabla_{\Gamma} \widehat{u}_h^\ell \cdot \nabla_{\Gamma} \widehat{w}_h^\ell \, \nabla_{\Gamma} \cdot \widehat{v}_h^\ell.
        \end{aligned}
    \end{equation*}
    \item Approximate the finite element functions with regular functions
    \begin{equation*}
        \int_{\Gamma} \nabla_{\Gamma} \widehat{u}_h^\ell \cdot \nabla_{\Gamma} \widehat{w}_h^\ell \, \nabla_{\Gamma} \cdot \widehat{v}_h^\ell - \int_{\Gamma} \nabla_{\Gamma} u^* \cdot \nabla_{\Gamma} w^* \, \nabla_{\Gamma} \cdot v^*,
    \end{equation*}
    and thus use the $H^{1/2}$-$W^{1/2, \infty}$-$H^{-1/2}$ estimate to derive
    \begin{equation*}
        \begin{aligned}
        \int_{\Gamma} \nabla_{\Gamma} u^* \cdot \nabla_{\Gamma} w^* \, \nabla_{\Gamma} \cdot v^* &\leq C \left\| \nabla_\Gamma u^* \right\|_{H^{1/2}(\Gamma)} \left\| \nabla_\Gamma w^* \right\|_{W^{1/2, \infty}(\Gamma)} \left\| \nabla_\Gamma v^* \right\|_{H^{-1/2}(\Gamma)} \\
        &\leq C \left\| u^* \right\|_{H^{3/2}(\Gamma)} \left\| w^* \right\|_{W^{3/2, \infty}(\Gamma)} \left\| v^* \right\|_{H^{1/2}(\Gamma)}\\
        &\leq C \left\| \widehat{u}_h^\ell \right\|_{H_h^{3/2}(\Gamma)} \left\| w^* \right\|_{W^{3/2, \infty}(\Gamma)} \left\| \widehat{v}_h^\ell \right\|_{H_h^{1/2}(\Gamma)}.
        \end{aligned}
    \end{equation*}
    \item Pullback the $H_h^s$ norm by norm equivalence between $\Gamma$ and $\Gamma_h[\mathbf{x}^*]$
    \begin{equation*}
        \left\| \widehat{u}_h^\ell \right\|_{H_h^{s}(\Gamma)}  \leq C \left\| \widehat{u}_h\right\|_{H_h^{s}(\Gamma_h[\mathbf{x}^*])}.
    \end{equation*}
\end{enumerate}
Each step will be discussed in detail in the following.

We first examine the perturbation of functions on a fixed surface.
\begin{lemma}[Function perturbation]\label{lemma:comparison of different functions on the same surface}
    Let $\Gamma$ be either a regular surface or a discrete surface. Consider functions $u_1, u_2 \in L^2(\Gamma)$ and $v_1, \ldots, v_n \in L^\infty(\Gamma)$ together with their perturbations $u_1^*, u_2^* \in L^2(\Gamma)$ and $v_1^*, \ldots, v_n^* \in L^{\infty}(\Gamma)$. Assume there is a constant $C > 0$ such that $\left\|u_i^*\right\|_{L^2(\Gamma)} \leq C \left\|u_i\right\|_{L^2(\Gamma)}$ for $i = 1, 2$, and $\left\|v_j\right\|_{L^{\infty}(\Gamma)} \leq C$, $\left\|v_j^*\right\|_{L^{\infty}(\Gamma)} \leq C$ for $j = 1, \ldots, n$.
    
    Then for any multilinear form $T$ with constant coefficients and smooth function $g = g(v_1, \ldots, v_n)$, the following function perturbation estimate holds:
    \begin{equation}
        \begin{aligned}
        &\int_{\Gamma} T(u_1, u_2, g) - \int_{\Gamma} T(u_1^*, u_2^*, g^*) \\
        & \leq C \left(\left\|u_1 - u_1^*\right\|_{L^2(\Gamma)} \left\|u_2\right\|_{L^2(\Gamma)} + \left\|u_2 - u_2^*\right\|_{L^2(\Gamma)} \left\|u_1\right\|_{L^2(\Gamma)}\right.\\
        & \qquad + \left. \left\|u_1\right\|_{L^2(\Gamma)} \left\|u_2\right\|_{L^2(\Gamma)}\sum_{i=1}^n \left\|v_i - v_i^*\right\|_{L^{\infty}(\Gamma)}  \right).
        \end{aligned}
    \end{equation}
    Where $g^* := g(v_1^*, \ldots, v_n^*)$.
\end{lemma}

\begin{proof}
    Denote $u_i^\theta = u_i^* + \theta \left(u_i - u_i^*\right)$, $i = 1, 2$, and $v_j^\theta = v_j^* + \theta \left(v_j - v_j^*\right)$, $j = 1, \ldots, n$. $g^\theta = g(v_1^\theta, \ldots, v_n^\theta)$. From the assumption, for $i=1, 2$ and $j=1, \ldots, n$, we have
    \begin{equation*}
        \begin{aligned}
        &\left\|u_i^\theta\right\|_{L^2(\Gamma)} \leq C \left\|u_i\right\|_{L^2(\Gamma)},  \quad \left\|v_j^\theta\right\|_{L^\infty(\Gamma)} \leq C, \quad \partial_{\theta} u_i^\theta = u_i - u_i^*, \quad \partial_{\theta} v_j^\theta = v_j - v_j^*, \\
        &\left\|\partial_{\theta} g^\theta\right\|_{L^\infty(\Gamma)} \leq \sum_{j=1}^n \left\|\partial_{j} g^\theta\right\|_{L^\infty(\Gamma)} \left\|\partial_{\theta} v_j^\theta\right\|_{L^\infty(\Gamma)} \leq C \sum_{j=1}^n \left\|v_j - v_j^*\right\|_{L^\infty(\Gamma)}.
        \end{aligned}
    \end{equation*}
    Using the Leibniz formula and the $L^{2}$-$L^2$-$L^{\infty}$ estimate, we have
    \begin{equation*}
        \begin{aligned}
        &\int_{\Gamma} T(u_1, u_2, g) - \int_{\Gamma} T(u_1^*, u_2^*, g^*) \\
        & = \int_0^1 d\theta \int_{\Gamma} \frac{d}{d\theta} T(u_1^\theta, u_2^\theta, g^\theta)  \\
        & = \int_0^1 d\theta \int_{\Gamma} \left(T(u_1 - u_1^*, u_2^\theta, g^\theta) + T(u_1^\theta, u_2 - u_2^*, g^\theta) + T(u_1^\theta, u_2^\theta, \partial_{\theta} g^\theta)\right) \\
        & \leq C \left(\left\|u_1 - u_1^*\right\|_{L^2(\Gamma)} \left\|u_2\right\|_{L^2(\Gamma)} + \left\|u_2 - u_2^*\right\|_{L^2(\Gamma)} \left\|u_1\right\|_{L^2(\Gamma)}\right.\\
        & \qquad + \left. \left\|u_1\right\|_{L^2(\Gamma)} \left\|u_2\right\|_{L^2(\Gamma)}\sum_{i=1}^n \left\|v_i - v_i^*\right\|_{L^{\infty}(\Gamma)}  \right).
        \end{aligned}
    \end{equation*}
\end{proof}

We then consider the perturbation of two discrete surfaces, and the perturbation of a regular surface with its interpolated surface.

\begin{lemma}[Surface perturbation I]\label{lemma:comparison of discrete surface and interpolated surface}
    Let $\Gamma_h[\mathbf{x}]$ and $\Gamma_h[\mathbf{x}^*]$ be two discrete surfaces with geometric perturbation satisfying $\left\|e_x\right\|_{W^{1,\infty}(\Gamma_h[\mathbf{x}^*])} \leq \frac{1}{4}$. Consider functions $u_1, u_2, v_1, \ldots, v_n \in \mathscr{S}_h[\mathbf{x}]$. Assume there is a constant $C > 0$ such that $\left\|v_j\right\|_{W^{1,\infty}(\Gamma_h[\mathbf{x}])} \leq C$ for $j = 1, \ldots, n$.

    Then for any multilinear form $T$ with constant coefficients, $f_i(u_i, \nabla_{\Gamma_h[\mathbf{x}]} u_i) \in \{u_i, \nabla_{\Gamma_h[\mathbf{x}]} u_i\}$ for $i = 1, 2$, and smooth function $g = g(v_1, \nabla_{\Gamma_h[\mathbf{x}]} v_1, \ldots, \nabla_{\Gamma_h[\mathbf{x}]} v_n)$, the following surface perturbation estimate holds:
    \begin{equation}\label{eq:comparison of discrete surface and interpolated surface}
        \begin{aligned}
        &\left|\int_{\Gamma_h[\mathbf{x}]} T(f_1, f_2, g) - \int_{\Gamma_h[\mathbf{x}^*]} T(\widehat{f}_1, \widehat{f}_2, \widehat{g})\right| \\
        &\leq C \left\|e_x\right\|_{W^{1, \infty}(\Gamma_h[\mathbf{x}^*])} \left\|\widehat{f}_1\right\|_{L^2(\Gamma_h[\mathbf{x}^*])} \left\|\widehat{f}_2\right\|_{L^2(\Gamma_h[\mathbf{x}^*])},
        \end{aligned}
    \end{equation}
    where $f_i$, $\widehat{g}$, $\widehat{f}_i$, are defined similarly for $i = 1, 2$.
\end{lemma}
\begin{proof}
    We adopt the notations $\Gamma_h^\theta= \Gamma_h[\mathbf{x}^* + \theta \mathbf{e}]$, $w_h^\theta = \sum w_j \phi_j[\mathbf{x}^* + \theta \mathbf{e}]$ in \cite{KLLP2017}, where $\mathbf{e} = \mathbf{x} - \mathbf{x}^*$. Since $\left\|e_x\right\|_{W^{1,\infty}(\Gamma_h[\mathbf{x}]^*)} \leq \frac{1}{4}$, there is a constant $C > 0$ such that \cite{KLLP2017}
    \begin{equation*}
        \left\|w_h^\theta\right\|_{W^{k, p}(\Gamma_h^\theta)} \leq C \left\|w_h\right\|_{W^{k, p}(\Gamma_h[\mathbf{x}^*])}, \quad k = 0, 1, \, 2 \leq p \leq \infty.
    \end{equation*}
    Therefore, from the definition of $f_i$ and $g$, we have
    \begin{equation*}
        \begin{aligned}
        &\left\|f_i^\theta\right\|_{L^2(\Gamma_h^\theta)} \leq C \left\|\widehat{f}_i\right\|_{L^2(\Gamma_h[\mathbf{x}^*])} , \quad \left\|f_i^\theta\right\|_{L^\infty(\Gamma_h^\theta)} \leq C \left\|\widehat{f}_i\right\|_{L^\infty(\Gamma_h[\mathbf{x}^*])}, \quad i = 1, 2, \\
        &\left\|g^\theta\right\|_{L^\infty(\Gamma_h^\theta)} \leq C \left\|\widehat{g}\right\|_{L^\infty(\Gamma_h[\mathbf{x}^*])}.
        \end{aligned}
    \end{equation*}
    Here $f_1^\theta=f_1(u_1^\theta, \nabla_{\Gamma_h^\theta} u_1^\theta)$, $f_2^\theta$ and $g^\theta$ are defined similarly.

    From \cite[Section 7.3]{MCF}, for any $w_h \in \mathscr{S}_h[\mathbf{x}]$ and smooth function $f(w_h, \nabla_{\Gamma_h[\mathbf{x}]} w_h)$, we have
    \begin{equation*}
        \partial_{\theta}^\bullet w_h^\theta = 0, \quad \partial_{\theta}^\bullet \nabla_{\Gamma_h^\theta} w_h^\theta = -\left(\nabla_{\Gamma_h^\theta} e_x^\theta - n_h^\theta (n_h^\theta)^T (\nabla_{\Gamma_h^\theta} e_x^\theta)^T\right) \nabla_{\Gamma_h^\theta} w_h^\theta.
    \end{equation*}
    \begin{equation*}
        \partial_{\theta}^\bullet f(w_h^\theta, \nabla_{\Gamma_h^\theta} w_h^\theta) = \partial_1 f(w_h^\theta, \nabla_{\Gamma_h^\theta} w_h^\theta) \partial_{\theta}^\bullet w_h^\theta + \partial_2 f(w_h^\theta, \nabla_{\Gamma_h^\theta} w_h^\theta) \partial_{\theta}^\bullet \nabla_{\Gamma_h^\theta} w_h^\theta.
    \end{equation*}
    Where $n_h^\theta$ is the unit normal vector of $\Gamma_h^\theta$. Therefore, we obtain
    \begin{equation*}
        \begin{aligned}
        \left\|\partial_{\theta}^\bullet f_i^\theta\right\|_{L^2(\Gamma_h^\theta)} &\leq C \left\|e_x^\theta\right\|_{W^{1, \infty}(\Gamma_h^\theta)} \left\|f_i^\theta\right\|_{L^2(\Gamma_h^\theta)} \\
        & \leq C \left\|e_x\right\|_{W^{1, \infty}(\Gamma_h[\mathbf{x}^*])} \left\|\widehat{f}_i\right\|_{L^2(\Gamma_h[\mathbf{x}^*])}.
        \end{aligned}
    \end{equation*}
    \begin{equation*}
        \begin{aligned}
        \left\|\partial_{\theta}^\bullet g^\theta\right\|_{L^\infty(\Gamma_h^\theta)} &\leq C \left\|e_x^\theta\right\|_{W^{1, \infty}(\Gamma_h^\theta)}\sum_{i=1}^n \left\|\partial_{2i-1} g^\theta\right\|_{L^\infty(\Gamma_h^\theta)} \left\|v_i^\theta\right\|_{W^{1, \infty}(\Gamma_h^\theta)} \\
        & \leq C \left\|e_x\right\|_{W^{1, \infty}(\Gamma_h[\mathbf{x}^*])}.
        \end{aligned}
    \end{equation*}

    Using the Leibniz formula and the $L^{2}$-$L^2$-$L^{\infty}$ estimate, we deduce that
    \begin{equation*}
        \begin{aligned}
        &\int_{\Gamma_h[\mathbf{x}]} T\left(f_1, f_2, g\right) - \int_{\Gamma_h[\mathbf{x}^*]} T\left(\widehat{f}_1, \widehat{f}_2, \widehat{g}\right) \\
        %&= \int_0^1 \frac{d}{d\theta} \int_{\Gamma_h^\theta} T\left(f_1^\theta, f_2^\theta, g^\theta\right) d\theta \\
        & = \int_0^1 \int_{\Gamma_h^\theta} T\left(\partial_{\theta}^\bullet f_1^\theta, f_2^\theta,  g^\theta\right) d\theta   + \int_0^1 \int_{\Gamma_h^\theta} T\left(f_1^\theta, \partial_{\theta}^\bullet f_2^\theta,  g^\theta\right) d\theta \\
        & \qquad + \int_0^1 \int_{\Gamma_h^\theta} T\left(f_1^\theta, f_2^\theta, \partial_{\theta}^\bullet g^\theta\right) d\theta + \int_0^1 \int_{\Gamma_h^\theta} T\left(f_1^\theta, f_2^\theta, g^\theta\right) \nabla_{\Gamma_h^\theta} \cdot e_x^\theta d\theta \\
        & \leq C \left\|e_x\right\|_{W^{1, \infty}(\Gamma_h[\mathbf{x}^*])} \left\|\widehat{f}_1\right\|_{L^2(\Gamma_h[\mathbf{x}^*])} \left\|\widehat{f}_2\right\|_{L^2(\Gamma_h[\mathbf{x}^*])},
        \end{aligned}
    \end{equation*}
    which is the desired result.
\end{proof}

\begin{lemma}[Surface perturbation II]\label{lemma:comparison of interpolated surface and exact surface}
    Let $\Gamma$ be a regular surface and $\Gamma_h[\mathbf{x}^*]$ its discrete approximation of degree $k$. Consider functions $u_1, u_2, v_1, \ldots, v_n \in \mathscr{S}_h[\mathbf{x}^*]$. Assume there is a constant $C > 0$ such that $\left\|v_j\right\|_{W^{1,\infty}(\Gamma_h[\mathbf{x}^*])} \leq C$ for $j = 1, \ldots, n$.

    Then for any multilinear form $T$ with constant coefficients, $f_i(u_i, \nabla_{\Gamma_h[\mathbf{x}^*]} u_i) \in \{u_i, \nabla_{\Gamma_h[\mathbf{x}^*]} u_i\}$ for $i = 1, 2$, and smooth function $g = g(v_1, \nabla_{\Gamma_h[\mathbf{x}^*]} v_1, \ldots, \nabla_{\Gamma_h[\mathbf{x}^*]} v_n)$, the following surface perturbation estimate holds:
    \begin{equation}\label{eq:comparison of interpolated surface and exact surface}
        \begin{aligned}
        &\left|\int_{\Gamma_h[\mathbf{x}^*]} T(f_1, f_2, g) - \int_{\Gamma} T(f_1^\ell, f_2^\ell, g^\ell)\right| \\
        &\leq C h^k \left\|f_1\right\|_{L^2(\Gamma_h[\mathbf{x}^*])} \left\|f_2\right\|_{L^2(\Gamma_h[\mathbf{x}^*])},
        \end{aligned}
    \end{equation}
    where $f_i$, $g^\ell$, $f_i^\ell$, are defined similarly for $i = 1, 2$. 
\end{lemma}

\begin{proof}
    Using the norm equivalence by lift \cite[Proposition 7.8]{ElliottRanner_unified}, we have that for any $m = 0, 1$ and $2 \leq p \leq \infty$, there exist constants $C_1, C_2 > 0$ such that
    \begin{equation*}
        C_1 \left\|w\right\|_{W^{m, p}(\Gamma_h[\mathbf{x}^*])} \leq \left\|w^\ell\right\|_{W^{m, p}(\Gamma)} \leq C_2 \left\|w\right\|_{W^{m, p}(\Gamma_h[\mathbf{x}^*])}, \quad \forall w \in W^{m, p}(\Gamma_h[\mathbf{x}^*]).
    \end{equation*}
    In particular, this implies $\left\|v_j^\ell\right\|_{W^{1, \infty}(\Gamma)} \leq C$ for all $j = 1, \ldots, n$.

    By change of variables, the integral $\int_{\Gamma_h[\mathbf{x}^*]} T(f_1, f_2, g)$ can be written as an integration over $\Gamma$ \cite{elliott2024sfem}:
    \begin{equation*}
        \int_{\Gamma} \frac{1}{\mu_h} T\left(f_1(u_1^\ell, \boldsymbol{B}_h \nabla_\Gamma u_1^\ell), f_2(u_2^\ell, \boldsymbol{B}_h \nabla_\Gamma u_2^\ell), g(v_1^\ell, \boldsymbol{B}_h \nabla_\Gamma v_1^\ell, \ldots, v_n^\ell, \boldsymbol{B}_h \nabla_\Gamma v_n^\ell)\right).
    \end{equation*}
    Here, $\mu_h$ and $\boldsymbol{B}_h$ are defined in \cite{elliott2024sfem} and satisfy the properties
    \begin{equation*}
        \left\|1-\mu_h\right\|_{L^\infty(\Gamma)} \leq C h^{k+1}, \quad \left\|\boldsymbol{B}_h\right\|_{L^\infty(\Gamma)} \leq C, \quad \left\|\boldsymbol{B}_h - \boldsymbol{P}_\Gamma\right\|_{L^\infty(\Gamma)} \leq C h^k,
    \end{equation*}
    where $\boldsymbol{P}_\Gamma \nabla_\Gamma w = \nabla_\Gamma w$ for all $w \in H^1(\Gamma)$.

    Next, we examine the property of $f_i$ and $g$. If $f_1(u_1, \nabla_\Gamma u_1) = u_1$, then
    \begin{equation*}
        \left\|f_1^\ell - f_1(u_1^\ell, \boldsymbol{B}_h \nabla_\Gamma u_1^\ell)\right\|_{L^2(\Gamma)} = 0 \leq C h^k \left\|f_1\right\|_{L^2(\Gamma_h[\mathbf{x}^*])}.
    \end{equation*}
    If $f_1(u_1, \nabla_\Gamma u_1) = \nabla_\Gamma u_1$, then
    \begin{align*}
        \left\|f_1^\ell - f_1(u_1^\ell, \boldsymbol{B}_h \nabla_\Gamma u_1^\ell)\right\|_{L^2(\Gamma)} &= \left\|\boldsymbol{P}_\Gamma \nabla_\Gamma u_1^\ell - \boldsymbol{B}_h \nabla_\Gamma u_1^\ell\right\|_{L^2(\Gamma)} \\
        &\leq C h^k \left\|\nabla_\Gamma u_1^\ell\right\|_{L^2(\Gamma)} \leq C h^k \left\|f_1\right\|_{L^2(\Gamma_h[\mathbf{x}^*])}.
    \end{align*}
    Therefore, $\left\|f_1^\ell - f_1(u_1^\ell, \boldsymbol{B}_h \nabla_\Gamma u_1^\ell)\right\|_{L^2(\Gamma)} \leq C h^k \left\|f_1\right\|_{L^2(\Gamma_h[\mathbf{x}^*])}$. The same estimate holds for $f_2$.

    For the terms in $g$, we have
    \begin{equation*}
        \left\|\nabla_\Gamma v_i^\ell - \boldsymbol{B}_h \nabla_\Gamma v_i^\ell\right\|_{L^\infty(\Gamma)} = \left\|(\boldsymbol{P}_\Gamma - \boldsymbol{B}_h) \nabla_\Gamma v_i^\ell\right\|_{L^\infty(\Gamma)} \leq C h^k.
    \end{equation*}
    Therefore, $\left\|\nabla_\Gamma v_i^\ell\right\|_{L^\infty(\Gamma)} \leq C$ and $\left\|\boldsymbol{B}_h \nabla_\Gamma v_i^\ell\right\|_{L^\infty(\Gamma)} \leq C$.

    Combining these estimates with Lemma \ref{lemma:comparison of different functions on the same surface} and the norm equivalence by lift, we finally obtain
    \begin{equation*}
        \begin{aligned}
        &\int_{\Gamma_h[\mathbf{x}^*]} T(f_1, f_2, g) - \int_{\Gamma} T(f_1^\ell, f_2^\ell, g^\ell) \\
        & = \int_{\Gamma_h[\mathbf{x}^*]} (1-\mu_h)T(f_1, f_2, g) + \left(\int_{\Gamma_h[\mathbf{x}^*]} \mu_h T(f_1, f_2, g) - \int_{\Gamma} T(f_1^\ell, f_2^\ell, g^\ell)\right) \\
        & \leq C \left\|1-\mu_h\right\|_{L^\infty(\Gamma)} \left\|f_1\right\|_{L^2(\Gamma_h[\mathbf{x}^*])} \left\|f_2\right\|_{L^2(\Gamma_h)} \left\|g\right\|_{L^\infty(\Gamma_h[\mathbf{x}^*])} \\
        & \qquad + C \left\|f_1^\ell - f_1(u_1^\ell, \boldsymbol{B}_h \nabla_\Gamma u_1^\ell)\right\|_{L^2(\Gamma)} \left\|f_2^\ell\right\|_{L^2(\Gamma)} \\
        & \qquad + C \left\|f_1^\ell\right\|_{L^2(\Gamma)} \left\|f_2^\ell - f_2(u_2^\ell, \boldsymbol{B}_h \nabla_\Gamma u_2^\ell)\right\|_{L^2(\Gamma)}  \\
        & \qquad + C \left\|f_1^\ell\right\|_{L^2(\Gamma)} \left\|f_2^\ell\right\|_{L^2(\Gamma)} \sum_{i=1}^n \left\|\nabla_\Gamma v_i^\ell - \boldsymbol{B}_h \nabla_\Gamma v_i^\ell\right\|_{L^\infty(\Gamma)} \\
        & \leq C (h^{k+1} + h^k) \left\|f_1\right\|_{L^2(\Gamma_h[\mathbf{x}^*])} \left\|f_2\right\|_{L^2(\Gamma_h[\mathbf{x}^*])} \\
        & \leq C h^k \left\|f_1\right\|_{L^2(\Gamma_h[\mathbf{x}^*])} \left\|f_2\right\|_{L^2(\Gamma_h[\mathbf{x}^*])}.
        \end{aligned}
    \end{equation*}
\end{proof}

Next, we establish the $H^{1/2}$-$H^{-1/2}$-$W^{1/2, \infty}$ estimate for finite element functions on the smooth surface, which plays the same role as the $L^2$-$L^2$-$L^\infty$ estimate in \cite{KLLP2017}.

\begin{lemma}[$H^{1/2}$-$H^{-1/2}$-$W^{1/2, \infty}$ estimate]\label{lemma:comparison of interpolated surface and exact surface 2}
    Let $\Gamma$ be a regular surface.  Consider functions $u_1, u_2\in \mathscr{S}_h$, $v_1, \ldots, v_n \in \mathscr{S}_h$ with their perturbations $v_1^*, \ldots, v_n^* \in W^{3/2, \infty}(\Gamma)$. Assume there are constants $C > 0, \alpha_i>1/2$, such that $\left\|v_i\right\|_{W^{1,\infty}(\Gamma)} \leq C$, $\left\|v_i^*\right\|_{W^{3/2,\infty}(\Gamma)} \leq C$, and $\left\|v_i - v_i^*\right\|_{W^{1, \infty(\Gamma)}} \leq C h^{\alpha_i}$, for all $i = 1, \ldots, n$. 

    Then for any multilinear form $T$ with constant coefficients, and smooth function $g$ defined as in the previous lemma, we have
    \begin{subequations}\label{eq:comparison of interpolated surface and exact surface 2}
        \begin{align}
        & \int_{\Gamma} T(u_1, u_2, g) \leq C \left\|u_1\right\|_{H_h^{1/2}(\Gamma)} \left\|u_2\right\|_{H_h^{-1/2}(\Gamma)}, \label{eq:comparison of interpolated surface and exact surface 2, case 1} \\
        & \int_{\Gamma} T(\nabla_\Gamma u_1, u_2, g) \leq C \left\| u_1\right\|_{H_h^{3/2}(\Gamma)} \left\|u_2\right\|_{H_h^{-1/2}(\Gamma)},\label{eq:comparison of interpolated surface and exact surface 2, case 2} \\
        & \int_{\Gamma} T(u_1, \nabla_\Gamma u_2, g) \leq C \left\|u_1\right\|_{H_h^{1/2}(\Gamma)} \left\| u_2\right\|_{H_h^{1/2}(\Gamma)},\label{eq:comparison of interpolated surface and exact surface 2, case 3} \\
        & \int_{\Gamma} T(\nabla_\Gamma u_1, \nabla_\Gamma u_2, g) \leq C \left\|u_1\right\|_{H_h^{3/2}(\Gamma)} \left\| u_2\right\|_{H_h^{1/2}(\Gamma)}. \label{eq:comparison of interpolated surface and exact surface 2, case 4}
        \end{align}
    \end{subequations}
\end{lemma}

\begin{proof}
    We only prove the last inequality, the other inequalities can be proved similarly. Since $g$ is a smooth function, we have
    \begin{subequations}
        \begin{align*}
        & \left\|g(w_1^*, \nabla_{\Gamma} w_1^*, \ldots, w_n^*, \nabla_{\Gamma} w_n^*) - g(0)\right\|_{L^{\infty}(\Gamma)} \leq C \prod_{i=1}^n \left\|w_i^*\right\|^2_{W^{1, \infty}(\Gamma)},  \forall w_i^* \in W^{1, \infty}(\Gamma),  \\
        & \left\|g(w_1^*, \nabla_{\Gamma} w_1^*, \ldots, w_n^*, \nabla_{\Gamma} w_n^*) - g(0)\right\|_{W^{1, \infty}(\Gamma)} \leq C\prod_{i=1}^n \left\|w_i^*\right\|^2_{W^{2, \infty}(\Gamma)},  \forall w_i^* \in W^{2, \infty}(\Gamma).
        \end{align*}
    \end{subequations}
    By interpolation and the fact $\left\|v_i^*\right\|_{W^{3/2, \infty}(\Gamma)} \leq C$, $|g(0)| \leq C$, we obtain that
    \begin{equation*}
        \left\|g(v_1^*, \nabla_{\Gamma} v_1^*, \ldots, v_n^*, \nabla_{\Gamma} v_n^*)\right\|_{W^{1/2, \infty}(\Gamma)} \leq C.
    \end{equation*}
    
    Next, let $u_1^* = (I-\Delta_\Gamma)^{-1}(I-\Delta_\Gamma)u_1$. Then Lemma \ref{lemma:PDE result} gives that $\left\|u_1^*\right\|_{H^{3/2}(\Gamma)} \leq C \left\|u_1\right\|_{H_h^{3/2}(\Gamma)}$ and
    \begin{equation*}
        \left\|u_1 - u_1^*\right\|_{L^2(\Gamma)} + h \left\|\nabla_\Gamma u_1 - \nabla_\Gamma u_1^*\right\|_{L^2(\Gamma)} \leq C h^{2} \left\|u_1\right\|_{H_h^2(\Gamma)}.
    \end{equation*}
    For $T(\nabla_\Gamma u_1^*, \nabla_\Gamma u_2, g(v_1^*, \nabla_{\Gamma} v_1^*, \ldots, v_n^*, \nabla_{\Gamma} v_n^*))$, using the $H^{1/2}$-$H^{-1/2}$-$W^{1/2, \infty}$ estimate, we have
    \begin{equation*}
        \begin{aligned}
        & \int_{\Gamma}T(\nabla_\Gamma u_1^*, \nabla_\Gamma u_2, g(v_1^*, \nabla_{\Gamma} v_1^*, \ldots, v_n^*, \nabla_{\Gamma} v_n^*)) \\
        & \leq C \left\|\nabla_\Gamma u_1^*\right\|_{H^{1/2}(\Gamma)} \left\|\nabla_\Gamma u_2\right\|_{H^{-1/2}(\Gamma)} \left\|g(v_1^*, \nabla_{\Gamma} v_1^*, \ldots, v_n^*, \nabla_{\Gamma} v_n^*)\right\|_{W^{1/2, \infty}(\Gamma)} \\
        & \leq C \left\|u_1^*\right\|_{H^{3/2}(\Gamma)} \left\|u_2\right\|_{H^{1/2}(\Gamma)}  \\
        & \leq C \left\|u_1\right\|_{H_h^{3/2}(\Gamma)} \left\|u_2\right\|_{H_h^{1/2}(\Gamma)} .
        \end{aligned}
    \end{equation*}
    Here we use the fact that $\left\|\nabla_\Gamma u_2\right\|_{H^{-1/2}(\Gamma)} \leq C \left\|u_2\right\|_{H^{1/2}(\Gamma)}$ together with the estimate $\left\|u_2\right\|_{H^{1/2}(\Gamma)} \leq C \left\|u_2\right\|_{H_h^{1/2}(\Gamma)}$ by Lemma \ref{lem:equivalence of norms, continuous and discrete}. For \eqref{eq:comparison of interpolated surface and exact surface 2, case 1} and \eqref{eq:comparison of interpolated surface and exact surface 2, case 3}, we can simply use Lemma \ref{lem:equivalence of norms, continuous and discrete, negative order} to get $\left\|u_2\right\|_{H^{-1/2}(\Gamma)} \leq C \left\|u_2\right\|_{H_h^{-1/2}(\Gamma)}$.

    Using Lemma \ref{lemma:comparison of different functions on the same surface}, we have
    \begin{equation*}
        \begin{aligned}
        & \int_{\Gamma} T(\nabla_\Gamma u_1, \nabla_\Gamma u_2, g) - \int_{\Gamma} T(\nabla_\Gamma u_1^*, \nabla_\Gamma u_2, g(v_1^*, \nabla_{\Gamma} v_1^*, \ldots, v_n^*, \nabla_{\Gamma} v_n^*)) \\
        & \leq C \left\|\nabla_\Gamma u_1 - \nabla_\Gamma u_1^*\right\|_{L^2(\Gamma)} \left\|\nabla_\Gamma u_2\right\|_{L^2(\Gamma)} \\
        & \qquad + C \left\|\nabla_\Gamma u_1\right\|_{L^2(\Gamma)} \left\|\nabla_\Gamma u_2\right\|_{L^2(\Gamma)} \sum_{i=1}^n \left\|v_i - v_i^*\right\|_{W^{1, \infty}(\Gamma)}  \\
        & \leq  C h^{1}\left\|u_1\right\|_{H_h^{2}(\Gamma)} h^{-1/2} \left\|u_2\right\|_{H_h^{1/2}(\Gamma)} + C \left\|u_1\right\|_{H_h^{3/2}(\Gamma)} h^{-1/2} \left\|u_2\right\|_{H_h^{1/2}(\Gamma)} \sum_{i=1}^n h^{\alpha_i} \\
        & \leq C \left\|u_1\right\|_{H_h^{3/2}(\Gamma)} \left\|u_2\right\|_{H_h^{1/2}(\Gamma)}.
        \end{aligned}
    \end{equation*}
    Combining the above two estimates, we have the desired result.
\end{proof}

We then show the $H_h^s$ norm equivalence to pull back the $H_h^s$ norm on the smooth surface $\Gamma$ to its interpolated surface $\Gamma_h[\mathbf{x}^*]$.
\begin{lemma}[$H_h^s$ norm equivalence]\label{lemma:norm equivalence for Hhs on different surfaces}
    Let $\Gamma$ be a regular surface and $\Gamma_h[\mathbf{x}^*]$ its discrete approximation of degree $k \geq 3$. Then for $s = -1/2, 1/2, 1, 3/2, 2$, there exist constants $C_1, C_2 > 0$ such that
    \begin{equation}
        C_1 \left\|u_h\right\|_{H_h^s(\Gamma_h[\mathbf{x}^*])} \leq \left\|u_h^\ell\right\|_{H_h^s(\Gamma)} \leq C_2 \left\|u_h\right\|_{H_h^s(\Gamma_h[\mathbf{x}^*])}, \quad \forall u_h \in \mathscr{S}_h[\mathbf{x}^*].
    \end{equation}
\end{lemma}
\begin{proof}
    We first prove the case $s = 2$. Let $A_h = I - \Delta_{h, \Gamma}: \mathscr{S}_h \to \mathscr{S}_h$, and we denote $B_h: \mathscr{S}_h \to \mathscr{S}_h$ as 
    \begin{equation*}
        B_h u_h^\ell = \left((I-\Delta_{h, \Gamma_h[\mathbf{x}^*]}) u_h\right)^\ell, \quad \forall u_h^\ell \in \mathscr{S}_h.
    \end{equation*}
    It is easy to see that $A_h, B_h$ are linear and self-adjoint. We first show that
    \begin{equation}
        \left\|(A_h - B_h) u_h^\ell\right\|_{L^2(\Gamma)} \leq C h^{k-1} \left\|u_h^\ell\right\|_{L^2(\Gamma)}.
    \end{equation}
    The difference between $A_h$ and $B_h$ can be estimated as
    \begin{equation*}
        \begin{aligned}
        &\int_{\Gamma} u_h^\ell (A_h - B_h) w_h^\ell\\
         & = \left(\int_{\Gamma_h[\mathbf{x}^*]} u_h (I - \Delta_{h, \Gamma_h[\mathbf{x}^*]}) w_h - \int_{\Gamma} u_h^\ell B_h w_h^\ell \right) \\
        & \quad + \left(\int_{\Gamma} \nabla_\Gamma u_h^\ell \cdot \nabla_\Gamma w_h^\ell + u_h^\ell w_h^\ell - \int_{\Gamma_h[\mathbf{x}^*]} \nabla_{\Gamma_h[\mathbf{x}^*]} u_h \cdot \nabla_{\Gamma_h[\mathbf{x}^*]} w_h + u_h w_h \right) \\
        & \leq C h^{k+1} \left(\left\|u_h^\ell\right\|_{L^2(\Gamma)} \left\|(I - \Delta_{h, \Gamma_h[\mathbf{x}^*]}) w_h\right\|_{L^2(\Gamma_h[\mathbf{x}^*])} + \left\|u_h^\ell\right\|_{H^1(\Gamma)} \left\|w_h\right\|_{H^1(\Gamma)}\right)\\
        & \leq C h^{k-1} \left\|u_h^\ell\right\|_{L^2(\Gamma)} \left\|w_h^\ell\right\|_{L^2(\Gamma)}, \quad \forall u_h^\ell, w_h^\ell \in \mathscr{S}_h.
        \end{aligned}
    \end{equation*}
    The first inequality comes from Lemma 5.6 in \cite{highorderESFEM}, the second is inverse estimate.
    
    Since $\left\|u_h^\ell\right\|_{H_h^2(\Gamma)} = \left\|A_h u_h^\ell\right\|_{L^2(\Gamma)}$, and $\left\|u_h\right\|_{H_h^2(\Gamma_h[\mathbf{x}^*])} = \left\|(B_h u_h^\ell)^{-\ell}\right\|_{L^2(\Gamma_h[\mathbf{x}^*])}$, we have the following norm equivalence:
    \begin{subequations}
        \begin{align*}
            & \left\|u_h\right\|_{H_h^2(\Gamma_h[\mathbf{x}^*])} \leq C \left\|A_h u_h^\ell\right\|_{L^2(\Gamma)} +  C \left\|(A_h - B_h)u_h^\ell\right\|_{L^2(\Gamma)} \leq C \left\|u_h^\ell\right\|_{H_h^2(\Gamma)}, \\
            & \left\|u_h^\ell\right\|_{H_h^2(\Gamma)} \leq \left\|(A_h - B_h)u_h^\ell\right\|_{L^2(\Gamma)} + \left\|B_h u_h^\ell\right\|_{L^2(\Gamma)} \leq C \left\|u_h\right\|_{H_h^2(\Gamma_h[\mathbf{x}^*])}.
        \end{align*}
    \end{subequations}
    
    For $s = 1$, we have the following identity for $A_h^{1/2} - B_h^{1/2}$:
    \begin{equation*}
        A_h^{1/2}(A_h^{1/2} - B_h^{1/2}) + (A_h^{1/2} - B_h^{1/2})(B_h^{1/2}) = A_h - B_h.
    \end{equation*}
    This is a Sylvester equation. By \eqref{eq:explicit representation of solution to Sylvester equation}, we have 
    \begin{equation*}
        \left\|\left(A_h^{1/2} - B_h^{1/2}\right) u_h^\ell\right\|_{L^2(\Gamma)}  = \left\|\int_0^\infty e^{-tA_h^{1/2}} (A_h - B_h) e^{-tB_h^{1/2}} u_h^\ell dt\right\|_{L^2(\Gamma)}, \forall u_h^\ell \in \mathscr{S}_h.
    \end{equation*}
    We know that $\left\|A_h u_h^\ell\right\|_{L^2(\Gamma)} \geq \left\|u_h^\ell\right\|_{L^2(\Gamma)}$, and $\left\|B_h u_h^\ell\right\|_{L^2(\Gamma)} \geq \left\|A_h u_h^\ell\right\|_{L^2(\Gamma)} - \left\|(A_h - B_h)u_h^\ell\right\|_{L^2(\Gamma)} \geq \frac{1}{2} \left\|u_h^\ell\right\|_{L^2(\Gamma)}$. Thus, we have
    \begin{equation}
        \left\|\left(A_h^{1/2} - B_h^{1/2}\right) u_h^\ell\right\|_{L^2(\Gamma)} \leq \int_0^\infty C e^{-t} h^{k-1} e^{-\frac{1}{\sqrt{2}}t} \left\|u_h^\ell\right\|_{L^2(\Gamma)} dt \leq C h^{k-1} \left\|u_h^\ell\right\|_{L^2(\Gamma)}.
    \end{equation}
    Similarly, we have
    \begin{equation}
        \left\|\left(A_h^{1/4} - B_h^{1/4}\right) u_h^\ell\right\|_{L^2(\Gamma)} \leq C h^{k-1} \left\|u_h^\ell\right\|_{L^2(\Gamma)}.
    \end{equation}
    The identity $A_h^{3/4} - B_h^{3/4} = A_h^{1/2} (A_h^{1/4} - B_h^{1/4}) + (A_h^{1/4} - B_h^{1/4})(B_h^{1/4})$ yields
    \begin{equation}
        \begin{aligned}
        &\left\|\left(A_h^{3/4} - B_h^{3/4}\right) u_h^\ell\right\|_{L^2(\Gamma)}\\
        &  \leq C \left\|(A_h^{1/2} - B_h^{1/2})u_h^\ell\right\|_{H_h^1(\Gamma)} + C h^{k-1}\left\|B_h^{1/4}u_h^\ell\right\|_{L^2(\Gamma)} \\
        & \leq C h^{-1}\left\|(A_h^{1/2} - B_h^{1/2})u_h^\ell\right\|_{L^2(\Gamma)} + C h^{k-1}\left\|u_h\right\|_{H_h^{1/2}(\Gamma_h[\mathbf{x}^*])} \\
        & \leq C h^{k-2} \left\|u_h^\ell\right\|_{L^2(\Gamma)}.
        \end{aligned}
    \end{equation}
    Using the same argument, we get the norm equivalence for $s = 1, 1/2, 3/2$.

    Finally, for $s = -1/2$, by definition we have
    \begin{equation*}
        \begin{aligned}
        \left\|u_h^\ell\right\|_{H_h^{-1/2}(\Gamma)} &\leq \sup_{0\neq w_h^\ell \in \mathscr{S}_h} \frac{\int_{\Gamma_h[\mathbf{x}^*]} u_h w_h}{\left\|w_h\right\|_{H_h^{1/2}(\Gamma_h[\mathbf{x}^*])}} \frac{\left\|w_h\right\|_{H_h^{1/2}(\Gamma_h[\mathbf{x}^*])}}{\left\|w_h^\ell\right\|_{H_h^{1/2}(\Gamma)}} \\
        & \qquad + \sup_{0\neq w_h^\ell \in \mathscr{S}_h} \frac{\int_{\Gamma} u_h^\ell w_h^\ell - \int_{\Gamma_h[\mathbf{x}^*]} u_h^\ell w_h^\ell}{\left\|w_h^\ell\right\|_{H_h^{1/2}(\Gamma)}} \\
        & \leq C \left\|u_h\right\|_{H_h^{-1/2}(\Gamma_h[\mathbf{x}^*])} + C h^{k+1} \sup_{0\neq w_h^\ell \in \mathscr{S}_h} \frac{\left\|u_h^\ell\right\|_{L^2(\Gamma)}\left\|w_h^\ell\right\|_{L^2(\Gamma)}}{\left\|w_h^\ell\right\|_{H_h^{1/2}(\Gamma)}} \\
        & \leq C \left\|u_h\right\|_{H_h^{-1/2}(\Gamma_h[\mathbf{x}^*])}.
        \end{aligned}
    \end{equation*}
    The last inequality is inverse estimate $\left\|u_h\right\|_{L^2(\Gamma_h[\mathbf{x}^*])} \leq C h^{-1/2} \left\|u_h\right\|_{H_h^{-1/2}(\Gamma_h[\mathbf{x}^*])}$. The other direction is similar, therefore we obtain the norm equivalence.
\end{proof}

\begin{remark}
    As a consequence of this lemma, we know that for any $u_h \in \mathscr{S}_h[\mathbf{x}^*]$, there exist constants $C_1, C_2 > 0$ such that
    \begin{equation}
        C_1 \left\|u_h \right\|_{H_h^{1/2}(\Gamma_h[\mathbf{x}^*])} \leq \left\|u_h\right\|_{H_h^{1/2}(\Gamma_h[\mathbf{x}^*])} \leq C_2 \left\|u_h \right\|_{H_h^{1/2}(\Gamma_h[\mathbf{x}^*])}.
    \end{equation}
    It combines three norm equivalences: $H^{1/2}$ norm equivalence between $\Gamma$ and $\Gamma_h[\mathbf{x}^*]$ (standard $H^1, L^2$ norm equivalence with interpolation), the $H^{1/2}, H_h^{1/2}$ equivalence on $\Gamma$ (Lemma \ref{lem:equivalence of norms, continuous and discrete}), and the $H_h^{1/2}$ norm equivalence between surfaces (Lemma \ref{lemma:norm equivalence for Hhs on different surfaces}).
\end{remark}

Combining the above lemmas, we arrive at the desired $H^{1/2}$-$H^{-1/2}$-$W^{1/2, \infty}$ estimate for a discrete surface $\Gamma_h[\mathbf{x}]$.
\begin{theorem}\label{theorem:estimate of the multilinear form}
    Let $\Gamma$ be a regular surface, $\Gamma_h[\mathbf{x}^*]$ its discrete approximation of degree $k\geq 3$, and $\Gamma_h[\mathbf{x}]$ a discrete surface with geometric perturbation satisfying $\left\|e_x\right\|_{W^{1,\infty}(\Gamma_h[\mathbf{x}^*])} \leq \frac{1}{4}$ and $\left\|e_x\right\|_{W^{1,\infty}(\Gamma_h[\mathbf{x}])} \leq C h^{\alpha_0}$ for some $C > 0, \alpha_0 > 1/2$. Consider functions $u_1, u_2 \in \mathscr{S}_h[\mathbf{x}]$, and $v_1, \ldots, v_n \in \mathscr{S}_h[\mathbf{x}]$ together with their perturbations $v_1^*, \ldots, v_n^* \in W^{3/2, \infty}(\Gamma)$. Assume there are constants $C > 0, \alpha_i>1/2$ such that $\left\|v_i\right\|_{W^{1,\infty}(\Gamma_h[\mathbf{x}])} \leq C$, $\left\|v_i^*\right\|_{W^{3/2,\infty}(\Gamma)} \leq C$, and $\left\|\widehat{v}_i^\ell - v_i^*\right\|_{W^{1, \infty(\Gamma)}} \leq C h^{\alpha_i}$, for all $i = 1, \ldots, n$. 
    
    Then for any multilinear form $T$ and any smooth function $g$ defined as in Lemma \ref{lemma:comparison of discrete surface and interpolated surface}, we have
    \begin{subequations}\label{eq:estimate of the multilinear form}
        \begin{align}
        & \int_{\Gamma_h[\mathbf{x}]} T(u_1, u_2, g) \leq C \left\|\widehat{u}_1\right\|_{H_h^{1/2}(\Gamma_h[\mathbf{x}^*])} \left\|\widehat{u}_2\right\|_{H_h^{-1/2}(\Gamma_h[\mathbf{x}^*])}, \label{eq:estimate of the multilinear form, case 1} \\
        & \int_{\Gamma_h[\mathbf{x}]} T(\nabla_{\Gamma_h[\mathbf{x}]} u_1, u_2, g) \leq C \left\|\widehat{u}_1\right\|_{H_h^{3/2}(\Gamma_h[\mathbf{x}^*])} \left\|\widehat{u}_2\right\|_{H_h^{-1/2}(\Gamma_h[\mathbf{x}^*])},\label{eq:estimate of the multilinear form, case 2} \\
        & \int_{\Gamma_h[\mathbf{x}]} T(u_1, \nabla_{\Gamma_h[\mathbf{x}]} u_2, g) \leq C \left\|\widehat{u}_1\right\|_{H_h^{1/2}(\Gamma_h[\mathbf{x}^*])} \left\|\widehat{u}_2\right\|_{H_h^{1/2}(\Gamma_h[\mathbf{x}^*])},\label{eq:estimate of the multilinear form, case 3} \\
        & \int_{\Gamma_h[\mathbf{x}]} T(\nabla_{\Gamma_h[\mathbf{x}]} u_1, \nabla_{\Gamma_h[\mathbf{x}]} u_2, g) \leq C \left\|\widehat{u}_1\right\|_{H_h^{3/2}(\Gamma_h[\mathbf{x}^*])} \left\|\widehat{u}_2\right\|_{H_h^{1/2}(\Gamma_h[\mathbf{x}^*])}. \label{eq:estimate of the multilinear form, case 4}
        \end{align}
    \end{subequations}
\end{theorem}
\begin{proof}
    We only prove \eqref{eq:estimate of the multilinear form, case 1}. By applying Lemma \ref{lemma:comparison of discrete surface and interpolated surface}, Lemma \ref{lemma:comparison of interpolated surface and exact surface}, Lemma \ref{lemma:comparison of interpolated surface and exact surface 2}, and Lemma \ref{lemma:norm equivalence for Hhs on different surfaces}, together with the inverse estimate \eqref{eq:inverse estimate}, we have
    \begin{equation*}
        \begin{aligned}
        \int_{\Gamma_h[\mathbf{x}]} T(u_1, u_2, g) & = \left(\int_{\Gamma_h[\mathbf{x}]} T(u_1, u_2, g) - \int_{\Gamma_h[\mathbf{x}^*]} T(\widehat{u}_1, \widehat{u}_2, \widehat{g})\right) \\
        &\qquad  + \left(\int_{\Gamma_h[\mathbf{x}^*]} T(\widehat{u}_1, \widehat{u}_2, \widehat{g}) - \int_{\Gamma} T(\widehat{u}_1^\ell, \widehat{u}_2^\ell, \widehat{g})\right) + \int_{\Gamma} T(\widehat{u}_1^\ell, \widehat{u}_2^\ell, \widehat{g}) \\
        & \leq C \left(\left\|e_x\right\|_{W^{1,\infty}(\Gamma_h[\mathbf{x}^*])} + h^k \right) \left\|\widehat{u}_1\right\|_{L^2(\Gamma_h[\mathbf{x}^*])} \left\|\widehat{u}_2\right\|_{L^2(\Gamma_h[\mathbf{x}^*])} \\
        & \qquad + C \left\|\widehat{u}_1^\ell\right\|_{H_h^{1/2}(\Gamma)} \left\|\widehat{u}_2^\ell\right\|_{H_h^{-1/2}(\Gamma)} \\
        & \leq C \left(h^{\alpha_0 - 1/2} + h^{k-1/2} + 1\right) \left\|\widehat{u}_1\right\|_{H_h^{1/2}(\Gamma_h[\mathbf{x}^*])} \left\|\widehat{u}_2\right\|_{H_h^{-1/2}(\Gamma_h[\mathbf{x}^*])}.
        \end{aligned}
    \end{equation*}
    The other inequalities can be proved similarly.
\end{proof}

\subsection{Time derivative of the discrete Laplacian}
In stability proofs, it is typically necessary to redistribute temporal derivatives to other functions by utilizing the fundamental identity $\partial^\bullet u_h v_h = \partial^\bullet (u_h v_h) - u_h \partial^\bullet v_h$. This technique is commonly referred to as the derivative trick. The extension of the derivative trick to fractional frameworks requires the estimation of $\partial^\bullet (I - \Delta_h)^s$. For our practical applications, we focus on the critical cases $s = 1$, $s = 1/2$, and $s = -1/2$ in the following lemmas.

The case $s = 1$ follows immediately from the standard derivative trick combined with Theorem \ref{theorem:estimate of the multilinear form}.
\begin{lemma}\label{lemma:time derivative of discrete Laplacian of order 1}
    Let $\Gamma$ be a regular surface with velocity $v \in W^{2, \infty}(\Gamma)$, and $\Gamma_h[\mathbf{x}^*]$ be its discrete approximation of degree $k \geq 3$ with a discrete velocity $v_h \in \mathscr{S}_h$. Suppose there exist constants $C>0, \alpha_0>1$ such that
    \begin{equation}
        \left\|v_h^\ell - v\right\|_{W^{1, \infty}(\Gamma)} \leq C h^{\alpha_0}.
    \end{equation}
    \begin{equation}
        \left\|v\right\|_{W^{2, \infty}(\Gamma)} \leq C, \quad \left\|v_h\right\|_{W^{1, \infty}(\Gamma_h[\mathbf{x}^*])} \leq C.
    \end{equation}
    Then for any $u_h, w_h \in \mathscr{S}_h$, we have
    \begin{subequations}
        \begin{align}
            & \int_{\Gamma_h[\mathbf{x}^*]} u_h (\partial^\bullet (I - \Delta_h)) w_h \leq C \left\|u_h\right\|_{H_h^1(\Gamma_h[\mathbf{x}^*])} \left\|w_h\right\|_{H_h^1(\Gamma_h[\mathbf{x}^*])}, \label{eq: time derivative of discrete Laplacian of order 1, 1} \\
            & \int_{\Gamma_h[\mathbf{x}^*]} u_h (\partial^\bullet (I - \Delta_h)) w_h \leq C \left\|u_h\right\|_{H_h^{3/2}(\Gamma_h[\mathbf{x}^*])} \left\|w_h\right\|_{H_h^{1/2}(\Gamma_h[\mathbf{x}^*])}, \label{eq: time derivative of discrete Laplacian of order 1, 2} \\
            & \int_{\Gamma_h[\mathbf{x}^*]} u_h (\partial^\bullet (I - \Delta_h)) w_h \leq C \left\|u_h\right\|_{H_h^{2}(\Gamma_h[\mathbf{x}^*])} \left\|w_h\right\|_{L^2(\Gamma_h[\mathbf{x}^*])}. \label{eq: time derivative of discrete Laplacian of order 1, 3}
        \end{align}
    \end{subequations}
\end{lemma}
\begin{proof}
    We only prove the second inequality \eqref{eq: time derivative of discrete Laplacian of order 1, 2}. The other inequalities can be proved similarly. Since $I - \Delta_h$ is self-adjoint, using the Leibniz formula and Lemma 4.1 in \cite{KLLP2017}, we have
    \begin{equation}\label{eq: time derivative of discrete Laplacian of order 3/2, 1, leibniz rule}
        \begin{aligned}
        &\int_{\Gamma_h[\mathbf{x}^*]} u_h (\partial^\bullet (I - \Delta_h)) w_h \\
        & = - \int_{\Gamma_h[\mathbf{x}^*]} \left((I-\Delta_h)u_h\right) w_h \nabla_{\Gamma_h[\mathbf{x}^*]} \cdot v_h + \frac{d}{dt} \int_{\Gamma_h[\mathbf{x}^*]} u_h (I - \Delta_h) w_h  \\
         & = - \int_{\Gamma_h[\mathbf{x}^*]} \left((I-\Delta_h)u_h\right) w_h \nabla_{\Gamma_h[\mathbf{x}^*]} \cdot v_h \\
        & \qquad + \int_{\Gamma_h[\mathbf{x}^*]} \nabla_{\Gamma_h[\mathbf{x}^*]} u_h \cdot \left(D_{\Gamma_h[\mathbf{x}^*]}(v_h) \nabla_{\Gamma_h[\mathbf{x}^*]} w_h\right) + \int_{\Gamma_h[\mathbf{x}^*]}  u_h w_h \nabla_{\Gamma_h[\mathbf{x}^*]} \cdot v_h
        \end{aligned}
    \end{equation}
    where $D_{\Gamma_h[\mathbf{x}^*]}(v_h) = (\nabla_{\Gamma_h[\mathbf{x}^*]} \cdot v_h) I_3 - \nabla_{\Gamma_h[\mathbf{x}^*]} v_h - (\nabla_{\Gamma_h[\mathbf{x}^*]} v_h)^T$.

    Using \eqref{eq:estimate of the multilinear form, case 1} in Theorem \ref{theorem:estimate of the multilinear form} with $u_1 = w_h, u_2 = (I- \Delta_h)u_h$ and $g = \text{tr}(\nabla_{\Gamma_h[\mathbf{x}^*]} v_h) =\nabla_{\Gamma_h[\mathbf{x}^*]} \cdot v_h$, the first term in \eqref{eq: time derivative of discrete Laplacian of order 3/2, 1, leibniz rule} can be estimated as
    \begin{equation*}
        \begin{aligned}
        &- \int_{\Gamma_h[\mathbf{x}^*]} \left((I-\Delta_h)u_h\right) w_h \nabla_{\Gamma_h[\mathbf{x}^*]} \cdot v_h\\
        & \leq C \left\|w_h\right\|_{H_h^{1/2}(\Gamma_h[\mathbf{x}^*])} \left\|(I-\Delta_h)u_h\right\|_{H^{-1/2}_h(\Gamma_h[\mathbf{x}^*])}\\
        & \leq C \left\|u_h\right\|_{H_h^{3/2}(\Gamma_h[\mathbf{x}^*])} \left\|w_h\right\|_{H_h^{1/2}(\Gamma_h[\mathbf{x}^*])}
        \end{aligned}
    \end{equation*}
    For the second term, Using \eqref{eq:estimate of the multilinear form, case 4} in Theorem \ref{theorem:estimate of the multilinear form} with $u_1 = u_h, u_2 = w_h$ and $g = D_{\Gamma_h[\mathbf{x}^*]}(v_h)$, we have
    \begin{equation*}
        \begin{aligned}
            &\int_{\Gamma_h[\mathbf{x}^*]} \nabla_{\Gamma_h[\mathbf{x}^*]} u_h \cdot \left(D_{\Gamma_h[\mathbf{x}^*]}(v_h) \nabla_{\Gamma_h[\mathbf{x}^*]} w_h\right)  \\
            &\leq C \left\|u_h\right\|_{H_h^{3/2}(\Gamma_h[\mathbf{x}^*])} \left\|w_h\right\|_{H_h^{1/2}(\Gamma_h[\mathbf{x}^*])}.
        \end{aligned}
    \end{equation*}
    The last term can be estimated similarly. 
\end{proof}

For $s = 1/2$, the analysis depends upon leveraging the Sylvester equation structure of $\partial^\bullet (I - \Delta_h)^{1/2}$.
\begin{lemma}\label{lemma:time derivative of discrete Laplacian of order 1/2}
    Under the same assumptions as in Lemma \ref{lemma:time derivative of discrete Laplacian of order 1}, for any $u_h, w_h \in \mathscr{S}_h$, we have
    \begin{subequations}
        \begin{align}
            & \int_{\Gamma_h[\mathbf{x}^*]} u_h (\partial^\bullet (I - \Delta_h)^{1/2}) w_h \leq C \left\|u_h\right\|_{H_h^{1/2}(\Gamma_h[\mathbf{x}^*])} \left\|w_h\right\|_{H_h^{1/2}(\Gamma_h[\mathbf{x}^*])}, \label{eq: time derivative of discrete Laplacian of order 1/2, 1} \\
            & \int_{\Gamma_h[\mathbf{x}^*]} u_h (\partial^\bullet (I - \Delta_h)^{1/2}) w_h \leq C \left\|u_h\right\|_{H_h^{1}(\Gamma_h[\mathbf{x}^*])} \left\|w_h\right\|_{L^2(\Gamma_h[\mathbf{x}^*])}, \label{eq: time derivative of discrete Laplacian of order 1/2, 2} \\
            & \int_{\Gamma_h[\mathbf{x}^*]} u_h (\partial^\bullet (I - \Delta_h)^{1/2}) w_h \leq C \left\|u_h\right\|_{H_h^{3/2}(\Gamma_h[\mathbf{x}^*])} \left\|w_h\right\|_{H_h^{-1/2}(\Gamma_h[\mathbf{x}^*])}. \label{eq: time derivative of discrete Laplacian of order 1/2, 3}
        \end{align}
    \end{subequations}
\end{lemma}
\begin{proof}
    We only prove the second inequality \eqref{eq: time derivative of discrete Laplacian of order 1/2, 2}. The other inequalities can be proved similarly. Taking the time derivative of $\left(I - \Delta_h\right)^{1/2}\left(I - \Delta_h\right)^{1/2} = I - \Delta_h$ on both sides, we have
    \begin{equation*}
        \left(\partial^\bullet \left(I - \Delta_h\right)^{1/2}\right) \left(I - \Delta_h\right)^{1/2} + \left(I - \Delta_h\right)^{1/2} \left(\partial^\bullet \left(I - \Delta_h\right)^{1/2}\right) = \partial^\bullet (I - \Delta_h).
    \end{equation*}
    This is a Sylvester equation. Since $\left(I - \Delta_h\right)^{1/2}$ is positive and self-adjoint, by \eqref{eq:explicit representation of solution to Sylvester equation}, we obtain
    \begin{equation}
        \begin{aligned}
        &\int_{\Gamma_h[\mathbf{x}^*]} u_h (\partial^\bullet (I - \Delta_h)^{1/2}) w_h \\
        & =  \int_0^\infty \int_{\Gamma_h[\mathbf{x}^*]} u_h e^{-t (I - \Delta_h)^{1/2}} \left(\partial^\bullet (I - \Delta_h)\right) e^{-t (I - \Delta_h)^{1/2}} w_h dt.
        \end{aligned}
    \end{equation}
    Next, using the fact that $I - \Delta_h$ is positive and self-adjoint, we have
    \begin{equation}
        \begin{aligned}
        &\int_0^\infty \left\|e^{-t (I - \Delta_h)^{1/2}}u_h\right\|_{H_h^s(\Gamma_h[\mathbf{x}^*])}^2 dt \\
        & = \int_{\Gamma_h[\mathbf{x}^*]} \int_0^\infty \left(e^{-t (I - \Delta_h)^{1/2}}u_h\right) \left(I - \Delta_h\right)^{s} \left(e^{-t (I - \Delta_h)^{1/2}}u_h\right) dt \\
        & = \int_{\Gamma_h[\mathbf{x}^*]} u_h \left( \int_0^\infty e^{-2t(I - \Delta_h)^{1/2}} dt \right) (I - \Delta_h)^s u_h  \\
        & = \int_{\Gamma_h[\mathbf{x}^*]} u_h (2(I - \Delta_h)^{1/2})^{-1} (I - \Delta_h)^s u_h = \frac{1}{2} \left\|u_h\right\|_{H_h^{s - 1/2}(\Gamma_h[\mathbf{x}^*])}^2.
        \end{aligned}
    \end{equation}
    Combining these two equations, and using \eqref{eq: time derivative of discrete Laplacian of order 1, 2} in Lemma \ref{lemma:time derivative of discrete Laplacian of order 1} yield that
    \begin{equation}
        \begin{aligned}
            &\int_{\Gamma_h[\mathbf{x}^*]} u_h (\partial^\bullet (I - \Delta_h)^{1/2}) w_h \\
            &  \leq C \int_0^\infty \left\|e^{-t (I - \Delta_h)^{1/2}}u_h\right\|_{H_h^{3/2}(\Gamma_h[\mathbf{x}^*])} \left\|e^{-t (I - \Delta_h)^{1/2}}w_h\right\|_{H_h^{1/2}(\Gamma_h[\mathbf{x}^*])} dt \\
            & \leq C \sqrt{\int_0^\infty \left\|e^{-t (I - \Delta_h)^{1/2}}u_h\right\|_{H_h^{3/2}(\Gamma_h[\mathbf{x}^*])}^2 dt} \sqrt{\int_0^\infty \left\|e^{-t (I - \Delta_h)^{1/2}}w_h\right\|_{H_h^{1/2}(\Gamma_h[\mathbf{x}^*])}^2 dt} \\
            & \leq C \left\|u_h\right\|_{H_h^{1}(\Gamma_h[\mathbf{x}^*])} \left\|w_h\right\|_{L^2(\Gamma_h[\mathbf{x}^*])}.
        \end{aligned}
    \end{equation}
\end{proof}

The case $s = -1/2$ follows from the identity $\left(I - \Delta_h\right)^{-1/2}\left(I - \Delta_h\right)^{1/2} = I$ combined with the established result for $s = 1/2$ in Lemma \ref{lemma:time derivative of discrete Laplacian of order 1/2}.
\begin{lemma}\label{lemma:time derivative of discrete Laplacian of order -1/2}
    Under the same assumptions as in Lemma \ref{lemma:time derivative of discrete Laplacian of order 1}, for any $u_h, w_h \in \mathscr{S}_h$, we have
    \begin{subequations}
        \begin{align}
            & \int_{\Gamma_h[\mathbf{x}^*]} u_h (\partial^\bullet (I - \Delta_h)^{-1/2}) w_h \leq C \left\|u_h\right\|_{H_h^{-1/2}(\Gamma_h[\mathbf{x}^*])} \left\|w_h\right\|_{H_h^{-1/2}(\Gamma_h[\mathbf{x}^*])}, \label{eq: time derivative of discrete Laplacian of order -1/2, 1} \\
            & \int_{\Gamma_h[\mathbf{x}^*]} u_h (\partial^\bullet (I - \Delta_h)^{-1/2}) w_h \leq C \left\|u_h\right\|_{H_h^{-1}(\Gamma_h[\mathbf{x}^*])} \left\|w_h\right\|_{L^2(\Gamma_h[\mathbf{x}^*])}. \label{eq: time derivative of discrete Laplacian of order -1/2, 2} 
        \end{align}
    \end{subequations}
\end{lemma}
\begin{proof}
    We only prove the second inequality \eqref{eq: time derivative of discrete Laplacian of order -1/2, 2}.  Taking the time derivative of $\left(I - \Delta_h\right)^{1/2}\left(I - \Delta_h\right)^{-1/2} = I$ on both sides, we have
    \begin{equation*}
        \left(\partial^\bullet \left(I - \Delta_h\right)^{1/2}\right) \left(I - \Delta_h\right)^{-1/2} + \left(I - \Delta_h\right)^{1/2} \left(\partial^\bullet \left(I - \Delta_h\right)^{-1/2}\right) = 0.
    \end{equation*}
    Therefore, the time derivative of $\left(I - \Delta_h\right)^{-1/2}$ is given by
    \begin{equation}\label{eq: time derivative of discrete Laplacian of order -1/2, representation}
        \partial^\bullet \left(I - \Delta_h\right)^{-1/2} = - \left(I - \Delta_h\right)^{-1/2} \left(\partial^\bullet \left(I - \Delta_h\right)^{1/2}\right) \left(I - \Delta_h\right)^{-1/2}.
    \end{equation}
    Using this and \eqref{eq: time derivative of discrete Laplacian of order 1/2, 2} in Lemma \ref{lemma:time derivative of discrete Laplacian of order 1/2}, we have
    \begin{equation*}
        \begin{aligned}
            &\int_{\Gamma_h[\mathbf{x}^*]} u_h (\partial^\bullet (I - \Delta_h)^{-1/2}) w_h \\
            & = \int_{\Gamma_h[\mathbf{x}^*]} u_h\left(I - \Delta_h\right)^{-1/2} \left(\partial^\bullet \left(I - \Delta_h\right)^{1/2}\right) \left(I - \Delta_h\right)^{-1/2} w_h \\
            & \leq C \left\|\left(I - \Delta_h\right)^{-1/2} u_h\right\|_{L^2(\Gamma_h[\mathbf{x}^*])} \left\|\left(I - \Delta_h\right)^{-1/2} w_h\right\|_{H^1(\Gamma_h[\mathbf{x}^*])} \\
            & \leq C \left\|u_h\right\|_{H_h^{-1}(\Gamma_h[\mathbf{x}^*])} \left\|w_h\right\|_{L^2(\Gamma_h[\mathbf{x}^*])}.
        \end{aligned}
    \end{equation*}
\end{proof}

\section{Proof of the main result}
In this section, we give a detailed proof of the main result Theorem \ref{thm:ESFEM_error_bound}. For simplicity, we restrict to $\alpha = \beta = 1$. Throughout this section, we also adopt the notation $\dot{f}$ for the discrete material derivative $\partial_h^{\bullet} f$.
\subsection{Error equations and consistency errors}
The notations and results in this subsection are adopted from \cite[Section~5]{EKL24}.

Let $(X, v, u, w)$ denote the exact solutions with $C^{k+1}\cap W^{2, \infty}$ regularity, where $w=(n, H)$.  Following \cite{EKL24}, their finite element projections are $(X_h^*, v_h^*, u_h^*, w_h^*)$, where $w_h^*=(n_h^*, H_h^*)$.
\begin{alignat*}{5}
    X_h^* = &\ \tilde{I}_h^\Omega(X), &&& &\\
    v_h^* = (v_{\Omega, h}^*, v_{\Gamma, h}^*) = &\ \tilde{I}_h^\Omega(v) , \qquad & \text{therefore} & \qquad & \gamma_h^*(v_{\Omega, h}^*) = v_{\Gamma,h}^* = &\ \tilde{I}_h^\Gamma(\gamma(v)), \\
    u_h^* = &\ \tilde{R}_h^{\text{Ritz}}(u) , \qquad & \text{therefore} & \qquad & \gamma_h^*(u_h^*) = &\ \tilde{R}_h^\Gamma (\gamma(u)), \\
    w_h^* = (n_h^*,H_h^*)  = &\ (\tilde{R}_h^\Gamma n , \tilde{R}_h^\Gamma H), \quad &&&  & Q_h^* = \tilde{I}_h^\Gamma(Q),
\end{alignat*}
where the definitions of $\tilde{I}_h^\Omega$, $\tilde{I}_h^\Gamma$, $\tilde{R}_h^{\text{Ritz}}$ and $\tilde{R}_h^\Gamma$ are given in \cite{EKL24}.

Let us introduce the error vectors between the semi-discrete solutions and their corresponding finite element interpolations. For the error of flow map, we define
\begin{equation*}
    e_x:= \sum_{j=1}^N ({\bf x}_j - {\bf x}^*_j) \phi_j[{\bf x}^*] = \Phi_h - \text{id} \in \mathscr{V}_h[{\bf x}^*].
\end{equation*}
For function $u$ defined on $\Omega(t)/\Gamma(t)$, we have: 
\begin{equation*}
    e_u:= \widehat{u}_h - u_h^* \in \mathscr{V}_h[{\bf x}^*],
\end{equation*}
the errors $e_v = (e_{v_{\Omega}}, e_{v_{\Gamma}}), e_w=(e_{n}, e_H)$ are defined analogously.
\subsubsection{Error equation for the Robin boundary value problem}
Insert the interpolated solution $u_h^*$, $X_h^*$, $v_h^*$, $w_h^*$ into the Robin boundary value problem, we know that there exists a defect $d_u \in \mathscr{V}_h[{\bf x}^*]$ such that

\begin{equation}
    \label{eq:robin-bvp, weak form, interpolated solution} 
    \begin{aligned}
    \mathbf{L}(\mathbf{x}^*)\mathbf{u}^* = \mathbf{f}_u(\mathbf{x}^*, \mathbf{H}^*) + \mathbf{M}_{\bar{\Omega}}(\mathbf{x}^*)\mathbf{d}_{\mathbf{u}}.
    \end{aligned}
\end{equation}

Subtracting \eqref{eq:robin-bvp, weak form, interpolated solution} from \eqref{eq:robin-bvp, weak form}, we obtain the error $\mathbf{e}_\mathbf{u}$ satisfies
\begin{equation}\label{eq:robin-bvp, error equation}
    \begin{aligned}
        &\mathbf{L}(\mathbf{x})\mathbf{e}_\mathbf{u} = \left(\mathbf{f}_u(\mathbf{x}, \mathbf{e}_\mathbf{H} + \mathbf{H}^*) - \mathbf{f}_u(\mathbf{x}^*, \mathbf{H}^*) \right) - \left(\mathbf{L}(\mathbf{x}) - \mathbf{L}(\mathbf{x}^*)\right)\mathbf{u}^* - \mathbf{M}_{\bar{\Omega}}(\mathbf{x}^*)\mathbf{d}_{\mathbf{u}}.
    \end{aligned}
\end{equation}

Adapt the same technique as in \cite[Subsection~5.3.1]{EKL24}, we know the defect bounds for $\mathbf{d}_{\mathbf{u}}$ and its time derivative $\dot{\mathbf{d}}_{\mathbf{u}}$ are
\begin{equation*}
    \left\|\mathbf{d}_{\mathbf{u}}\right\|_{\mathbf{M}_{\bar{\Omega}}[{\bf x}^*(t)]} = \mathcal{O}\left(h^{k+1}\right), \quad \left\|\dot{\mathbf{d}}_{\mathbf{u}}\right\|_{\mathbf{M}_{\bar{\Omega}}[{\bf x}^*(t)]} = \mathcal{O}\left(h^{k+1}\right).
\end{equation*}
Here $\left\|\mathbf{d}_{\mathbf{u}}\right\|_{\mathbf{M}_{\bar{\Omega}}[{\bf x}^*(t)]} = \left\|d_u \right\|_{L^2(\Gamma_h[{\bf x}^*(t)])}$.

\subsubsection{Error equation for the forced mean curvature flow}
Insert the interpolated solution $u_h^*$, $X_h^*$, $v_h^*$, $w_h^*$ into the forced mean curvature flow, we know that there exists defects $d_H \in \mathscr{S}_h[{\bf x}^*]$ and $d_{n}, d_{v_{\Gamma}} \in [\mathscr{S}_h[{\bf x}^*]]^3$ such that 
\begin{equation}
    \label{eq:forced-mcf, semi-discretization, interpolated solution}
    \begin{aligned}
        & \mathbf{M}(\mathbf{x}^*) \dot{\mathbf{n}}^* + \mathbf{A}(\mathbf{x}^*) \mathbf{n}^* = \mathbf{f}_{n}(\mathbf{x}^*, \mathbf{n}^*) - \mathbf{D}(\mathbf{x}^*) \boldsymbol{\gamma} \mathbf{u}^* + \mathbf{M}_{\Gamma}(\mathbf{x}^*) \mathbf{d}_{\mathbf{n}}, \\[1em]
        & \mathbf{M}(\mathbf{x}^*) \dot{\mathbf{H}}^* + \mathbf{A}(\mathbf{x}^*) \mathbf{H}^* = \mathbf{f}_{H}(\mathbf{x}^*, \mathbf{n}^*, -\mathbf{H}^* + \mathbf{u}^*) + \mathbf{A}(\mathbf{x}^*) \boldsymbol{\gamma} \mathbf{u}^* + \mathbf{M}_{\Gamma}(\mathbf{x}^*) \mathbf{d}_{\mathbf{H}}, \\[1em]
        & v_{\Gamma_h}^* = \sum \left(V_{j, h}^* n_{j, h}^*\right) \psi_j[{\bf x}^*] 
        + d_{v_{\Gamma}},  \qquad  \text{with } V_h^* = - H_h^* + u_h^*.
    \end{aligned}
    \end{equation}
Here $A_h^* = \nabla_{\Gamma_h[{\bf x}^*]} n_h^*$.

Subtracting \eqref{eq:forced-mcf, semi-discretization, interpolated solution} from \eqref{eq:forced-mcf, semi-discretization, n} and \eqref{eq:forced-mcf, semi-discretization, H}, we obtain the error $\mathbf{e}_\mathbf{n}$ and $\mathbf{e}_\mathbf{H}$ satisfies
\begin{equation}\label{eq:forced-mcf, error equation for nu}
    \begin{aligned}
        \mathbf{M}(\mathbf{x}^*) \dot{\mathbf{e}}_\mathbf{n} + \mathbf{A}(\mathbf{x}^*) \mathbf{e}_\mathbf{n} & = - \left(\mathbf{M}(\mathbf{x}) - \mathbf{M}(\mathbf{x}^*)\right) \dot{\mathbf{n}}^* - \left(\mathbf{A}(\mathbf{x}) - \mathbf{A}(\mathbf{x}^*)\right) \mathbf{n}^* \\
        & \qquad - \left(\mathbf{M}(\mathbf{x}) - \mathbf{M}(\mathbf{x}^*)\right) \dot{\mathbf{e}}_\mathbf{n} - \left(\mathbf{A}(\mathbf{x}) - \mathbf{A}(\mathbf{x}^*)\right) \mathbf{e}_\mathbf{n} \\
        & \qquad + \left(\mathbf{f}_{n}(\mathbf{x}, \mathbf{n}) - \mathbf{f}_{n}(\mathbf{x}^*, \mathbf{n}^*) \right) - \mathbf{D}(\mathbf{x}) \boldsymbol{\gamma} \mathbf{e}_\mathbf{u} \\
        & \qquad - \left(\mathbf{D}(\mathbf{x}) - \mathbf{D}(\mathbf{x}^*)\right) \boldsymbol{\gamma} \mathbf{u}^* - \mathbf{M}_{\Gamma}(\mathbf{x}^*) \mathbf{d}_{\mathbf{n}}.
    \end{aligned}
\end{equation}
\begin{equation}\label{eq:forced-mcf, error equation for H}
    \begin{aligned}
        \mathbf{M}(\mathbf{x}^*) \dot{\mathbf{e}}_\mathbf{H} + \mathbf{A}(\mathbf{x}^*) \mathbf{e}_\mathbf{H} & = - \left(\mathbf{M}(\mathbf{x}) - \mathbf{M}(\mathbf{x}^*)\right) \dot{\mathbf{H}}^* - \left(\mathbf{A}(\mathbf{x}) - \mathbf{A}(\mathbf{x}^*)\right) \mathbf{H}^* \\
        & \qquad - \left(\mathbf{M}(\mathbf{x}) - \mathbf{M}(\mathbf{x}^*)\right) \dot{\mathbf{e}}_\mathbf{H} - \left(\mathbf{A}(\mathbf{x}) - \mathbf{A}(\mathbf{x}^*)\right) \mathbf{e}_\mathbf{H} \\
        & \qquad + \left(\mathbf{f}_{H}(\mathbf{x}, \mathbf{n}^*, -\mathbf{H}^* + \mathbf{u}^*) - \mathbf{f}_{H}(\mathbf{x}^*, \mathbf{n}^*, -\mathbf{H}^* + \mathbf{u}^*) \right) \\
        & \qquad + \mathbf{A}(\mathbf{x}) \boldsymbol{\gamma} \mathbf{e}_\mathbf{u} + \left(\mathbf{A}(\mathbf{x}) - \mathbf{A}(\mathbf{x}^*)\right)  \boldsymbol{\gamma} \mathbf{u}^* - \mathbf{M}_{\Gamma}(\mathbf{x}^*) \mathbf{d}_{\mathbf{H}}.
    \end{aligned}
\end{equation}
And the error $e_{v_{\Gamma}}$ is
\begin{equation}\label{eq:forced-mcf, error equation for v_Gamma}
    e_{v_{\Gamma}} = \sum \left(V_{j, h}n_{j, h} - V_{j, h}^* n_{j, h}^* \right) \psi_j[{\bf x}^*] + d_{v_{\Gamma}}.
\end{equation}

From \cite[Subsection~5.3.2]{EKL24}, we know the defect bounds for $d_{n}, d_{H}, d_{v_{\Gamma}}$ are
\begin{equation*}
    \left\| d_{n} \right\|_{L^2(\Gamma_h[{\bf x}^*])} = \mathcal{O}\left(h^{k+1}\right), \, \left\| d_{H} \right\|_{L^2(\Gamma_h[{\bf x}^*])} = \mathcal{O}\left(h^{k+1}\right), \, \left\| d_{v_{\Gamma}} \right\|_{H^1(\Gamma_h[{\bf x}^*])} = \mathcal{O}\left(h^{k}\right).
\end{equation*}
%Follow the same technique as in \cite[Lemma 8.1]{MCF}, the discrete material derivative of the defect $d_{n}, d_{H}$ are also of order $\mathcal{O}\left(h^{k+1}\right)$.
%\begin{equation*}
%    \left\| \dot{d}_{n} \right\|_{L^2(\Gamma_h[{\bf x}^*])} = \mathcal{O}\left(h^{k+1}\right), \, \left\| \dot{d}_{H} \right\|_{L^2(\Gamma_h[{\bf x}^*])} = \mathcal{O}\left(h^{k+1}\right).
%\end{equation*}
And the initial error $e_H(0), e_n(0)$ in $L^2$ are also of order $\mathcal{O}\left(h^{k+1}\right)$.

\subsubsection{Error equation for the harmonic velocity extension}
Insert the interpolated solution $X_h^*$, $v_h^*$, into the harmonic velocity extension, we know that there exists defects $d_{v_{\Omega}} \in [\mathscr{V}_h[{\bf x}^*]]^3$ such that for all $ \phi^{v}_h \in [\mathscr{V}_h[{\bf x}^*]]^3$
\begin{equation}
    \label{eq:harmonic-velocity-extension, semi-discretization, interpolated solution}
    \begin{aligned}
        \mathbf{A}_{\Omega\Omega}(\mathbf{x}^*) \mathbf{v}_{\Omega}^* = - \mathbf{A}_{\Omega\Gamma}(\mathbf{x}^*) \mathbf{v}_{\Gamma}^* + \mathbf{M}_{\Omega\Omega}(\mathbf{x}^*) \mathbf{d}_{\mathbf{v}_{\Omega}},
    \end{aligned}
\end{equation}
with $\gamma_h(e_{v_{\Omega}}) = e_{v_{\Gamma}}$.

From \cite[Subsection~5.3.3]{EKL24} and the inverse estimate, we know the defect bound for $d_{v_{\Omega}}$ is
\begin{equation*}
    \left\|\mathbf{d}_{\mathbf{v}}\right\|_{\star, {\mathbf{x}^*}} = \mathcal{O}\left(h^{k}\right).
\end{equation*}
Here $\left\|\mathbf{d}\right\|_{\star, {\mathbf{x}}} = \sup\limits_{0\neq \phi_h \in \mathscr{V}_h^0[{\bf x}]} \frac{\int_{\Omega_h[\mathbf{x}]} d \, \phi_h}{\left\|\phi_h\right\|_{H^1(\Omega_h[\mathbf{x}])}}$.

Subtracting \eqref{eq:harmonic-velocity-extension, semi-discretization, interpolated solution} from \eqref{eq:harmonic extension, semi-discretization}, we obtain the error $\mathbf{e}_{\mathbf{v}_{\Omega}}$ satisfies
\begin{equation}
    \label{eq:harmonic-velocity-extension, error equation}
    \begin{aligned}
        \mathbf{A}_{\Omega\Omega}(\mathbf{x}^*) \mathbf{e}_{\mathbf{v}_{\Omega}} + \mathbf{A}_{\Omega\Gamma}(\mathbf{x}^*) \mathbf{e}_{\mathbf{v}_{\Gamma}} & = - \left(\mathbf{A}_{\Omega\Omega}(\mathbf{x}) - \mathbf{A}_{\Omega\Omega}(\mathbf{x}^*)\right) \left(\mathbf{e}_{\mathbf{v}_{\Omega}} + \mathbf{v}_{\Omega}^*\right) \\
        & \qquad - \left(\mathbf{A}_{\Omega\Gamma}(\mathbf{x}) - \mathbf{A}_{\Omega\Gamma}(\mathbf{x}^*)\right) \left(\mathbf{e}_{\mathbf{v}_{\Gamma}} + \mathbf{v}_{\Gamma}^*\right) \\
        & \qquad - \mathbf{M}_{\Omega\Omega}(\mathbf{x}^*) \mathbf{d}_{\mathbf{v}_{\Omega}}.
    \end{aligned}
\end{equation}

\subsubsection{The $W^{1, \infty}(\Gamma)$ norm of the boundary errors}
To apply Theorem \ref{theorem:estimate of the multilinear form} and Lemma \ref{lemma:time derivative of discrete Laplacian of order 1}, we need to bound the $W^{1, \infty}(\Gamma)$ norm of the boundary errors.

For the interpolation error by $\tilde{I}_h^\Gamma$, from \cite[(5.1)]{EKL24}, it holds that
\begin{equation}
    \left\|z - (\tilde{I}_h^\Gamma z)^\ell\right\|_{H^1(\Gamma)} \leq C h^{k} \left\|z\right\|_{H^{k+1}(\Gamma)}, \quad \forall z \in H^{k+1}(\Gamma).
\end{equation}
The inverse estimate and \cite[(4.4.22)]{brenner2008mathematical} for the standard interpolation $I_h$ yield that
\begin{equation}\label{eq: tmp, third reason for k geq 3}
    \begin{aligned}
        \left\|z - (\tilde{I}_h^\Gamma z)^\ell\right\|_{W^{1, \infty}(\Gamma)} &\leq \left\|z - I_h z\right\|_{W^{1, \infty}(\Gamma)} + \left\|I_h z - (\tilde{I}_h^\Gamma z)^\ell\right\|_{W^{1, \infty}(\Gamma)} \\
        & \leq C h^{k-1} \left\|z\right\|_{H^{k+1}(\Gamma)} + C h^{-1} \left\|I_h z  - (\tilde{I}_h^\Gamma z)^\ell\right\|_{H^1(\Gamma)} \\
        & \leq C h^{k-1} \left\|z\right\|_{H^{k+1}(\Gamma)}.
    \end{aligned}
\end{equation}
Since $k-1 > 1$, we can apply Theorem \ref{theorem:estimate of the multilinear form} and Lemma \ref{lemma:time derivative of discrete Laplacian of order 1} with $x^* = x$ and $v^* = v$.

For the Ritz projection error by $\tilde{R}_h^\Gamma$, from Theorem 6.2 and Theorem 6.3 in \cite{highorderESFEM}, it holds that
\begin{equation}
    \left\|z - (\tilde{R}_h^\Gamma z)^\ell\right\|_{H^1(\Gamma)} \leq C h^{k} \left\|z\right\|_{H^{k+1}(\Gamma)}, \quad \forall z \in H^{k+1}(\Gamma).
\end{equation}
\begin{equation}
    \left\|\partial^\bullet z - \partial^\bullet(\tilde{R}_h^\Gamma z)^\ell\right\|_{H^1(\Gamma)} \leq C h^{k} \left(\left\|\partial^\bullet z\right\|_{H^{k+1}(\Gamma)} + \left\|z\right\|_{H^{k+1}(\Gamma)}\right), \, \forall z , \partial^\bullet z \in H^{k+1}(\Gamma).
\end{equation}
Using the same argument, we know that $\left\|z - (\tilde{R}_h^\Gamma z)^\ell\right\|_{W^{1, \infty}(\Gamma)} \leq C h^{k-1} \left\|z\right\|_{H^{k+1}(\Gamma)}$ and $\left\|\partial^\bullet z - \partial^\bullet(\tilde{R}_h^\Gamma z)^\ell\right\|_{W^{1, \infty}(\Gamma)} \leq C h^{k-1} \left(\left\|\partial^\bullet z\right\|_{H^{k+1}(\Gamma)} + \left\|z\right\|_{H^{k+1}(\Gamma)}\right)$. Since $k-1 > 1/2$, we can apply Theorem \ref{theorem:estimate of the multilinear form} with $u^* = u$, $(\partial^\bullet u)^* = \partial^\bullet u$, $w^* = w$, and $(\partial^\bullet w)^* = \partial^\bullet w$.

\subsection{Stability bound}
Here we establish the stability bound for the error of the semi-discretized solution. The proof follows the structure of \cite[Section~4.4]{EKL24}. The main difficulty is to derive the $H_h^{1/2}$ version of the forced mean curvature flow, which aims to ensure the compatibility of the $H^1$ bulk error and the $H_h^{1/2}$ boundary error.

\begin{proposition} 
    Assume that the exact solutions $(X, v, u, w)$ have $C^{k+1}\cap W^{2, \infty}$ regularity for $0\leq t\leq T$. Let $\delta, \varepsilon^0$ be a bound of the defects and a bound of the initial errors, given as follows:
    \begin{equation}
        \begin{aligned}
            \delta & = \max_{0\leq t\leq T} \left(\left\|\mathbf{d}_{\mathbf{u}}\right\|_{\mathbf{M}_{\bar{\Omega}}[{\bf x}^*(t)]} + \left\|\dot{\mathbf{d}}_{\mathbf{u}}(t)\right\|_{\mathbf{M}_{\bar{\Omega}}[{\bf x}^*(t)]} + \left\|\mathbf{d}_{\mathbf{v}}\right\|_{\star, {\mathbf{x}^*(t)}}\right. \\
            &  \qquad \left.+  \left\|d_H(t)\right\|_{H_h^{1/2}(\Gamma_h[\mathbf{x}^*(t)])}  + \left\|d_n(t)\right\|_{H_h^{1/2}(\Gamma_h[\mathbf{x}^*(t)])} + \left\|d_{v_{\Gamma}}(t)\right\|_{H_h^{1/2}(\Gamma_h[\mathbf{x}^*(t)])} \right). 
        \end{aligned}
    \end{equation}
    \begin{equation}
        \begin{aligned}
            \varepsilon^0 & = \left\|e_H(0)\right\|_{H_h^{-1/2}(\Gamma_h[\mathbf{x}^*(0)])} + \left\|e_H(0)\right\|_{H_h^{1/2}(\Gamma_h[\mathbf{x}^*(0)])} \\
            & \qquad + \left\|e_n(0)\right\|_{H_h^{-1/2}(\Gamma_h[\mathbf{x}^*(0)])} + \left\|e_n(0)\right\|_{H_h^{1/2}(\Gamma_h[\mathbf{x}^*(0)])}.
        \end{aligned}
    \end{equation}
    Suppose that $\delta, \varepsilon^0$ are bounded by
    \begin{equation}
        \delta + \varepsilon^0 \leq C h^{\kappa}, \qquad \text{for some } \kappa \text{ with } 2< \kappa \leq k.
    \end{equation}
    Then there exists $\bar{h}>0$ such that for all $h\leq \bar{h}$ and $t\in[0,T]$, the following stability bound holds:
    \begin{equation}
        \label{eq:stability bound}
        \begin{aligned}
            &\left\| e_u \right\|_{H^1(\Omega_h[{\bf x}^*(t)])} + \left\| e_v \right\|_{H^1(\Omega_h[{\bf x}^*(t)])} + \left\| e_x \right\|_{H^1(\Omega_h[{\bf x}^*(t)])} \\
            &\qquad + \left\| e_H \right\|_{H_h^{1/2}(\Gamma_h[{\bf x}^*(t)])} + \left\| e_n \right\|_{H_h^{1/2}(\Gamma_h[{\bf x}^*(t)])} \leq C (\delta + \varepsilon^0).
        \end{aligned}
    \end{equation}
	Here $C$ is independent of $h$ and $t$. 
\end{proposition}

\begin{proof}
    We employ a bootstrap argument. Define $t^*\in [0, T]$ as the maximal time up to which the following errors remain sufficiently small:
    \begin{equation}\label{eq:error bounds, assumptions}
        \begin{aligned}
            & \left\|e_u\right\|_{W^{1, \infty}(\Omega_h[{\bf x}^*(t)])} \leq h^{(\kappa - 2)/2 + 1/2} \\
            & \left\|e_v\right\|_{W^{1, \infty}(\Omega_h[{\bf x}^*(t)])} \leq h^{(\kappa - 2)/2 + 1/2} \\
            & \left\|e_x\right\|_{W^{1, \infty}(\Omega_h[{\bf x}^*(t)])} \leq h^{(\kappa - 2)/2 + 1/2} \\
            & \left\|\gamma_h(e_u)\right\|_{W^{1, \infty}(\Gamma_h[{\bf x}^*(t)])} \leq h^{(\kappa - 2)/2 + 1/2} \\
            & \left\|\gamma_h(e_v)\right\|_{W^{1, \infty}(\Gamma_h[{\bf x}^*(t)])} \leq h^{(\kappa - 2)/2 + 1/2} \\
            & \left\|\gamma_h(e_x)\right\|_{W^{1, \infty}(\Gamma_h[{\bf x}^*(t)])} \leq h^{(\kappa - 2)/2 + 1/2} \\
            & \left\|e_H\right\|_{W^{1, \infty}(\Gamma_h[{\bf x}^*(t)])} \leq h^{(\kappa - 2)/2 + 1/2} \\
            & \left\|e_n\right\|_{W^{1, \infty}(\Gamma_h[{\bf x}^*(t)])} \leq h^{(\kappa - 2)/2 + 1/2}
        \end{aligned} \qquad \textrm{ for } \quad t\in[0,t^*].
    \end{equation}
    We will employ a bootstrap argument to establish our result. First, we demonstrate that \eqref{eq:error bounds, assumptions} holds for $t\in[0,t^*]$, and then prove that $t^*$ extends to the full time interval, i.e., $t^*=T$, which completes the proof. $t^*>0$ is guaranteed by the inverse estimate applied to the initial conditions.
    
    For notational convenience, we omit $\gamma_h$ in boundary norms, so that $\left\|f\right\|_{H_h^{1/2}(\Gamma_h)}$ is understood as $\left\|\gamma_h(f)\right\|_{H_h^{1/2}(\Gamma_h)}$ for functions $f$ defined on the domain. 

    \textit{(A) Error bounds for the Robin boundary value problem: } We test \eqref{eq:robin-bvp, error equation} with $\mathbf{e}_{\mathbf{u}}$ and obtain
    \begin{equation}
        \begin{aligned}
            \mathbf{e}_{\mathbf{u}}^T\mathbf{L}(\mathbf{x}) \mathbf{e}_{\mathbf{u}} &= \mathbf{e}_{\mathbf{u}}^T\left(\mathbf{f}_u(\mathbf{x}, \mathbf{e}_\mathbf{H} + \mathbf{H}^*) - \mathbf{f}_u(\mathbf{x}^*, \mathbf{H}^*) \right) \\
            & \qquad  - \mathbf{e}_{\mathbf{u}}^T\left(\mathbf{L}(\mathbf{x}) - \mathbf{L}(\mathbf{x}^*)\right)\mathbf{u}^* - \mathbf{e}_{\mathbf{u}}^T\mathbf{M}_{\bar{\Omega}}(\mathbf{x}^*)\mathbf{d}_{\mathbf{u}}.
        \end{aligned}
    \end{equation}

    For the left-hand side, using the norm equivalence, we know that
    \begin{equation*}
        \begin{aligned}
            \mathbf{e}_{\mathbf{u}}^T\mathbf{L}(\mathbf{x}) \mathbf{e}_{\mathbf{u}} \geq C\left\|e_u\right\|_{H^1(\Omega_h[{\bf x}^*])}^2.
        \end{aligned}
    \end{equation*}

    For the first term on the right-hand side, it can be written as
    \begin{equation*}
        \begin{aligned}
            &-\int_0^1 d\theta \int_{\Omega_h^\theta}e_u^\theta \nabla \cdot e_x^\theta + \int_0^1 d\theta \int_{\Gamma_h^\theta} e_u^\theta \left(H_h^\theta + Q_h^\theta\right) \left(\nabla_{\Gamma_h^\theta} \cdot e_x^\theta\right) \\
            & \qquad + \int_{\Gamma_h[\mathbf{x}]} e_u \left(e_H + e_Q\right) + \int_0^1 d\theta \int_{\Gamma_h^\theta} e_u^\theta (e_H^\theta + e_Q^\theta).
        \end{aligned}
    \end{equation*}
    The first term is bounded by the $L^2$-$L^2$ estimate. For the rest terms, applying Theorem \ref{theorem:estimate of the multilinear form} and the trace inequality $\left\|e_u\right\|_{H^{1/2}(\Gamma_h[{\bf x}^*])} \leq C \left\|e_u\right\|_{H^1(\Omega_h[{\bf x}^*])}$, we obtain
    \begin{equation*}
        \begin{aligned}
            & \mathbf{e}_{\mathbf{u}}^T\left(\mathbf{f}_u(\mathbf{x}, \mathbf{e}_\mathbf{H} + \mathbf{H}^*) - \mathbf{f}_u(\mathbf{x}^*, \mathbf{H}^*) \right) \\
            & \leq C \left\|e_u\right\|_{L^2(\Omega_h[{\bf x}^*])} \left\|e_x\right\|_{H^1(\Omega_h[{\bf x}^*])} \\
            & \qquad  + C \left\|e_u\right\|_{H_h^{1/2}(\Gamma_h[{\bf x}^*])} \left(\left\|e_H + e_Q\right\|_{H_h^{-1/2}(\Gamma_h[{\bf x}^*])} + \left\|e_x\right\|_{H_h^{1/2}(\Gamma_h[{\bf x}^*])} \right) \\
            & \leq C \left\|e_u\right\|_{H^1(\Omega_h[{\bf x}^*])} \left(\left\|e_H + e_Q\right\|_{H_h^{-1/2}(\Gamma_h[{\bf x}^*])} + \left\|e_x\right\|_{H^1(\Omega_h[{\bf x}^*])} \right).
        \end{aligned}
    \end{equation*}
    For the second term on the right-hand side, its bulk term can be estimated with \cite[Lemma 5.1]{Edelmann_harmonicvelo}, its surface term can be estimated with \cite[Lemma 4.1]{KLLP2017} and Theorem \ref{theorem:estimate of the multilinear form} as
    \begin{equation*}
        \begin{aligned}
            -\mathbf{e}_{\mathbf{u}}^T\left(\boldsymbol{\gamma}^T\mathbf{M}_{\Gamma}(\mathbf{x})\boldsymbol{\gamma} - \boldsymbol{\gamma}^T\mathbf{M}_{\Gamma}(\mathbf{x}^*)\boldsymbol{\gamma}\right)\mathbf{u}^*& = - \int_0^1 d\theta \int_{\Gamma_h^\theta} e_u^\theta u_h^{\theta, *} \nabla_{\Gamma_h^\theta} \cdot e_x^\theta \\
            & \leq C \left\|e_u\right\|_{H_h^{1/2}(\Gamma_h[{\bf x}^*])} \left\|e_x\right\|_{H_h^{1/2}(\Gamma_h[{\bf x}^*])} \\
            & \leq C \left\|e_u\right\|_{H^1(\Omega_h[{\bf x}^*])} \left\|e_x\right\|_{H^1(\Omega_h[{\bf x}^*])}.
        \end{aligned}
    \end{equation*}
    Therefore, we obtain
    \begin{equation*}
        \begin{aligned}
            - \mathbf{e}_{\mathbf{u}}^T\left(\mathbf{L}(\mathbf{x}) - \mathbf{L}(\mathbf{x}^*)\right)\mathbf{u}^* \leq C \left\|e_u\right\|_{H^1(\Omega_h[{\bf x}^*])} \left\|e_x\right\|_{H^1(\Omega_h[{\bf x}^*])}.
        \end{aligned}
    \end{equation*}
    The last term is simply bounded by $C \left\|e_u\right\|_{L^2(\Omega_h[{\bf x}^*])} \left\|\mathbf{d}_{\mathbf{u}}\right\|_{\mathbf{M}_{\bar{\Omega}}[{\bf x}^*]}$. 
    
    Combining these bounds with the fact $\left\|e_Q\right\|_{L^2(\Gamma_h[{\bf x}^*])} \leq C \left\|e_x\right\|_{L^2(\Gamma_h[{\bf x}^*])}$ from \cite{EKL24}, we deduce that
    \begin{equation}\label{eq:error bounds, e_u, H1}
        \begin{aligned}
        \left\|e_u\right\|_{H^1(\Omega_h[{\bf x}^*])} \leq C \left(\left\|e_x\right\|_{H^1(\Omega_h[{\bf x}^*])} + \left\|e_H\right\|_{H_h^{-1/2}(\Gamma_h[{\bf x}^*])} + \delta \right).
        \end{aligned}
    \end{equation}
    
    Using the same arguments as in bounding $e_u$ before, we obtain the following estimate of $\dot{e}_u$ as
    \begin{equation}\label{eq:error bounds, e_u, H1, material derivative}
        \begin{aligned}
        \left\|\dot{e}_u\right\|_{H^1(\Omega_h[{\bf x}^*])} &\leq C \left( \left\|e_u\right\|_{H^1(\Omega_h[{\bf x}^*])} + \left\|e_x\right\|_{H^1(\Omega_h[{\bf x}^*])} + \left\|e_v\right\|_{H^1(\Omega_h[{\bf x}^*])} \right.\\
        &\qquad \left. + \left\|e_H\right\|_{H_h^{-1/2}(\Gamma_h[{\bf x}^*])} + \left\|\dot{e}_H\right\|_{H_h^{-1/2}(\Gamma_h[{\bf x}^*])} + \delta \right)
        \end{aligned}
    \end{equation}

    \textit{(B) Error bounds for the forced mean curvature flow: } 
    In contrast to the $H^1$ estimate, we test \eqref{eq:forced-mcf, error equation for nu} with $\mathbf{w}:=(I - \Delta_h)^{-1/2} \dot{e}_n$ and obtain
    \begin{equation*}
        \begin{aligned}
            \mathbf{w}^T\mathbf{M}(\mathbf{x}^*) \dot{\mathbf{e}}_\mathbf{n} + \mathbf{w}^T\mathbf{A}(\mathbf{x}^*) \mathbf{e}_\mathbf{n} & = - \mathbf{w}^T\left(\mathbf{M}(\mathbf{x}) - \mathbf{M}(\mathbf{x}^*)\right) \dot{\mathbf{n}}^* - \mathbf{w}^T\left(\mathbf{A}(\mathbf{x}) - \mathbf{A}(\mathbf{x}^*)\right) \mathbf{n}^* \\
            & \quad - \mathbf{w}^T\left(\mathbf{M}(\mathbf{x}) - \mathbf{M}(\mathbf{x}^*)\right) \dot{\mathbf{e}}_\mathbf{n} - \mathbf{w}^T\left(\mathbf{A}(\mathbf{x}) - \mathbf{A}(\mathbf{x}^*)\right) \mathbf{e}_\mathbf{n} \\
            & \quad + \mathbf{w}^T\left(\mathbf{f}_{n}(\mathbf{x}, \mathbf{n}) - \mathbf{f}_{n}(\mathbf{x}^*, \mathbf{n}^*) \right) - \mathbf{w}^T\mathbf{D}(\mathbf{x}) \boldsymbol{\gamma} \mathbf{e}_\mathbf{u} \\
            & \quad - \mathbf{w}^T\left(\mathbf{D}(\mathbf{x}) - \mathbf{D}(\mathbf{x}^*)\right) \boldsymbol{\gamma} \mathbf{u}^* - \mathbf{w}^T\mathbf{M}_{\Gamma}(\mathbf{x}^*) \mathbf{d}_{\mathbf{n}}.
        \end{aligned}
    \end{equation*}
    
    For the first term on the left-hand side, it is easy to see that
    \begin{equation*}
        \begin{aligned}
            \mathbf{w}^T\mathbf{M}(\mathbf{x}^*) \dot{\mathbf{e}}_\mathbf{n} = \int_{\Gamma_h[\mathbf{x}^*]} \dot{e}_n (I - \Delta_h)^{-1/2} \dot{e}_n = \left\|\dot{e}_n\right\|_{H^{-1/2}_h(\Gamma_h[\mathbf{x}^*])}^2.
        \end{aligned}
    \end{equation*}
    The second term can be estimated by using \eqref{eq:estimate of the multilinear form, case 1} and the derivative trick \eqref{eq: time derivative of discrete Laplacian of order 1/2, 2} as follows:
    \begin{equation*}
        \begin{aligned}
            \mathbf{w}^T\mathbf{A}(\mathbf{x}^*) \mathbf{e}_\mathbf{n} &= \int_{\Gamma_h[\mathbf{x}^*]} \nabla_{\Gamma_h[\mathbf{x}^*]}((I - \Delta_h)^{-1/2} \dot{e}_n) \cdot \nabla_{\Gamma_h[\mathbf{x}^*]} e_n \\
            & = \int_{\Gamma_h[\mathbf{x}^*]} \dot{e}_n (I - \Delta_h)^{1/2} e_n - \int_{\Gamma_h[\mathbf{x}^*]} \dot{e}_n (I - \Delta_h)^{-1/2} e_n \\
            & \geq \frac{1}{2} \frac{d}{dt} \left\| e_n \right\|_{H_h^{1/2}(\Gamma_h[\mathbf{x}^*])}^2 - \int_{\Gamma_h[\mathbf{x}^*]} e_n \left(\partial^\bullet \left(I - \Delta_h\right)^{1/2}\right) e_n \\
            & \qquad - \int_{\Gamma_h[\mathbf{x}^*]} e_n \left(I - \Delta_h\right)^{1/2} e_n \nabla_{\Gamma_h[\mathbf{x}^*]} \cdot v_h \\
            & \qquad - C \left\|\dot{e}_n\right\|_{H_h^{-1/2}(\Gamma_h[\mathbf{x}^*])} \left\|e_n\right\|_{H_h^{-1/2}(\Gamma_h[\mathbf{x}^*])} \\
            & \geq \frac{1}{2} \frac{d}{dt} \left\| e_n \right\|_{H_h^{1/2}(\Gamma_h[\mathbf{x}^*])}^2 - C \left\|e_n\right\|_{H_h^{1/2}(\Gamma_h[\mathbf{x}^*])}^2 \\
            & \qquad  - C \left\|\dot{e}_n\right\|_{H_h^{-1/2}(\Gamma_h[\mathbf{x}^*])} \left\|e_n\right\|_{H_h^{-1/2}(\Gamma_h[\mathbf{x}^*])}.
        \end{aligned}
    \end{equation*}

    For the first and seventh terms on the right-hand side, \cite[Lemma 4.1]{KLLP2017}, \cite[(7.26)]{MCF}, and \eqref{eq:estimate of the multilinear form, case 3} yield that
    \begin{equation*}
        \begin{aligned}
            &-\mathbf{w}^T\left(\mathbf{M}(\mathbf{x}) - \mathbf{M}(\mathbf{x}^*)\right) \dot{\mathbf{n}}- \mathbf{w}^T\left(\mathbf{D}(\mathbf{x}) - \mathbf{D}(\mathbf{x}^*)\right) \boldsymbol{\gamma} \mathbf{u}^*\\
            & = \int_0^1 d\theta \int_{\Gamma_h^\theta} w_h^\theta \cdot\left( -\dot{n}_h^\theta - \nabla_{\Gamma_h^\theta} u_h^{*, \theta}\right) \nabla_{\Gamma_h^\theta} \cdot e_x^\theta \\
            & \qquad + \int_0^1 d\theta \int_{\Gamma_h^\theta} w_h^\theta \left(\nabla_{\Gamma_h^\theta} e_x^\theta - n_h^\theta(n_h^\theta)^T(\nabla_{\Gamma_h^\theta} e_x^\theta)^T\right) \nabla_{\Gamma_h^\theta} u_h^{*, \theta}   \\
            & \leq C \left\|w_h\right\|_{H_h^{1/2}(\Gamma_h[\mathbf{x}^*])} \left\|e_x\right\|_{H_h^{1/2}(\Gamma_h[\mathbf{x}^*])} \\
            & \leq C \left\|\dot{e}_n\right\|_{H^{-1/2}_h(\Gamma_h[\mathbf{x}^*])} \left\|e_x\right\|_{H_h^{1/2}(\Gamma_h[\mathbf{x}^*])}.
        \end{aligned}
    \end{equation*}
    Here we use the fact $\nabla_{\Gamma_h^\theta} \text{id}_h^\theta = I_3 - n_h^\theta(n_h^\theta)^T$ and $\text{id}_h \in \mathscr{S}_h[\mathbf{x}^*]$.

    The third and fourth terms can be bounded by the standard $L^2$-$L^2$-$L^\infty$ estimate and the inverse estimate with the fact $\left\|e_n\right\|_{W^{1, \infty}(\Gamma_h[\mathbf{x}^*])} \leq h^{(\kappa - 2)/2 + 1/2}$ as follows:
    \begin{equation}\label{eq: standard 22inf estimate}
        \begin{aligned}
            &\mathbf{w}^T\left(\mathbf{M}(\mathbf{x}) - \mathbf{M}(\mathbf{x}^*)\right) \dot{\mathbf{e}}_\mathbf{n} + \mathbf{w}^T\left(\mathbf{A}(\mathbf{x}) - \mathbf{A}(\mathbf{x}^*)\right) \mathbf{e}_\mathbf{n} \\
            & \leq C \left\|w_h\right\|_{H^1(\Gamma_h[\mathbf{x}^*])} \left\|e_x\right\|_{H^1(\Gamma_h[\mathbf{x}^*])} \left\|e_n\right\|_{W^{1, \infty}(\Gamma_h[\mathbf{x}^*])} \\
            & \leq C h^{(\kappa - 2)/2 + 1/2}\left\|\dot{e}_n\right\|_{L^2(\Gamma_h[\mathbf{x}^*])} \left(h^{-1/2}\left\|e_x\right\|_{H_h^{1/2}(\Gamma_h[\mathbf{x}^*])}\right) \\
            & \leq C h^{(\kappa - 2)/2} \left\|\dot{e}_n\right\|_{H_h^{1/2}(\Gamma_h[\mathbf{x}^*])} \left\|e_x\right\|_{H_h^{1/2}(\Gamma_h[\mathbf{x}^*])}
        \end{aligned}
    \end{equation}

    Denote the nodal vector of $(I - \Delta_h)^{-1/2} e_n$ by $\mathbf{C} \mathbf{e}_{\mathbf{n}}$, we know that $\mathbf{w} = \mathbf{C} \dot{\mathbf{e}}_{\mathbf{n}}$. For the second term, use the derivative trick in \cite[(7.24)]{MCF} and \eqref{eq:estimate of the multilinear form, case 4} to obtain
    \begin{equation*}
        \begin{aligned}
            &- \mathbf{w}^T\left(\mathbf{A}(\mathbf{x}) - \mathbf{A}(\mathbf{x}^*)\right) \mathbf{n}^* \\
            & = - \frac{d}{dt} \left((\mathbf{C} \mathbf{e}_{\mathbf{n}})^T (\mathbf{A}(\mathbf{x}) - \mathbf{A}(\mathbf{x}^*)) \mathbf{n}^* \right) + \left( \frac{d}{dt} \mathbf{C} \mathbf{e}_{\mathbf{n}} \right)^T (\mathbf{A}(\mathbf{x}) - \mathbf{A}(\mathbf{x}^*)) \mathbf{n}^* \\
            & \qquad + \left( \mathbf{C} \mathbf{e}_{\mathbf{n}} \right)^T \left( \frac{d}{dt} (\mathbf{A}(\mathbf{x}) - \mathbf{A}(\mathbf{x}^*)) \right) \mathbf{n}^*  + \left( \mathbf{C} \mathbf{e}_{\mathbf{n}} \right)^T \left( \mathbf{A}(\mathbf{x}) - \mathbf{A}(\mathbf{x}^*) \right) \dot{\mathbf{n}}^*\\
            & \leq - \frac{d}{dt} \left((\mathbf{C} \mathbf{e}_{\mathbf{n}})^T (\mathbf{A}(\mathbf{x}) - \mathbf{A}(\mathbf{x}^*)) \mathbf{n}^* \right) \\
            & \qquad + C \left\| \left(\partial^\bullet (I - \Delta_h)^{-1/2}\right) e_n \right\|_{H_h^{3/2}(\Gamma_h[\mathbf{x}^*])} \left\| e_x \right\|_{H_h^{1/2}(\Gamma_h[\mathbf{x}^*])} \\
            & \qquad + C \left\| \left( I - \Delta_h \right)^{-1/2} e_n \right\|_{H_h^{3/2}(\Gamma_h[\mathbf{x}^*])} \left( \left\|e_x\right\|_{H_h^{1/2}(\Gamma_h[\mathbf{x}^*])} + \left\|e_v\right\|_{H_h^{1/2}(\Gamma_h[\mathbf{x}^*])} \right).
        \end{aligned}
    \end{equation*}
    It is easy to see that $\left\| \left( I - \Delta_h \right)^{-1/2} e_n \right\|_{H_h^{3/2}(\Gamma_h[\mathbf{x}^*])} = \left\|e_n \right\|_{H_h^{1/2}(\Gamma_h[\mathbf{x}^*])}$. For the time derivative term, using \eqref{eq: time derivative of discrete Laplacian of order -1/2, representation} and \eqref{eq: time derivative of discrete Laplacian of order 1/2, 3}, we have
    \begin{equation*}
        \begin{aligned}
            &\left\| \left(\partial^\bullet (I - \Delta_h)^{-1/2}\right) e_n \right\|_{H_h^{3/2}(\Gamma_h[\mathbf{x}^*])} \\
            & = \sup_{0 \neq z_h \in \mathscr{S}_h[\mathbf{x}^*]} \frac{\int_{\Gamma_h[\mathbf{x}^*]} \left((I - \Delta_h)^{3/4}\left(\partial^\bullet (I - \Delta_h)^{-1/2}\right) e_n\right)  z_h}{\left\|z_h\right\|_{L^2(\Gamma_h[\mathbf{x}^*])}} \\
            & = \sup_{0 \neq z_h \in \mathscr{S}_h[\mathbf{x}^*]} \frac{-\int_{\Gamma_h[\mathbf{x}^*]} \left((I - \Delta_h)^{1/4}\left(\left(\partial^\bullet \left(I - \Delta_h\right)^{1/2}\right) \left(I - \Delta_h\right)^{-1/2}\right) e_n\right)  z_h}{\left\|z_h\right\|_{L^2(\Gamma_h[\mathbf{x}^*])}} \\
            & \leq C \sup_{0 \neq z_h \in \mathscr{S}_h[\mathbf{x}^*]} \frac{\left\| (I - \Delta_h)^{-1/2} e_n \right\|_{H_h^{3/2}(\Gamma_h[\mathbf{x}^*])} \left\| (I - \Delta_h)^{1/4} z_h \right\|_{H_h^{-1/2}(\Gamma_h[\mathbf{x}^*])}}{\left\|z_h\right\|_{L^2(\Gamma_h[\mathbf{x}^*])}} \\
            & \leq C \left\| e_n \right\|_{H_h^{1/2}(\Gamma_h[\mathbf{x}^*])}.
        \end{aligned}
    \end{equation*}
    Combining these estimates yields that
    \begin{equation*}
        \begin{aligned}
            &-\mathbf{w}^T\left(\mathbf{A}(\mathbf{x}) - \mathbf{A}(\mathbf{x}^*)\right) \mathbf{n}^* \\
            & \leq - \frac{d}{dt} \left((\mathbf{C} \mathbf{e}_{\mathbf{n}})^T (\mathbf{A}(\mathbf{x}) - \mathbf{A}(\mathbf{x}^*)) \mathbf{n}^* \right) \\
            & \qquad + C \left\|e_n\right\|_{H_h^{1/2}(\Gamma_h[\mathbf{x}^*])} \left( \left\|e_x\right\|_{H_h^{1/2}(\Gamma_h[\mathbf{x}^*])} + \left\|e_v\right\|_{H_h^{1/2}(\Gamma_h[\mathbf{x}^*])} \right).
        \end{aligned}
    \end{equation*}

    From page 24 in \cite{MCF}, the fifth term can be written as
    \begin{equation*}
        \begin{aligned}
            & \mathbf{w}^T\left(\mathbf{f}_{n}(\mathbf{x}, \mathbf{n}) - \mathbf{f}_{n}(\mathbf{x}^*, \mathbf{n}^*) \right) \\
            & = \int_0^1  d\theta \int_{\Gamma_h^\theta} w_h^\theta \cdot \left(f(n_{\Gamma_h^\theta}, \nabla_{\Gamma_h^\theta} n_{\Gamma_h^\theta}) \nabla_{\Gamma_h^\theta} \cdot e_x^\theta + \partial_1 f(n_{\Gamma_h^\theta}, \nabla_{\Gamma_h^\theta} n_{\Gamma_h^\theta}) e_n^\theta \right. \\
            & \qquad \left. + \partial_2 f(n_{\Gamma_h^\theta}, \nabla_{\Gamma_h^\theta} n_{\Gamma_h^\theta}) \left(\nabla_{\Gamma_h^\theta} e_n^\theta - \left( \nabla_{\Gamma_h^\theta} e_x^\theta - n_h^\theta(n_h^\theta)^T(\nabla_{\Gamma_h^\theta} e_x^\theta)^T\right) \nabla_{\Gamma_h^\theta} n_{\Gamma_h^\theta} \right) \right),
        \end{aligned}
    \end{equation*}
    with $n_{\Gamma_h^\theta} = n_h^{*, \theta} + \theta e_n^\theta$. Apply the identity $\nabla_{\Gamma_h^\theta} \text{id}_h^\theta = I_3 - n_h^\theta(n_h^\theta)^T$ again together with \eqref{eq:estimate of the multilinear form, case 1} and \eqref{eq:estimate of the multilinear form, case 3}, we further have
    \begin{equation*}
        \begin{aligned}
            & \mathbf{w}^T\left(\mathbf{f}_{n}(\mathbf{x}, \mathbf{n}) - \mathbf{f}_{n}(\mathbf{x}^*, \mathbf{n}^*) \right) \\
            & \leq C \left\|w_h\right\|_{H_h^{1/2}(\Gamma_h[\mathbf{x}^*])} \left(\left\|e_n\right\|_{H_h^{1/2}(\Gamma_h[\mathbf{x}^*])} + \left\|e_x\right\|_{H_h^{1/2}(\Gamma_h[\mathbf{x}^*])} \right) \\
            & \leq C \left\|\dot{e}_n\right\|_{H_h^{-1/2}(\Gamma_h[\mathbf{x}^*])} \left(\left\|e_n\right\|_{H_h^{1/2}(\Gamma_h[\mathbf{x}^*])} + \left\|e_x\right\|_{H_h^{1/2}(\Gamma_h[\mathbf{x}^*])} \right).
        \end{aligned}
    \end{equation*}

    For the rest two terms, \eqref{eq:estimate of the multilinear form, case 3} and standard $L^2$-$L^2$ estimate yield that
    \begin{equation*}
        \begin{aligned}
            &- \mathbf{w}^T\mathbf{D}(\mathbf{x}) \boldsymbol{\gamma} \mathbf{e}_\mathbf{u}  - \mathbf{w}^T\mathbf{M}_{\Gamma}(\mathbf{x}^*) \mathbf{d}_{\mathbf{n}} \\
            & \qquad \leq C \left\|\dot{e}_n\right\|_{H_h^{-1/2}(\Gamma_h[\mathbf{x}^*])} \left(\left\|e_u\right\|_{H_h^{1/2}(\Gamma_h[\mathbf{x}^*])} + \left\|d_n\right\|_{L^2(\Gamma_h[\mathbf{x}^*])}\right)
        \end{aligned}
    \end{equation*}

    Combining all the above estimates yields that
    \begin{equation*}
        \begin{aligned}
            & \left\|\dot{e}_n\right\|_{H^{-1/2}_h(\Gamma_h[\mathbf{x}^*])}^2 + \frac{1}{2} \frac{d}{dt} \left\| e_n \right\|_{H_h^{1/2}(\Gamma_h[\mathbf{x}^*])}^2 \\
            & \leq C \left\|e_n\right\|_{H_h^{1/2}(\Gamma_h[\mathbf{x}^*])}^2  + C \left\|\dot{e}_n\right\|_{H_h^{-1/2}(\Gamma_h[\mathbf{x}^*])} \left\|e_n\right\|_{H_h^{1/2}(\Gamma_h[\mathbf{x}^*])} \\
            & \qquad + C \left\|\dot{e}_n\right\|_{H_h^{-1/2}(\Gamma_h[\mathbf{x}^*])} \left(\left\|e_u\right\|_{H_h^{1/2}(\Gamma_h[\mathbf{x}^*])} + \left\|e_x\right\|_{H_h^{1/2}(\Gamma_h[\mathbf{x}^*])} + \left\|d_n\right\|_{L^2(\Gamma_h[\mathbf{x}^*])}\right) \\
            & \qquad + C \left\|e_n\right\|_{H_h^{1/2}(\Gamma_h[\mathbf{x}^*])} \left( \left\|e_x\right\|_{H_h^{1/2}(\Gamma_h[\mathbf{x}^*])} + \left\|e_v\right\|_{H_h^{1/2}(\Gamma_h[\mathbf{x}^*])} \right)\\
            & \qquad - \frac{d}{dt} \left((\mathbf{C} \mathbf{e}_{\mathbf{n}})^T (\mathbf{A}(\mathbf{x}) - \mathbf{A}(\mathbf{x}^*)) \mathbf{n}^* \right) .
        \end{aligned}
    \end{equation*}
    Using \eqref{eq:estimate of the multilinear form, case 1} and \eqref{eq: time derivative of discrete Laplacian of order -1/2, 1}, the first term can be further bounded as follows:
    \begin{equation}\label{eq:error bound, time derivative of e_n, H-1/2}
        \begin{aligned}
            &\frac{1}{2}\frac{d}{dt}  \left\| e_n \right\|_{H_h^{-1/2}(\Gamma_h[\mathbf{x}^*])}^2 \\
            & = \int_{\Gamma_h[\mathbf{x}^*]}  \left((I - \Delta_h)^{-1/2} e_n\right) \cdot \left( \dot{e}_n + \frac{1}{2} e_n \nabla_{\Gamma_h[\mathbf{x}^*]} \cdot v_h \right) \\
            & \qquad + \frac{1}{2} \int_{\Gamma_h[\mathbf{x}^*]} e_n \left( \partial^\bullet (I - \Delta_h)^{-1/2}\right) e_n  \\
            & \leq C \left\|(I - \Delta_h)^{-1/2} e_n \right\|_{H_h^{1/2}(\Gamma_h[\mathbf{x}^*])} \left(\left\| \dot{e}_n \right\|_{H_h^{-1/2}(\Gamma_h[\mathbf{x}^*])} + \left\| e_n \right\|_{H_h^{-1/2}(\Gamma_h[\mathbf{x}^*])} \right) \\
            & \qquad + C \left\| e_n \right\|_{H_h^{-1/2}(\Gamma_h[\mathbf{x}^*])}^2 \\
            & \leq \frac{1}{2} \left\| \dot{e}_n \right\|_{H_h^{-1/2}(\Gamma_h[\mathbf{x}^*])}^2 + C \left\| e_n \right\|_{H_h^{-1/2}(\Gamma_h[\mathbf{x}^*])}^2.
        \end{aligned}
    \end{equation}
    Therefore, we have the following error bound for $e_n$:
    \begin{equation}\label{eq:error bounds, e_n, H1/2}
        \begin{aligned}
            &\frac{1}{2}\frac{d}{dt}  \left\| e_n \right\|_{H_h^{-1/2}(\Gamma_h[\mathbf{x}^*])}^2 + \frac{1}{2} \frac{d}{dt} \left\| e_n \right\|_{H_h^{1/2}(\Gamma_h[\mathbf{x}^*])}^2 \\
            & \leq C \left( \left\|e_u\right\|_{H_h^{1/2}(\Gamma_h[\mathbf{x}^*])}^2 + \left\|e_x\right\|_{H_h^{1/2}(\Gamma_h[\mathbf{x}^*])}^2 + \left\|e_v\right\|_{H_h^{1/2}(\Gamma_h[\mathbf{x}^*])}^2\right) \\
            & \quad + C \left\|e_n\right\|_{H_h^{1/2}(\Gamma_h[\mathbf{x}^*])}^2 + C \delta^2  - \frac{d}{dt} \left((\mathbf{C} \mathbf{e}_{\mathbf{n}})^T (\mathbf{A}(\mathbf{x}) - \mathbf{A}(\mathbf{x}^*)) \mathbf{n}^* \right).
        \end{aligned}
    \end{equation}

    Using the same argument for $e_H$ and the same steps as in \cite{MCF_soldriven}, we can obtain that
    \begin{equation}\label{eq:error bounds, e_H, H1/2}
        \begin{aligned}
            & \left\| \dot{e}_H \right\|_{H_h^{-1/2}(\Gamma_h[\mathbf{x}^*])}^2 + \frac{1}{2} \frac{d}{dt} \left\| e_H \right\|_{H_h^{1/2}(\Gamma_h[\mathbf{x}^*])}^2 \\
            & \leq C \left( \left\|e_u\right\|_{H_h^{1/2}(\Gamma_h[\mathbf{x}^*])}^2  + \left\|e_x\right\|_{H_h^{1/2}(\Gamma_h[\mathbf{x}^*])}^2 + \left\|e_v\right\|_{H_h^{1/2}(\Gamma_h[\mathbf{x}^*])}^2\right) \\
            & \qquad + C \left\|e_n\right\|_{H_h^{1/2}(\Gamma_h[\mathbf{x}^*])}^2 + C \left\|e_H\right\|_{H_h^{1/2}(\Gamma_h[\mathbf{x}^*])}^2 + C \delta^2 \\
            & \qquad + \frac{1}{2} \rho \left\|\dot{e}_u\right\|_{H_h^{1/2}(\Gamma_h[\mathbf{x}^*])}^2 - \frac{d}{dt} \left((\mathbf{C} \mathbf{e}_{\mathbf{H}})^T (\mathbf{A}(\mathbf{x}) - \mathbf{A}(\mathbf{x}^*))( \mathbf{H}^* - \boldsymbol{\gamma} \mathbf{u}^*) \right).
        \end{aligned}
    \end{equation}
    Here $\rho > 0$ can be chosen sufficiently small.

    Follow the same argument of Lemma 5.3 in \cite{kovacs2021convergent}, $e_v$ is bounded by
    \begin{equation}\label{eq:error bounds, e_v_Gamma, H1/2}
        \begin{aligned}
            \left\| e_v \right\|_{H_h^{1/2}(\Gamma_h[\mathbf{x}^*])} &\leq   C \left\|e_H\right\|_{H_h^{1/2}(\Gamma_h[\mathbf{x}^*])} + C \left\|e_n\right\|_{H_h^{1/2}(\Gamma_h[\mathbf{x}^*])} \\
            & \qquad + C  \left\|e_u\right\|_{H_h^{1/2}(\Gamma_h[\mathbf{x}^*])}+ C \left\|d_v\right\|_{H_h^{1/2}(\Gamma_h[\mathbf{x}^*])} \\
            &\leq   C \left\|e_H\right\|_{H_h^{1/2}(\Gamma_h[\mathbf{x}^*])} + C \left\|e_n\right\|_{H_h^{1/2}(\Gamma_h[\mathbf{x}^*])} \\
            & \qquad + C  \left\|e_u\right\|_{H_h^{1/2}(\Gamma_h[\mathbf{x}^*])}+ C \delta.
        \end{aligned}
    \end{equation}

    For $e_x$, the trace inequality simply implies
    \begin{equation}\label{eq:error bounds, e_x_Gamma, H1/2}
        \begin{aligned}
            \left\| e_x \right\|_{H_h^{1/2}(\Gamma_h[\mathbf{x}^*])} \leq C \left\| e_x \right\|_{H^1(\Omega_h[\mathbf{x}^*])}
        \end{aligned}
    \end{equation}

    \textit{(C) Error bounds for the harmonic velocity extension: } This error bound using the $H^{1/2}(\Gamma_h[\mathbf{x}^*])$ norm is already established in (4.34), (4.35), and (4.36) in \cite{EKL24}. 
    \begin{equation}\label{eq:error bounds, e_v_Omega, H1}
        \begin{aligned}
            \left\| e_v \right\|_{H^1(\Omega_h[\mathbf{x}^*])} &\leq C \left\| e_v \right\|_{H_h^{1/2}(\Gamma_h[\mathbf{x}^*])} + C \left\| e_x \right\|_{H^1(\Omega_h[\mathbf{x}^*])} + \left\|\mathbf{d}_{\mathbf{v}}\right\|_{\star, {\mathbf{x}^*}} \\
            &\leq C \left\| e_v \right\|_{H_h^{1/2}(\Gamma_h[\mathbf{x}^*])} + C \left\| e_x \right\|_{H^1(\Omega_h[\mathbf{x}^*])} + \delta.
        \end{aligned}
    \end{equation}
    \begin{equation}\label{eq:error bounds, e_x_Omega, H1}
        \begin{aligned}
            \left\| e_x(t) \right\|_{H^1(\Omega_h[\mathbf{x}^*(t)])} \leq C \int_0^t \left\| e_v(s) \right\|_{H^1(\Omega_h[\mathbf{x}^*(s)])} \, ds.
        \end{aligned}
    \end{equation}
    The second inequality in \eqref{eq:error bounds, e_v_Omega, H1} is due to the norm equivalence of $\left\| \cdot \right\|_{H_h^{1/2}(\Gamma_h[\mathbf{x}^*])}$ and $\left\| \cdot \right\|_{H^{1/2}(\Gamma_h[\mathbf{x}^*])}$.

    \textit{(D) Combining the error bounds: } First, we eliminate $e_u$ by inserting the error bound for $e_u$ \eqref{eq:error bounds, e_u, H1} and $\dot{e}_u$ \eqref{eq:error bounds, e_u, H1, material derivative} into \eqref{eq:error bounds, e_n, H1/2}, \eqref{eq:error bounds, e_H, H1/2}, and \eqref{eq:error bounds, e_v_Gamma, H1/2}. Together with the trace inequality $\left\|e_u\right\|_{H_h^{1/2}(\Gamma_h[\mathbf{x}^*])} \leq C \left\|e_u\right\|_{H^1(\Omega_h[\mathbf{x}^*])}$, we obtain
    \begin{equation}\label{eq:error bounds, e_n, H1/2, step 1}
        \begin{aligned}
            &\frac{1}{2}\frac{d}{dt}  \left\| e_n \right\|_{H_h^{-1/2}(\Gamma_h[\mathbf{x}^*])}^2 + \frac{1}{2} \frac{d}{dt} \left\| e_n \right\|_{H_h^{1/2}(\Gamma_h[\mathbf{x}^*])}^2 \\
            & \leq C \left( \left\|e_H\right\|_{H_h^{-1/2}(\Gamma_h[\mathbf{x}^*])}^2 + \left\|e_x\right\|_{H^1(\Omega_h[\mathbf{x}^*])}^2 + \left\|e_v\right\|_{H^1(\Omega_h[\mathbf{x}^*])}^2\right) \\
            & \qquad + C \left\|e_n\right\|_{H_h^{1/2}(\Gamma_h[\mathbf{x}^*])}^2 + C \delta^2 - \frac{d}{dt} \left((\mathbf{C} \mathbf{e}_{\mathbf{n}})^T (\mathbf{A}(\mathbf{x}) - \mathbf{A}(\mathbf{x}^*)) \mathbf{n}^* \right).
        \end{aligned}
    \end{equation}
    \begin{equation}\label{eq:error bounds, e_H, H1/2, step 1}
        \begin{aligned}
            &\frac{1}{2}\frac{d}{dt}  \left\| e_H \right\|_{H_h^{-1/2}(\Gamma_h[\mathbf{x}^*])}^2 + \frac{1}{2} \frac{d}{dt} \left\| e_H \right\|_{H_h^{1/2}(\Gamma_h[\mathbf{x}^*])}^2 \\
            & \leq C \left( \left\|e_H\right\|_{H_h^{-1/2}(\Gamma_h[\mathbf{x}^*])}^2  + \left\|e_x\right\|_{H^1(\Omega_h[\mathbf{x}^*])}^2 + \left\|e_v\right\|_{H^1(\Omega_h[\mathbf{x}^*])}^2\right) \\
            & \quad + C \left\|e_n\right\|_{H_h^{1/2}(\Gamma_h[\mathbf{x}^*])}^2 + C \left\|e_H\right\|_{H_h^{1/2}(\Gamma_h[\mathbf{x}^*])}^2 + C \delta^2 \\
            & \qquad - \frac{d}{dt} \left((\mathbf{C} \mathbf{e}_{\mathbf{H}})^T (\mathbf{A}(\mathbf{x}) - \mathbf{A}(\mathbf{x}^*))( \mathbf{H}^* - \boldsymbol{\gamma} \mathbf{u}^*) \right).
        \end{aligned}
    \end{equation}
    \begin{equation}\label{eq:error bounds, e_v_Gamma, H1/2, step 1}
        \begin{aligned}
            \left\| e_v \right\|_{H_h^{1/2}(\Gamma_h[\mathbf{x}^*])} &\leq C \left\|e_H \right\|_{H_h^{-1/2}(\Gamma_h[\mathbf{x}^*])} + C \left\|e_H\right\|_{H_h^{1/2}(\Gamma_h[\mathbf{x}^*])} \\
            & \qquad  + C \left\|e_n\right\|_{H_h^{1/2}(\Gamma_h[\mathbf{x}^*])} + C  \left\|e_x\right\|_{H^1(\Omega_h[\mathbf{x}^*])}+ \delta.
        \end{aligned}
    \end{equation}
    Here for \eqref{eq:error bounds, e_H, H1/2, step 1}, the $\left\|\dot{e}_H\right\|_{H_h^{-1/2}(\Gamma_h[\mathbf{x}^*])}$ in $\left\| \dot{e}_u\right\|_{H^1(\Omega_h[\mathbf{x}^*])}$ is absorbed since $\rho$ is sufficiently small, then we use \eqref{eq:error bound, time derivative of e_n, H-1/2} to get $\frac{d}{dt} \left\| e_H\right\|_{H_h^{-1/2}(\Gamma_h[\mathbf{x}^*])}$ from $\left\| \dot{e}_H\right\|_{H_h^{-1/2}(\Gamma_h[\mathbf{x}^*])}$.

    We introduce the notation $e_w = (e_H, e_n)^T$. By adding \eqref{eq:error bounds, e_n, H1/2, step 1}, \eqref{eq:error bounds, e_H, H1/2, step 1} together, taking integration in time, and applying the Gronwall inequality, we obtain
    \begin{equation*}
        \begin{aligned}
            &\frac{1}{2} \left\| e_w(t) \right\|_{H_h^{-1/2}(\Gamma_h[\mathbf{x}^*(t)])}^2 + \frac{1}{2} \left\| e_w(t) \right\|_{H_h^{1/2}(\Gamma_h[\mathbf{x}^*(t)])}^2 \\
            & \leq \frac{1}{2} (\varepsilon^0)^2 + C \int_0^t \left(\left\|e_x(s)\right\|_{H^1(\Omega_h[\mathbf{x}^*(s)])}^2 + \left\|e_v(s)\right\|_{H^1(\Omega_h[\mathbf{x}^*(s)])}^2 + \delta^2\right)\, ds \\
            & \qquad - (\mathbf{C} \mathbf{e}_{\mathbf{n}})^T (\mathbf{A}(\mathbf{x}(t))  - \mathbf{A}(\mathbf{x}^*(t))) \mathbf{n}^*(t) \\
            & \qquad - (\mathbf{C} \mathbf{e}_{\mathbf{H}})^T (\mathbf{A}(\mathbf{x}(t)) - \mathbf{A}(\mathbf{x}^*(t))) \left(\mathbf{H}^*(t) - \boldsymbol{\gamma} \mathbf{u}^*(t)\right).
        \end{aligned}
    \end{equation*}
    For the third term, we use \eqref{eq:estimate of the multilinear form, case 4} and the Leibniz formula to get
    \begin{equation*}
        \begin{aligned}
            &(\mathbf{C} \mathbf{e}_{\mathbf{n}})^T (\mathbf{A}(\mathbf{x}(t))  - \mathbf{A}(\mathbf{x}^*(t))) \mathbf{n}^*(t) \\
            & \leq C \left\| (I - \Delta_h)^{-1/2} e_n(t) \right\|_{H_h^{3/2}(\Gamma_h[\mathbf{x}^*(t)])} \left\| e_x(t) \right\|_{H_h^{1/2}(\Gamma_h[\mathbf{x}^*(t)])} \\
            & \leq \rho \left\| e_n(t) \right\|_{H_h^{1/2}(\Gamma_h[\mathbf{x}^*(t)])}^2 + C(\rho) \left\| e_x(t) \right\|_{H^1(\Omega_h[\mathbf{x}^*(t)])}^2,
        \end{aligned}
    \end{equation*}
    where $\rho$ is sufficiently small. The last term can be estimated similarly. Absorbing $\rho \left\| e_w(t) \right\|_{H_h^{1/2}(\Gamma_h[\mathbf{x}^*(t)])}^2$ into the left-hand side, we obtain
    \begin{equation}\label{eq:error bounds, e_w, H1/2, step 2}
        \begin{aligned}
            & \left\| e_w(t) \right\|_{H_h^{-1/2}(\Gamma_h[\mathbf{x}^*(t)])}^2 + \left\| e_w(t) \right\|_{H_h^{1/2}(\Gamma_h[\mathbf{x}^*(t)])}^2 \\
            & \leq C (\varepsilon^0)^2 + C \int_0^t \left(\left\|e_x(s)\right\|_{H^1(\Omega_h[\mathbf{x}^*(s)])}^2 + \left\|e_v(s)\right\|_{H^1(\Omega_h[\mathbf{x}^*(s)])}^2 + \delta^2\right)\, ds \\
            & \qquad + \left\|e_x(t)\right\|_{H^1(\Omega_h[\mathbf{x}^*(t)])}^2 + C \int_0^t \delta^2 \, ds.
        \end{aligned}
    \end{equation}

    Next, we eliminate $e_v$ by inserting the error bound for $e_v$ \eqref{eq:error bounds, e_v_Gamma, H1/2, step 1} and \eqref{eq:error bounds, e_v_Omega, H1} into \eqref{eq:error bounds, e_w, H1/2, step 2} and \eqref{eq:error bounds, e_x_Omega, H1}. Applying the Gronwall inequality, we obtain
    \begin{equation}\label{eq:error bounds, e_w, H1/2, step 3}
        \begin{aligned}
            & \left\| e_w(t) \right\|_{H_h^{-1/2}(\Gamma_h[\mathbf{x}^*(t)])}^2 + \left\| e_w(t) \right\|_{H_h^{1/2}(\Gamma_h[\mathbf{x}^*(t)])}^2 \\
            & \leq C (\varepsilon^0)^2 + C \int_0^t \left(\left\|e_x(s)\right\|_{H^1(\Omega_h[\mathbf{x}^*(s)])}^2 + \delta^2\right)\, ds  + C\left\|e_x(t)\right\|_{H^1(\Omega_h[\mathbf{x}^*(t)])}^2.
        \end{aligned}
    \end{equation}
    \begin{equation}\label{eq:error bounds, e_x, H1, step 3}
        \begin{aligned}
            &\left\| e_x(t) \right\|_{H^1(\Omega_h[\mathbf{x}^*(t)])}\\
            &\leq C \int_0^t \left( \left\| e_w(s) \right\|_{H_h^{-1/2}(\Gamma_h[\mathbf{x}^*(s)])} + \left\| e_w(s) \right\|_{H_h^{1/2}(\Gamma_h[\mathbf{x}^*(s)])} + \delta \right) \, ds.
        \end{aligned}
    \end{equation}

    Finally, we eliminate $e_x$ by inserting the error bound for $e_x$ \eqref{eq:error bounds, e_x, H1, step 3} into \eqref{eq:error bounds, e_w, H1/2, step 3}. Applying the Gronwall inequality, we obtain
    \begin{equation}\label{eq:error bounds, e_w, H1/2, step 4}
        \begin{aligned}
            & \left\| e_w(t) \right\|_{H_h^{-1/2}(\Gamma_h[\mathbf{x}^*(t)])}^2 + \left\| e_w(t) \right\|_{H_h^{1/2}(\Gamma_h[\mathbf{x}^*(t)])}^2 \\
            & \leq C (\varepsilon^0)^2+ C \int_0^t \delta^2 \, ds \leq C \left(\delta + \varepsilon^0\right)^2,
        \end{aligned}
    \end{equation}
    which is the desired error bound \eqref{eq:stability bound}.

    To complete the proof, we need to show $t^* = T$. Using the inverse estimate and the trace inequality, noting that $\kappa -3/2 > (\kappa - 2)/2 + 1/2$ for $\kappa > 2$, we obtain
    \begin{equation*}
        \begin{aligned}
        &\left\|e_u\right\|_{W^{1, \infty}(\Omega_h[{\bf x}^*(t)])} + \left\|e_x\right\|_{W^{1, \infty}(\Omega_h[{\bf x}^*(t)])} + \left\|e_v\right\|_{W^{1, \infty}(\Omega_h[{\bf x}^*(t)])} \\
        &\leq C h^{-3/2}\left( \left\|e_u\right\|_{H^1(\Omega_h[{\bf x}^*(t)])} + \left\|e_x\right\|_{H^1(\Omega_h[{\bf x}^*(t)])} + \left\|e_v\right\|_{H^1(\Omega_h[{\bf x}^*(t)])} \right) \\
        & \leq C h^{\kappa - 3/2} \leq h^{(\kappa - 2)/2 + 1/2},
        \end{aligned}
    \end{equation*}
    and
    \begin{equation*}
        \begin{aligned}
        &\left\|e_u\right\|_{W^{1, \infty}(\Gamma_h[{\bf x}^*(t)])} + \left\|e_x\right\|_{W^{1, \infty}(\Gamma_h[{\bf x}^*(t)])} + \left\|e_v\right\|_{W^{1, \infty}(\Gamma_h[{\bf x}^*(t)])} \\
        & \qquad + \left\|e_H\right\|_{W^{1, \infty}(\Gamma_h[{\bf x}^*(t)])} + \left\|e_n\right\|_{W^{1, \infty}(\Gamma_h[{\bf x}^*(t)])} \\
        &\leq C h^{-3/2}\left( \left\|e_u\right\|_{H^{1/2}(\Gamma_h[{\bf x}^*(t)])} + \left\|e_x\right\|_{H^{1/2}(\Gamma_h[{\bf x}^*(t)])} + \left\|e_v\right\|_{H^{1/2}(\Gamma_h[{\bf x}^*(t)])} \right. \\
        & \qquad + \left.\left\|e_H\right\|_{H^{1/2}(\Gamma_h[{\bf x}^*(t)])} + \left\|e_n\right\|_{H^{1/2}(\Gamma_h[{\bf x}^*(t)])} \right)\\
        & \leq C h^{\kappa - 3/2} \leq h^{(\kappa - 2)/2 + 1/2},
        \end{aligned}
    \end{equation*}
    for sufficiently small $h$. Therefore, the assumptions in \eqref{eq:error bounds, assumptions} hold and we can extend $t^*$ to $T$. We know that the stability bound \eqref{eq:stability bound} holds for all $t \in [0, T]$.
\end{proof}

\subsection{Completion of the proof}
With the established consistency error and stability bounds, we can now prove the main theorem.
\begin{proof}
    We first consider the error in the flow map $X$. According to the stability bound \eqref{eq:stability bound} and the trace inequality, we know that
    \begin{equation*}
        \left\| e_x \right\|_{H^1(\Omega_h[{\bf x}^*])} \leq C h^{k}.
    \end{equation*}
    From the norm equivalence along the flow map \cite[Lemma4.2]{KLLP2017} and the fact that $\hat{X}_h - X_h^*  = e_x \circ X_h^*$, we obtain
    \begin{equation*}
        \left\| \hat{X}_h - X_h^* \right\|_{H^1(\Omega_h^0)} \leq C \left\| e_x \right\|_{H^1(\Omega_h[{\bf x}^*])} \leq C h^{k}.
    \end{equation*}
    Since the norms are equivalent under lifting, we have
    \begin{equation*}
        \left\| X^L_h - (X_h^*)^\ell \right\|_{H^1(\Omega^0)} = \left\| \hat{X}_h^\ell - (X_h^*)^\ell \right\|_{H^1(\Omega^0)} \leq C \left\| \hat{X}_h - X_h^* \right\|_{H^1(\Omega_h^0)}
    \end{equation*}
    Furthermore, from interpolation error estimates, we have
    \begin{equation*}
        \left\| X - (X_h^*)^\ell \right\|_{H^1(\Omega^0)} \leq C h^{k}
    \end{equation*}
    
    Finally, by the triangle inequality, we obtain
    \begin{equation*}
        \left\| X^L_h - X \right\|_{H^1(\Omega^0)} \leq \left\| X^L_h - (X_h^*)^\ell \right\|_{H^1(\Omega^0)} + \left\| (X_h^*)^\ell - X \right\|_{H^1(\Omega^0)} \leq C h^{k}
    \end{equation*}
    Thus, we have established the error estimate for the flow map $X$. 
    
    For $u$, $n$, $H$, $x$, and $v$, the error estimates are analogous, thereby completing the proof of Theorem \ref{thm:ESFEM_error_bound}.
\end{proof}

\begin{remark}
    It is noteworthy that the $k\geq 3$ is required in three places: (1) the $H_h^s$ norm equivalence for $s = 3/2$ in Lemma \ref{lemma:norm equivalence for Hhs on different surfaces}; (2) the inverse estimate for $L^2$-$L^2$-$L^\infty$ estimate as in \eqref{eq: standard 22inf estimate} and Lemma \ref{lemma:comparison of interpolated surface and exact surface 2} requires an additional $h^{1/2}$, but the inverse estimate with convergence rate only gives $h^{k - 3/2}$; (3) \eqref{eq: tmp, third reason for k geq 3} and the Lemma \ref{lemma:time derivative of discrete Laplacian of order 1} requires $k-1 = \alpha_0 > 1$.
\end{remark}
\section*{Acknowledgements}
The work of Yifei Li is funded by the Alexander von Humboldt Foundation. I would like to thank Professor Christian Lubich for his encouragement, insightful suggestions, and valuable discussions. I am also grateful to Professor Buyang Li and Professor Bal\'{a}zs Kov\'{a}cs for their engaging discussions and for providing suitable references.

%\begin{acknowledgements}
%\end{acknowledgements}


%\bibliography{bulksurface_literature}   % name your BibTeX data base

\begin{thebibliography}{10}
    \bibitem{arioli2009discrete} M. Arioli and D. Loghin, \textit{Discrete interpolation norms with applications}, SIAM J. Numer. Anal. \textbf{47}(4) (2009), 2924--2951.
    
    

    \bibitem{barrett2020parametric} J. W. Barrett, H. Garcke, and R. N\"{u}rnberg, \textit{Parametric Finite Element Approximations of Curvature-Driven Interface Evolutions}, Handbook of numerical analysis, volume 21, Elsevier (2020), 275--423.

    \bibitem{bhatia2013matrix} R. Bhatia, \textit{Matrix Analysis}, volume 169, Springer Science \& Business Media (2013).

    \bibitem{brenner2008mathematical} S. C. Brenner and L. R. Scott, \textit{The Mathematical Theory of Finite Element Methods}, Springer (2008).

    %\bibitem{deckelnick2000error-levelset} K. Deckelnick, \textit{Error bounds for a difference scheme approximating viscosity solutions of mean curvature flow}, Interfaces Free Bound. \textbf{2}(2) (2000), 117--142.

    \bibitem{DeckelnickDE2005} K. Deckelnick, G. Dziuk, and C. M. Elliott, \textit{Computation of geometric partial differential equations and mean curvature flow}, Acta Numerica \textbf{14} (2005), 139--232.

    \bibitem{Dziuk88} G. Dziuk, \textit{Finite elements for the Beltrami operator on arbitrary surfaces}, Partial differential equations and calculus of variations, Lecture Notes in Math., 1357, Springer, Berlin (1988), 142--155.

    \bibitem{DziukElliott_ESFEM} G. Dziuk and C. Elliott, \textit{Finite elements on evolving surfaces}, IMA J. Numer. Anal. \textbf{27}(2) (2007), 262--292.

    \bibitem{DziukElliott_SFEM} G. Dziuk and C. Elliott, \textit{Surface finite elements for parabolic equations}, J. Comput. Math. \textbf{25}(4) (2007), 385--407.

    %\bibitem{DziukElliott_L2} G. Dziuk and C. Elliott, \textit{$L^2$--estimates for the evolving surface finite element method}, Math. Comp. \textbf{82}(281) (2013), 1--24.

    \bibitem{Edelmann_harmonicvelo} D. Edelmann, \textit{Finite element analysis for a diffusion equation on a harmonically evolving domain}, IMA J. Numer. Anal. \textbf{42} (2022), 1866--1901.

    \bibitem{EKL24} D. Edelmann, B. Kov\'{a}cs, and C. Lubich, \textit{Numerical analysis of an evolving bulk–surface model of tumour growth}, IMA Journal of Numerical Analysis (2024), drae077.

    \bibitem{elliott2024sfem} C. M. Elliott and A. Mavrakis, \textit{SFEM for the Lagrangian formulation of the surface Stokes problem}, arXiv preprint arXiv:2410.19470 (2024).

    \bibitem{Elliott2010} C. Elliott, B. Stinner, V. Styles, and R. Welford, \textit{Numerical computation of advection and diffusion on evolving diffuse interfaces}, IMA J. Numer. Anal. \textbf{31}(3) (2010), 786--812.

    \bibitem{ElliottRanner_bulksurface} C. M. Elliott and T. Ranner, \textit{Finite element analysis for a coupled bulk-surface partial differential equation}, IMA J. Numer. Anal. \textbf{33}(2) (2013), 377--402.

    \bibitem{ElliottRanner_unified} C. M. Elliott and T. Ranner, \textit{A unified theory for continuous-in-time evolving finite element space approximations to partial differential equations in evolving domains}, IMA J. Numer. Anal. \textbf{41}(3) (2021), 1696--1845.

    \bibitem{EKS19} J. Eyles, J. R. King, and V. Styles, \textit{A tractable mathematical model for tissue growth}, Interfaces Free Bound. \textbf{21}(4) (2019), 463--493.

    \bibitem{hu2022optimal} H. Hu, B. Li, and J. Zou, \textit{Optimal convergence of the Newton iterative Crank--Nicolson finite element method for the nonlinear Schr{\"o}dinger equation}, Comput. Methods Appl. Math. \textbf{22}(3) (2022), 591--612.

    \bibitem{Huisken1984} G. Huisken, \textit{Flow by mean curvature of convex surfaces into spheres}, J. Differential Geometry \textbf{20}(1) (1984), 237--266.

    

    \bibitem{highorderESFEM} B. Kov\'{a}cs, \textit{High-order evolving surface finite element method for parabolic problems on evolving surfaces}, IMA J. Numer. Anal. \textbf{38}(1) (2018), 430--459.

    \bibitem{KovacsPower_quasilinear} B. Kov\'{a}cs and C. Power Guerra, \textit{Error analysis for full discretizations of quasilinear parabolic problems on evolving surfaces}, Numer. Methods Partial Differential Equations \textbf{32}(4) (2016), 1200--1231.


    \bibitem{KLLP2017} B. Kov\'{a}cs, B. Li, C. Lubich, and C. Power Guerra, \textit{Convergence of finite elements on an evolving surface driven by diffusion on the surface}, Numer. Math. \textbf{137}(3) (2017), 643--689.

    \bibitem{MCF} B. Kov\'{a}cs, B. Li, and C. Lubich, \textit{A convergent evolving finite element algorithm for mean curvature flow of closed surfaces}, Numer. Math. \textbf{143} (2019), 797--853.

    \bibitem{MCF_soldriven} B. Kov\'{a}cs, B. Li, and C. Lubich, \textit{A convergent algorithm for forced mean curvature flow driven by diffusion on the surfaces}, Interfaces and Free Boundaries \textbf{22}(4) (2020), 443--464.

    \bibitem{kovacs2021convergent} B. Kov\'{a}cs, B. Li, and C. Lubich, \textit{A convergent evolving finite element algorithm for Willmore flow of closed surfaces}, Numer. Math. \textbf{149}(3) (2021), 595--643.


    %\bibitem{LubichMansour_wave} C. Lubich and D. Mansour, \textit{Variational discretization of wave equations on evolving surfaces}, Math. Comp. \textbf{84}(292) (2015), 513--542.

    \bibitem{Walker2015} S. W. Walker, \textit{The Shape of Things: A Practical Guide to Differential Geometry and the Shape Derivative}, SIAM, Philadelphia (2015).
\end{thebibliography}

    %\begin{thebibliography}{}
    % etc
    %\end{thebibliography}
    
\end{document}